\documentclass[11pt, a4paper]{amsart}
\usepackage{mathbbol}
\usepackage{amsmath, amsthm, amssymb}
\usepackage{txfonts, mathrsfs}
\usepackage[usenames]{color}
\usepackage{psfrag}
\usepackage{graphicx}

\newtheorem{thm}{Theorem}[section]
\newcommand{\bt}{\begin{thm}}
\newcommand{\et}{\end{thm}}

\newtheorem{cor}[thm]{Corollary}
\newcommand{\bc}{\begin{cor}}
\newcommand{\ec}{\end{cor}}

\newtheorem{lem}[thm]{Lemma}
\newcommand{\bl}{\begin{lem}}
\newcommand{\el}{\end{lem}}

\newtheorem{prop}[thm]{Proposition}
\newcommand{\bp}{\begin{prop}}
\newcommand{\ep}{\end{prop}}

\newtheorem{defn}[thm]{Definition}
\newcommand{\bd}{\begin{defn}}      
\newcommand{\ed}{\end{defn}}

\newtheorem{rmrk}[thm]{Remark}
\newcommand{\br}{\begin{rmrk}}
\newcommand{\er}{\end{rmrk}}

\newtheorem{quest}[thm]{Question}
\newcommand{\bq}{\begin{quest}}
\newcommand{\eq}{\end{quest}}

\newtheorem{example}[thm]{Example}

\newcommand{\C}{\mathbb{C}}
\newcommand{\N}{\mathbb{N}}

\newcommand{\R}{\mathbb{R}}

\newdimen\vintkern\vintkern12pt
\def\vint{-\kern-\vintkern\int}

\newcommand{\hm}{{\mathcal H}}
\newcommand{\lm}{{\mathcal L}}

\newcommand{\dist}{\operatorname{dist}}
\newcommand{\diam}{\operatorname{diam}}
\newcommand{\trace}{\operatorname{tr}}
\newcommand{\length}{\ell}
\newcommand{\Area}{\operatorname{Area}}

\newcommand{\md}{\operatorname{md}}
\newcommand{\lip}{\operatorname{Lip}}
\newcommand{\bdry}{\partial}

\newcommand{\Vol}{\operatorname{Vol}}

\newcommand{\jac}{{\mathbf J}}
\newcommand{\ap}{\operatorname{ap}}
\newcommand{\apmd}{\ap\md}

\begin{document}
\bibliographystyle{amsalpha}

\title[Plateau's problem in metric spaces]{Area minimizing discs in metric spaces}

\author{Alexander Lytchak}


\address
  {Mathematisches Institut\\ Universit\"at K\"oln\\ Weyertal 86 -- 90\\ 50931 K\"oln, Germany}
\email{alytchak@math.uni-koeln.de}

\author{Stefan Wenger}

\address
  {Department of Mathematics\\ University of Fribourg\\ Chemin du Mus\'ee 23\\ 1700 Fribourg, Switzerland}
\email{stefan.wenger@unifr.ch}

\date{\today}

\thanks{S.~W.~was partially supported by Swiss National Science Foundation Grant 153599.}

\begin{abstract}
 We solve the classical problem of Plateau in the setting of proper metric spaces. Precisely, we prove that among all disc-type surfaces with
 prescribed Jordan boundary in a proper metric space there exists an area minimizing disc which moreover has a quasi-conformal parametrization. If the space supports a local quadratic isoperimetric inequality for curves we prove that such a solution is locally H\"older continuous in the interior and continuous up to the boundary. Our results generalize corresponding results of Douglas and  Morrey from the setting of Euclidean space and Riemannian manifolds to that of proper metric spaces.
 \end{abstract}

\maketitle

\bigskip

\section{Introduction and statement of main results}\label{sec:Intro}

\subsection{Introduction} 
The classical problem of Plateau asks to prove the existence of a minimal disc bounded by a given Jordan curve in Euclidean space. The first rigorous solutions of Plateau's problem for arbitrary Jordan curves were given independently by Douglas \cite{Dou31} and Rad\'o \cite{Rad30}.  In a major advance, Morrey \cite{Mor48} extended the solutions of Douglas and Rado to a large class of Riemannian manifolds. 
Beyond the setting considered in \cite{Mor48}, the existence and regularity of area minimizing discs is only known in a few classes of metric spaces. In \cite{Nik79}  Nikolaev considered the case of metric spaces of curvature bounded from above in the sense of Alexandrov. In \cite{MZ10} Mese-Zulkowski treated the case of some spaces of curvature bounded from below in the sense of Alexandrov. Finally, Overath and von der Mosel \cite{OvdM14} treated the case of $\R^3$ endowed with a Finsler metric. The purpose of the present paper is to generalize these results to metric spaces under minimal additional conditions. 

Before describing our results, we briefly mention that there are many other ways to pose and sometimes to solve a Plateau type problem, see e.g.~\cite{Dav14} for some of these ways. For instance, one may minimize area among surfaces of fixed topological type or among integral currents (generalized surfaces of arbitrary topological type). The theory of integral currents, developed by Federer-Fleming \cite{FF60} in the setting of Euclidean spaces, has been generalized by Ambrosio-Kirchheim \cite{AK00} to the setting of arbitrary complete metric spaces. Their theory allows to prove existence of mass minimizing integral currents in compact metric spaces and some locally non-compact ones, see \cite{AK00}, \cite{Wen05}, \cite{AS13}, \cite{Wen14}. In contrast to the well developed regularity theory for mass minimizing integral currents in Euclidean space, see~\cite{Alm00}, the regularity of minimal currents in metric spaces seems to be very difficult to approach, see \cite{ADLS15} for some progress. We will not further discuss or pursue these directions here and refer the reader to the articles above and the reference mentioned therein.

We return to the main subject of the present paper which concerns existence and regularity of area minimizing discs in the setting of metric spaces. Before describing our results in more detail in Section~\ref{sec:prec-stat-intro} we give a rough description of some of the highlights of our paper. The natural analog of smooth discs in metric spaces are Lipschitz discs. Since Lipschitz maps lack suitable compactness properties needed for proving the existence of area minimizers it is inevitable to increase the range of admissible discs. As in the classical setting, a natural class to work with is that of Sobolev maps. Various equivalent definitions of Sobolev maps from a Euclidean domain with values in a metric space exist. Their parametrized area or volume can be defined in analogy with the parametrized area of a Lipschitz map via integration of a suitable Jacobian. In Riemannian manifolds this yields the parametrized $2$-dimensional Hausdorff measure. In the realm of normed spaces, there exist several natural definitions of area coming from convex geometry. This yields different notions of parametrized areas of Lipschitz or Sobolev maps with values in metric spaces, one of which is the parametrized $2$-dimensional Hausdorff measure. Our results apply to many of these notions of parametrized area. For the sake of simplicity we will first formulate our results for the one coming from the $2$-dimensional Hausdorff measure. In our first main result we show that the classical Plateau problem has a solution in any proper metric space $X$. That is, among all Sobolev discs (maps from the disc to $X$) spanning a given Jordan curve in $X$ there exists one of minimal area. Moreover, this map can be chosen to be $\sqrt{2}$-quasi-conformal. This means, roughly speaking, that infinitesimal balls are mapped to ellipses of aspect ratio at most $\sqrt{2}$. Simple examples show that the constant $\sqrt{2}$ is optimal. For a large class of metric spaces, however, we can improve the constant and obtain a (weakly) conformal map.
Similarly to the classical solution in Euclidean space, energy minimizers play an important role in our approach. We show that these are always $\sqrt{2}$-quasi-conformal, which is again optimal. In the setting of metric spaces, however, energy minimizers need not be area minimizers anymore as we will show and thus the classical approach to solving Plateau's problem fails; see however
Section~\ref{sec:ET-case} and \cite{LW-energy-area}. We circumvent this by proving a general lower semi-continuity result for the area which also yields new proofs of the lower semi-continuity of various energies. In the second part of the paper we prove interior and boundary regularity of quasi-conformal area minimizers in any metric space admitting a local quadratic isoperimetric inequality for curves. More precisely we prove that a quasi-conformal area minimizer is continuous up to the boundary and locally H\"older continuous in the interior with H\"older exponent only depending on the isoperimetric and the quasi-conformality constants. Our exponent is in many cases optimal. 

We now pass to a precise description of the results mentioned above and to further results and applications.

\subsection{Precise statements of main results}\label{sec:prec-stat-intro}
We now give a more detailed description of some of the main results in our paper. Recall that there exist several equivalent definitions of Sobolev maps from Euclidean domains with values in a metric space, see e.g.~\cite{Amb90}, \cite{KS93}, \cite{Res97}, \cite{Res04}, \cite{Res06}, \cite{HKST01}, \cite{HKST15}, \cite{AT04}. We recall the definition of \cite{Res97} using compositions with real-valued Lip\-schitz functions.

Let $X=(X,d)$ be a complete metric space. In this introduction we will restrict ourselves to maps defined on the open unit disc $D$ in $\R^2$ with values in $X$. For $p>1$ the Sobolev space $W^{1,p}(D, X)$ may be defined as the space of measurable and essentially separably valued maps $u\colon D\to X$ for which there exists a non-negative function $h\in L^p(D)$ with the following property: for every $x\in X$ the function $u_x(z):= d(x, u(z))$ belongs to the classical Sobolev space $W^{1,p}(D)$ and its weak gradient satisfies $|\nabla u_x|\leq h$ almost everywhere in $D$. 
%
%
Sobolev maps with values in $X$ are almost everywhere approximately metrically differentiable, that is, at almost every point $z\in D$ there exists a unique seminorm on $\R^2$, denoted $\apmd u_z$, such that 
 \begin{equation*}
    \ap\lim_{z'\to z}\frac{d(u(z'), u(z)) - \apmd u_z(z'-z)}{|z'-z|} = 0,
 \end{equation*}
 see Proposition~\ref{prop:Sobolev-apmd-Lip} below or \cite{Kar07}. Using the approximate metric differentiability one obtains a natural notion of quasi-conformality and parametrized area of Sobolev maps. We say that a seminorm $s$ on $\R^2$ is $Q$-quasi-conformal if $s(v)\leq Q\cdot s(w)$ for all $v,w\in S^{1}$. Note that $s\equiv 0$ is allowed. A map $u\in W^{1,p}(D, X)$ is called $Q$-quasi-conformal if its approximate metric derivative $\apmd u_z$ is $Q$-quasi-conformal at almost every $z\in D$. If $Q=1$ then we call $u$ conformal. We emphasize that our notion of quasi-conformal map is different from the notion of quasi-conformal homeomorphism studied in the field of quasi-conformal mappings. In fact, our spaces $X$ in general have arbitrary dimension and topology and thus quasi-conformal maps in our sense will rarely be (local) homeomorphisms.
 %

As mentioned above, there are several natural notions of parametrized area in metric spaces. We will first introduce the one induced by the Hausdorff $2$-measure and state our results in this case before discussing to which extent they apply to other notions.
The parametrized Hausdorff area of a Sobolev map $u\in W^{1,2}(D, X)$ is defined by $$\Area(u):= \int_D \jac_2(\apmd u_z)\,d\lm^2(z),$$ where the Jacobian $\jac_2(s)$ of a seminorm  $s$ is given by the Hausdorff measure (with respect to the distance $s$ on $\R^2$) of the Euclidean unit square. In view of \cite{Iva08} and the area formula for Lipschitz maps \cite{AKrect00} this gives a natural definition of parametrized Hausdorff area. If $u$ is an injective Lipschitz map or, more generally, an injective Sobolev map satisfying Lusin's property (N) then $\Area(u)$ is simply the $2$-dimensional Hausdorff measure of the image of $u$.

Given a Jordan curve $\Gamma$ in $X$ we denote by $\Lambda(\Gamma, X)$ the family of Sobolev maps $u\in W^{1,2}(D, X)$ whose trace has a representative which is a weakly monotone parametrization of $\Gamma$. A special case of our main theorem concerning a solution of the problem of Plateau in metric spaces can be stated as follows.

\bt\label{thm:intro-exist-area-min}
 Let $X$ be a proper metric space and $\Gamma\subset X$ a Jordan curve such that $\Lambda(\Gamma, X)\not=\emptyset$. Then there exists $u\in \Lambda(\Gamma, X)$ which satisfies 
 \begin{equation}\label{eq:area-min-intro}
  \Area(u) = \inf\left\{\Area(u'): u'\in \Lambda(\Gamma, X)\right\}
 \end{equation}
  and which is $\sqrt{2}$-quasi-conformal.
\et

Here, a metric space is said to be proper if every closed ball of finite radius is compact.
In general, the quasi-conformality constant $\sqrt{2}$ in our theorem is optimal, see Remark~\ref{rem:qc-optimal-constant}. However, it can be improved to conformality for a large class of geometrically interesting spaces, see the paragraph below. 

In the classical proof of the solution of Plateau's problem in Euclidean space one first minimizes the Dirichlet energy in the class $\Lambda(\Gamma, \R^n)$ and then shows that an energy minimizer also minimizes area. The same reasoning cannot be used in the generality we work in. Indeed, we will prove that there exist metric spaces biLipschitz homeomorphic to the standard two-dimensional sphere in which an energy minimizer is not an area minimizer, see Proposition~\ref{prop:area-min-diff} and the remark following it. Nevertheless, energy minimizers still play an important role in our proof. Before explaining their role, let us recall that Korevaar-Schoen \cite{KS93} and Reshetnyak \cite{Res97} introduced different energies of a Sobolev map $u\in W^{1,p}(D, X)$. Using the approximate metric derivative Reshetnyak's energy, which we denote by $E_+^p(u)$, and Korevaar-Schoen's energy, which we denote by $E^p(u)$, take the form 
\begin{equation}\label{eq:intro-KS-Res-energies}
 E_+^p(u) = \int_D \mathcal{I}_+^p(\apmd u_z)\,d\lm^2(z)\quad \text{and} \quad E^p(u)= \int_D \mathcal{I}^p_{\rm avg}(\apmd u_z)\,d\lm^2(z),
 \end{equation}
where for a seminorm $s$ on $\R^2$ we define 
 \begin{equation}\label{eq:intro-def-energy-integrands}
 \mathcal{I}^p_+(s):= \max\left\{s(v)^p: v\in S^{1}\right\}\quad\text{and}\quad \mathcal{I}^p_{\rm avg}(s):= \pi^{-1}\int_{S^{1}}s(v)^p\,d\hm^{1}(v),
\end{equation}
see Proposition~\ref{prop:rep-energy} and the paragraph preceding Proposition~\ref{prop:ac-1-dim}. Reshetnyak's energy $E_+^p(u)$ is equal to the $p$-th power of the $L^p$-norm of the minimal weak upper gradient of $u$ in the sense of \cite{HKST15}. If $X=\R^N$ and $p=2$ then Korevaar-Schoen's energy corresponds to the classical Dirichlet energy.
One of the main ingredients in the proof of Theorem~\ref{thm:intro-exist-area-min} is the following result, which is of independent interest.

\bt\label{thm:intro-qc-energy-min}
 Let $X$ be a complete metric space. Suppose that $u\in W^{1,2}(D, X)$ is such that 
 \begin{equation}\label{eq-intro-energy-min-inner}
 E^2_+(u)\leq E^2_+(u\circ\psi)
\end{equation}
 for every biLipschitz homeomorphism $\psi\colon D\to D$. Then $u$ is $\sqrt{2}$-quasi-conformal.
\et

As is the case for Theorem~\ref{thm:intro-exist-area-min}, the quasi-conformality constant $\sqrt{2}$ is optimal but can be improved to conformality for a large class of geometrically interesting spaces. Reshetnyak's energy $E^2_+$ can be replaced by Korevaar-Schoen energy $E^2$, however, we only obtain the quasi-conformality constant $Q=2\sqrt{2} + \sqrt{6}$ in this case, see Theorem~\ref{thm:qc-domain-minimizers-KS-energy}, which is probably not optimal.

We turn to the question of regularity of area minimizing discs in metric spaces. Without any further assumptions on the underlying space $X$ one cannot expect an area minimizer even to be continuous, not even in the setting of Riemannian manifolds, see \cite{Mor48}. We will prove interior and boundary regularity under the condition of a local quadratic isoperimetric inequality.

\bd\label{def:isop-intro}
 A complete metric space $X$ is said to admit a uniformly local quadratic isoperimetric inequality if there exist $l_0, C>0$ such that for every Lipschitz curve $c\colon S^1\to X$  of length $\length_X(c)\leq l_0$ there exists $u\in W^{1,2}(D, X)$ with $$\Area(u) \leq C \length_X(c)^2$$ and such that $\trace(u)(t) = c(t)$ for almost every $t\in S^1$.  
\ed

Many interesting classes of spaces admit uniformly local quadratic isoperimetric inequalities. These include homogeneously regular Riemannian manifolds in the sense of \cite{Mor48}, compact Lipschitz manifolds and, in particular, all compact Finsler manifolds; moreover, complete ${\rm CAT}(\kappa)$ spaces for every $\kappa\in\R$, compact Alexandrov spaces, and all Banach spaces. Further examples include the Heisenberg groups $\mathbb{H}^n$ of topological dimension $2n+1$ for $n\geq 2$, endowed with a Carnot-Carath\'eodory distance. See Section~\ref{sec:higher-integrability-area-min} for more examples and for references.

A special case of our main result concerning interior and boundary regularity of area minimizing discs can be stated as follows.

\bt\label{thm:reg-intro}
  Let $X$ be a complete metric space admitting a uniformly local quadratic isoperimetric inequality with constant $C$. Let $\Gamma\subset X$ be a Jordan curve and suppose $u\in \Lambda(\Gamma, X)$ is $Q$-quasi-conformal and satisfies 
   \begin{equation*}
   \Area(u) = \inf\left\{\Area(v): v\in \Lambda(\Gamma, X)\right\}.
  \end{equation*}
 Then the following statements hold:
 \begin{enumerate}
  \item There exists $p>2$ such that $u\in W^{1,p}_{\rm loc}(D, X)$; in particular, $u$ has a continuous representative $\bar{u}$ which moreover satisfies Lusin's property (N).
  \item The representative $\bar{u}$ is locally $\alpha$-H\"older continuous with $\alpha = (4\pi Q^2C)^{-1}$ and extends continuously to all of $\overline{D}$.
  \item If $\Gamma$ is a chord-arc curve then $\bar{u}$ is H\"older continuous on all of $\overline{D}$.
  \end{enumerate}
\et

Note that we do not make any assumptions on local compactness on $X$. The H\"older exponent $\alpha$ is, in principle, allowed to be larger or equal to $1$, which can be used to prove that spaces with small isoperimetric constants have to be trees, see Corollary~\ref{cor:tree} below. Unlike in the classical setting, the H\"older exponent $\alpha$ in the above theorem is optimal, see Example~\ref{example:Hoelder-exp-opt}. We would like to mention the following refinement of statement (ii) of Theorem~\ref{thm:reg-intro} proved in Section~\ref{sec:higher-integrability-area-min}.   The classical proof of the Sobolev embedding theorems provides also in our case a very strong form of H\"older continuity. Namely, the upper bound on the distance between points in $u(D)$  leading to (ii) of Theorem~\ref{thm:reg-intro} is given by estimating  the length of the image of some  curve connecting the corresponding points in $D$.  This result provides, in particular, plenty of rectifiable curves in the image $u(D)$  and can be used to understand to some extent the \emph{intrinsic structure} of the minimal disc.  We refer to the continuation of the present paper in \cite{LW-intrinsic}, where this structure will be investigated in detail.

Our results can be improved in a large class of geometrically relevant spaces. We say that a metric space $X$ satisfies property (ET) if  for every $u\in W^{1,2}(D, X)$ the approximate metric derivative $\apmd u_z$ is induced by a possibly degenerate inner product at almost every $z\in D$. Examples of such spaces include Riemannian manifolds with continuous metric tensor, metric spaces of curvature bounded from above or below in the sense of Alexandrov, equiregular sub-Riemannian manifolds, and infinitesimally Hilbertian spaces with lower Ricci curvature bounds. We will show in Section~\ref{sec:ET-case} that under the additional assumption that $X$ satisfies property (ET) the maps in Theorems~\ref{thm:intro-exist-area-min} and \ref{thm:intro-qc-energy-min} may be taken to be conformal, that is, $1$-quasi-conformal. Moreover, in this case Theorem~\ref{thm:intro-qc-energy-min} also holds with Reshetnyak's energy $E^2_+$ replaced by Korevaar-Schoen's energy $E^2$. Finally, in such spaces energy minimizers are automatically area minimizers, see Theorem~\ref{thm:energy-min-is-area-min-ET}. In particular, Theorem~\ref{thm:intro-qc-energy-min} generalizes the classical result from Euclidean space to that of arbitrary complete metric spaces. Theorem~\ref{thm:intro-exist-area-min} in conjunction with Theorem~\ref{thm:reg-intro} generalizes Douglas' and Morrey's solutions of Plateau's problem from the setting of Euclidean space and homogeneously regular Riemannian manifolds to that of proper metric spaces admitting a uniformly local quadratic isoperimetric inequality.  It also generalizes the results \cite{Nik79}, \cite{MZ10}, \cite{OvdM14} mentioned at the beginning of our introduction. 

As a first application of the results described above we obtain a solution of the absolute Plateau problem described as follows. Let $\Gamma$ be a metric space homeomorphic to $S^1$ and of finite length, for example, a Jordan curve of finite length in some metric space. We want to minimize the area of Sobolev maps $u\in\Lambda(\Gamma, X)$ not only for a fixed metric space $X$ containing $\Gamma$ but over all such spaces. Precisely, set $$m(\Gamma):= \inf\{\Area(v): \text{$Y$ complete, $\iota\colon\Gamma\hookrightarrow Y$ isometric, $u\in\Lambda(\iota(\Gamma), Y)$}\}.$$ 
This value is closely related to Gromov's filling area in \cite{Gro83}. In our setting this infimum is indeed attained, due to the following solution of the absolute Plateau problem.

\bc\label{cor:Abs-Plateau}
 Let $\Gamma$ be a metric space homeomorphic to $S^1$ and of finite length. Then there exist  a compact metric space $X$, an isometric embedding $\iota\colon \Gamma\hookrightarrow X$, and a map $u\in \Lambda(\iota(\Gamma), X)$ such that $$\Area(u) = m(\Gamma).$$ Moreover, $u$ is $\sqrt{2}$-quasi-conformal and has a representative which is continuous on $\overline{D}$ and locally $\frac{1}{4}$-H\"older continuous on $D$.
\ec

We will discuss the exact relation with Gromov's filling area and the relations of solutions to the absolute Plateau problem with boundary minimal surfaces \cite{Iva08} in the sequel \cite{LW-intrinsic} of this paper. The corollary above can be reformulated by saying that area minimizing discs with prescribed  boundary exist in all $L^\infty$-spaces and, more generally, in every injective metric space, see Theorem~\ref{thm:Plateau-injective-metric}.

Another simple application of our results is the following:

\bc\label{cor:tree}
 Let $X$ be a proper, geodesic metric space admitting a global quadratic isoperimetric inequality with some constant $C$. If $C<\frac{1}{8\pi}$ then $X$ is a metric tree, that is, every geodesic triangle in $X$ is isometric to a tripod.
\ec

If $X$ satisfies property (ET) then the corollary holds with $\frac{1}{8\pi}$ replaced by $\frac{1}{4\pi}$. In view of the Euclidean plane the constant $\frac{1}{4\pi}$ is optimal. The corollary is not new. In fact, it follows from \cite{Wen07} and \cite{Wen08-sharp} that the corollary holds for any geodesic metric space with the sharp constant $\frac{1}{4\pi}$. The borderline case $C=\frac{1}{4\pi}$ characterizes proper ${\rm CAT}(0)$-spaces, as will be shown in \cite{LW-isoperimetric}.

We finally discuss to what extent our results hold when the parametrized Hausdorff area is replaced by the para\-metrized area induced by other notions of volume. Roughly speaking, a definition of volume in the sense of convex geometry assigns to each $2$-dimensional normed space a constant multiple of the Haar measure in a consistent way, see Section~\ref{sec:def-vol-normed}. For example, the Busemann definition of volume assigns the multiple of the Haar measure for which the unit ball in a $2$-dimensional normed space has measure $\pi$, thus giving rise to the $2$-dimensional Hausdorff measure. Other widely used definitions of volume are Benson's or Gromov's mass$^*$ definition of volume, the symplectic or Holmes-Thompson definition of volume, or Loewner's definition of volume studied in \cite{Iva08}. Given a definition of volume $\mu$ one obtains a Jacobian with respect to $\mu$ of a seminorm on $\R^2$ in the same way as above except that one replaces the Hausdorff measure with the volume $\mu$ in the definition of jacobian. By integrating this yields the $\mu$-area of a Sobolev map $u$ denoted $\Area_\mu(u)$. 
A definition of volume is said to induce quasi-convex $2$-volume densities if any affine disc in a finite dimensional normed space has minimal $\mu$-area among all smooth discs with the same boundary. All the examples of definitions of volume mentioned above induce quasi-convex $2$-volume densities. We will show that Theorem~\ref{thm:intro-exist-area-min} remains true when the parametrized Hausdorff area is replaced by the $\mu$-area for any definition of volume $\mu$ which induces quasi-convex $2$-volume densities, see Theorem~\ref{thm:existence-qc-area-min-general-mu}. The same is true for Corollaries~\ref{cor:Abs-Plateau} and \ref{cor:tree}, see Theorem~\ref{thm:tree-gen} and Corollary~\ref{cor:Abs-Plateau-general}. Moreover, Theorem~\ref{thm:reg-intro} holds for $\mu$-area minimizers for any definition of volume $\mu$. 
We would like to mention that different choices of areas give rise to different minimizers, unless the space $X$ has the  property (ET) mentioned above, see Section~\ref{sec:ET-case}.  For some of these choices of area the minimizer can be found by minimization of an appropriate energy (for instance $E^+$ or $E^2$). This will be discussed in \cite{LW-energy-area}.

\subsection{Outlines of proofs}
We next provide outlines of proofs for some of our theorems stated above. 

The two main ingredients in the proof of Theorem~\ref{thm:intro-exist-area-min} are Theorem~\ref{thm:intro-qc-energy-min} above and Theorem~\ref{thm:weak-lsc-qc-fctl} below, which provides a generalization to metric spaces of the classical weak sequential lower semi-continuity result for quasi-convex integrands. Theorem~\ref{thm:weak-lsc-qc-fctl} yields the lower semi-continuity of the $\mu$-area functional along sequences of uniformly bounded energy whenever $\mu$ is a definition of volume which induces quasi-convex $2$-volume densities. Using Theorem~\ref{thm:intro-qc-energy-min} together with this lower semi-continuity we obtain that for every $u\in\Lambda(\Gamma, X)$ there exists $v\in\Lambda(\Gamma,X)$ which is $\sqrt{2}$-quasi-conformal and which satisfies $\Area(v)\leq \Area(u)$. We use this to pass from an area minimizing sequence to an area minimizing sequence of uniformly bounded energy and, together with compactness and the lower semi-continuity of the $\mu$-area, we obtain Theorem~\ref{thm:intro-exist-area-min}. We mention that Theorem~\ref{thm:weak-lsc-qc-fctl} is of independent interest. Indeed, we will use it to provide in Corollaries~\ref{cor:lsc-KS-energy} and \ref{cor:lsc-+-energy} new proofs of the weak lower semi-continuity of the Korevaar-Schoen energy in \cite[Theorem 1.6.1]{KS93} and the Reshetnyak energy \cite[Theorem 4.2]{Res97}.

We briefly describe the strategy of proof of Theorem~\ref{thm:intro-qc-energy-min}. Apart from establishing the result in the biggest possible generality our proof also seems more natural and transparent than the proof of the classical result that an energy minimizing Sobolev map with values in Euclidean space or a Riemannian manifold is (weakly) conformal. Recall that the proof of this classical result basically follows from a computation of the derivative of the function $t\mapsto E^2(u\circ\varphi_t)$ for a suitable family of diffeomorphisms $\varphi_t$ of $D$, see e.g.~\cite[Chapter 4.5]{Dierkes-et-al10}. The classical construction of the variation $\varphi _t$ is global and requires the precise knowledge of the derivative of the energy and depends on the values of $u$ on all of $D$, also far from the points with non-conformal derivative. The computation of the derivative of $E^2(u\circ\varphi_t)$ still works when $u$ is a metric space valued Sobolev map such that $\apmd u_z$ comes from a possibly degenerate inner product almost everywhere. In particular, the classical proof can be generalized to the setting of metric spaces satisfying property (ET) described above. However, this kind of proof breaks down when non-Euclidean normed tangent spaces appear somewhere. Our proof is local and works by finding, modulo conformal gauges, a local variation $\varphi _t$  in the neighborhood of a point where the quasi-conformality claim does not hold.

The proof of Theorem~\ref{thm:reg-intro} follows more or less literally the classical approach going back to Morrey.

\subsection{Structure and content of the paper}
We provide a description of the content of each section of the paper and indicate some of the more general versions of the results stated above.

 In Section~\ref{sec:prelim} we establish basic notation and recall some results which will be used in the sequel. In particular, we recall the notion of a definition of volume from convex geometry and mention the primary examples of which the Hausdorff measure is one. We recall the concept of quasi-convex $n$-volume densities induced by a definition of volume, which will play a role in the solution of Plateau's problem. We furthermore define the notion of quasi-conformal seminorms on $\R^n$ and introduce generalizations to $\R^n$ of the functionals $\mathcal{I}^p_+$ and $\mathcal{I}^p_{\rm avg}$ on seminorms on $\R^n$ mentioned in \eqref{eq:intro-def-energy-integrands} above.

In Section~\ref{sec:background-Sobolev-theory} we recall the necessary background from the theory of Sobolev maps from a domain in Euclidean space $\R^n$ into complete metric spaces. We will follow the approach of Korevaar-Schoen \cite{KS93} and recall the equivalence with Reshetnyak's approach \cite{Res97}. The reason for choosing Korevaar-Schoen's approach is that their theory already contains many of the ingredients which we will need. We will furthermore provide a proof of the equivalence of Korevaar-Schoen's theory with a metric space version of Hajlasz's theory \cite{Haj96} yielding Lipschitz continuity on suitable subsets, see Proposition~\ref{prop:equiv-Hajlasz-KS}. This will be used in the subsequent section in order to establish the approximate metric differentiability almost everywhere of Sobolev maps. 

We begin Section~\ref{sec:Lip-prop-Sobolev} by recalling the definition of the approximate metric derivative $\apmd u_z$ of a map $u$ from a Euclidean domain $\Omega\subset\R^n$ to a metric space $X$. We then prove the approximate metric differentiability almost everywhere of Sobolev maps $u\in W^{1,p}(\Omega, X)$ and a strong first order approximation for distances of image points via the approximate metric derivative, see Proposition~\ref{prop:Sobolev-apmd-Lip}. As a by-product we then obtain representations of the Korevaar-Schoen energy $E^p(u)$ and the Reshetnyak energy $E^p_+(u)$ in terms of the functionals $\mathcal{I}^p_{\rm avg}$ and $\mathcal{I}^p_+$ akin to \eqref{eq:intro-KS-Res-energies} and furthermore a new proof of the main result of Logaritsch-Spadaro \cite{LS12} on the representation of the Korevaar-Schoen energy. See Proposition~\ref{prop:rep-energy} and the paragraph preceding Proposition~\ref{prop:ac-1-dim}. We furthermore introduce the notion of parametrized $\mu$-volume $\Vol_\mu(u)$  of a Sobolev map $u\in W^{1,n}(\Omega, X)$ induced by a definition of volume $\mu$, see Definition~\ref{def:vol-sobolev}. In the particular case of dimension $n=2$ the $\mu$-volume will be denoted by $\Area_\mu(u)$. 

The main result of Section~\ref{sec:weak-lsc-integrands} is Theorem~\ref{thm:weak-lsc-qc-fctl}. It establishes a generalization to metric spaces of the classical weak sequential lower semi-continuity of quasi-convex integrands. Using Theorem~\ref{thm:weak-lsc-qc-fctl} we provide new proofs of the weak lower semi-continuity of the Korevaar-Schoen and the Reshetnyak energies in Corollaries~\ref{cor:lsc-KS-energy} and \ref{cor:lsc-+-energy}. As a further and direct consequence we obtain in Corollary~\ref{cor:quasi-convex-volume-lsc} the weak lower semi-continuity of the volume functional $\Vol_\mu(\cdot)$ for any definition of volume $\mu$ inducing quasi-convex $n$-volume densities. This is later used in order to prove the existence of area minimizers.

In Section~\ref{sec:qc-energy-min} we prove Theorem~\ref{thm:intro-qc-energy-min} above and its analog for the Korevaar-Schoen energy, Theorem~\ref{thm:qc-domain-minimizers-KS-energy}.

The main result in Section~\ref{sec:existence-area-min} is Theorem~\ref{thm:existence-qc-area-min-general-mu}, which provides  our most general existence theorem for area minimizers with prescribed Jordan boundary. Precisely, it states that if $X$ and $\Gamma$ are as in Theorem~\ref{thm:intro-exist-area-min} above and if $\mu$ is a definition of volume which induces quasi-convex $2$-volume densities then there exists $u\in\Lambda(\Gamma, X)$ which minimizes the $\mu$-area $\Area_\mu$ among maps in $\Lambda(\Gamma, X)$. Theorem~\ref{thm:existence-qc-area-min-general-mu} in particular implies Theorem~\ref{thm:intro-exist-area-min} above.

In Section~\ref{sec:higher-integrability-area-min} we state and prove our most general version of our results concerning interior regularity of area minimizers, Theorem~\ref{thm:int-reg-summary}. We furthermore provide many examples of spaces satisfying a uniform local quadratic isoperimetric inequality. 

In Section~\ref{sec:cont-boundary} we prove the boundary regularity results, Theorems~\ref{thm:bdry-cont-u-classical} and \ref{thm:bdry-Hoelder}, which in particular imply the boundary regularity results in Theorem~\ref{thm:reg-intro}.

Section~\ref{sec:abs-Plateau-corollaries} contains the proofs of Corollaries~\ref{cor:Abs-Plateau} and \ref{cor:tree}, in fact, the more general versions for $\mu$-area for definitions of volume inducing quasi-convex volume densities. We also solve Plateau's problem in every injective metric space, Theorem~\ref{thm:Plateau-injective-metric}.

In the final Section~\ref{sec:ET-case} we introduce the property (ET) mentioned above and show that many geometrically interesting classes of spaces have this property. We then prove that in spaces satisfying property (ET) energy minimizers are conformal and area minimizers, see Theorems~\ref{thm:inner-var-min-conformal-XEucl} and \ref{thm:energy-min-is-area-min-ET}. We also show in Proposition~\ref{prop:area-min-diff} that in spaces without property (ET) area minimizers with respect to two different definitions of area are in general different. In particular, this implies that in this setting energy minimizers need not be area minimizers.

\bigskip

{\bf Acknowledgments:} We would like to thank Luigi Ambrosio, Heiko von der Mosel, and Stephan Stadler for helpful comments and conversations.

\section{Preliminaries}\label{sec:prelim}

\subsection{Basic notation}\label{sec:basic-notation}

The following notation will be used throughout the paper. 
The Euclidean norm of a vector $v\in\R^n$ is denoted by $|v|$; the open unit disc in $\R^2$ is the set
$$D:= \{v\in\R^2: |v|<1\}.$$
Given open sets $U\subset V\subset\R^n$ we write $U\subset\subset V$ to mean that $\overline{U}\subset V$. Lebesgue measure on $\R^n$ is denoted by $\lm^n$. We denote by $\omega_n$ the Lebesgue measure of the unit ball in $\R^n$. The indicator function of a set $A\subset \R^n$ will be denoted by $1_A$. 
An open subset $\Omega\subset\R^n$ is called Lipschitz domain if the boundary of $\Omega$ can be locally written as the graph of a Lipschitz function defined on an open ball of $\R^{n-1}$. 

Let $X=(X,d)$ be a metric space. The open ball in $X$ of radius $r$ and center $x_0\in X$ is denoted by $$B_X(x_0,r) = \{x\in X: d(x_0, x)<r\}$$
or simply by $B(x_0,r)$ if there is no danger of ambiguity.
A Jordan curve in $X$ is a subset $\Gamma\subset X$ which is homeomorphic to $S^1$. Given a Jordan curve $\Gamma\subset X$, a continuous map $c\colon S^1\to X$ is called weakly monotone parametrization of $\Gamma$ if $c$ is the uniform limit of some homeomorphisms $c_i\colon S^1\to\Gamma$.
For $m\geq 0$ the $m$-dimensional Hausdorff measure on $X$ is denoted by $\hm^m_X$ or simply by $\hm^m$ if there is no danger of ambiguity. The normalizing constant is chosen in such a way that on Euclidean $\R^m$ the Hausdorff measure $\hm^m$ equals the Lebesgue measure.

\subsection{Rectifiable curves}
Let $X=(X,d)$ be a metric space. The length of a continuous curve $c\colon I\to X$, defined on an interval $I\subset\R$, is given by
\begin{equation*}
 \length_X(c):= \sup\left\{\sum_{i=1}^k d(c(t_i), c(t_{i+1})): k\in \N, t_i\in I, t_1<t_2<\dots<t_{k+1}\right\}.
\end{equation*}
We allow $I$ to be open, closed, or half-open. The definition naturally extends to continuous curves defined on $S^1$ and, more generally, on connected $1$-dimensional manifolds. A continuous curve of finite length is called rectifiable.

We will need the following elementary lemma which is akin to the lemma on the existence of a parametrization proportional to arc-length.

\bl\label{lem:reparam-Lip}
 Let $c\colon [a,b]\to X$ be a rectifiable curve. Then there exists a sense-preserving homeomorphism $\psi\colon [a,b]\to [a,b]$ such that $\psi$ and $c\circ\psi$ are Lipschitz.
\el

An analogous statement holds when $[a,b]$ is replaced by $S^1$.

\begin{proof}
 We may assume that $a=0$ and $b=1$ and that furthermore $l:= \length_X(c)>0$. Define $\varrho:[0,1]\to[0,1]$ by $$\varrho(t):= \frac{1}{2l} \cdot\length_X(c|_{[0,t]}) + \frac{t}{2}.$$ Then $\varrho$ is a homeomorphism and its inverse $\varrho^{-1}$ is $2$-Lipschitz. We set $\psi:= \varrho^{-1}$. It is straight-forward to check that $c\circ\psi$ is $2l$-Lipschitz.
\end{proof}

\subsection{Seminorms on $\R^n$}

For $n\geq 1$ let $\mathfrak{S}_n$ denote the set of all seminorms on $\R^n$. Endow $\mathfrak{S}_n$ with the metric $$d_{\mathfrak{S}_n}(s, s'):= \max\{|s(v) - s'(v)|: v\in S^{n-1}\}.$$ 
Then $(\mathfrak{S}_n, d_{\mathfrak{S}_n})$ is a proper metric space and may be viewed as a subset of $C^0(S^{n-1}, \R)$, where the latter is endowed with the supremum norm. For $p\geq 1$ we define continuous functions $\mathcal{I}^p_+: \mathfrak{S}_n\to[0,\infty)$ and $\mathcal{I}^p_{\rm avg}: \mathfrak{S}_n\to[0,\infty)$ by
\begin{equation*}
 \mathcal{I}^p_+(s):= \max\{s(v)^p: v\in S^{n-1}\}
\end{equation*}
and
\begin{equation*}
 \mathcal{I}^p_{\rm avg}(s):= \omega_n^{-1}\int_{S^{n-1}}s(v)^p\,d\hm^{n-1}(v),
\end{equation*}
where $\omega_n$ denotes the Lebesgue measure of the unit ball in $\R^n$. These functions will be used extensively later in the paper. We have the following easy fact.

\bl\label{lem:int-seminorm}
 If $s$ is a seminorm on $\R^n$ and $p\geq 1$ then
 \begin{equation*}
  n^{-1}\mathcal{I}^p_{\rm avg}(s)\leq \mathcal{I}^p_+(s) \leq \lambda \mathcal{I}^p_{\rm avg}(s),
 \end{equation*}
 where $\lambda>0$ is a constant depending only on $n$ and $p$. If $n=p=2$ then $\lambda$ can be chosen to be $1$.
\el

\begin{proof}
The first inequality follows from the fact that $n^{-1}\cdot \omega _n =\hm^{n-1} (S^{n-1})$. As for the second inequality, let $v_0\in S^{n-1}$ be such that $s(v_0)^p = \mathcal{I}^p_+(s)$. Define a seminorm $s'$ on $\R^n$ by $s'(rv_0 + w):= |r|$ for all $r\in\R$ and every $w\in\R^n$ orthogonal to $v_0$. It follows that $s(v) \geq s'(v)\cdot s(v_0)$ for every $v\in\R^n$ and hence
$$s(v_0)^p\cdot \mathcal{I}^p_{\rm avg}(s') \leq \mathcal{I}^p_{\rm avg}(s).$$ From this the second inequality follows with $\lambda:= \left(\mathcal{I}^p_{\rm avg}(s')\right)^{-1}$. Note here that $\lambda$ only depends on $n$ and $p$. If $n=p=2$ then $$\lambda^{-1}= \mathcal{I}^2_{\rm avg}(s')= \pi^{-1}\int_0^{2\pi}\cos(\theta)^2\,d\theta = 1.$$
This completes the proof.
\end{proof}

\bd\label{def:def-qc-seminorm}
 A seminorm $s$ on $\R^n$ is called $Q$-quasi-conformal, $Q\geq 1$, if $$s(v)\leq Q\cdot s(w)\quad\text{for all $v,w\in S^{n-1}$.}$$
\ed

If $s$ is $1$-quasi-conformal then $s$ will be called conformal. Note that according to our definition the trivial seminorm $s=0$ is conformal. A seminorm $s$ on $\R^n$ is conformal if and only if $n^{-1}\mathcal{I}^p_{\rm avg}(s)= \mathcal{I}^p_+(s)$.

\bl\label{lem:qc-transform-energy-seminorm}
 Let $s$ be a seminorm on $\R^n$ and $Q\geq 1$. Let $T:\R^n\to\R^n$ be linear and bijective. If the norm $v\mapsto |T(v)|$ is $Q$-quasi-conformal then 
 \begin{equation}\label{eq:+-ineqs-qc}
  Q^{-(n-1)} \cdot \mathcal{I}^n_+(s) \leq |\det T|^{-1}\cdot \mathcal{I}^n_+(s\circ T)\leq Q^{n-1} \cdot \mathcal{I}^n_+(s)
 \end{equation}
and
  \begin{equation}\label{eq:avg-ineqs-qc}
 Q^{-2(n-1)} \cdot \mathcal{I}^n_{\rm avg}(s) \leq |\det T|^{-1}\cdot \mathcal{I}^n_{\rm avg}(s\circ T) \leq Q^{2(n-1)} \cdot \mathcal{I}^n_{\rm avg}(s).
 \end{equation}
\el

\begin{proof}
 Writing $T$ as the product as $T=A\cdot D\cdot P$, where $D$ is a diagonal matrix and $A$ and $P$ are orthogonal transformations we obtain that 
 \begin{equation}\label{eq:T-trans-Qqc-energy}
  Q^{-(n-1)} \|T\|^n \leq |\det T|\leq Q^{n-1} t^n,
 \end{equation}
 where $\|T\|$ denotes the operator norm of $T$ and $t:= \min_{v\in S^{n-1}}|T(v)|$.
The inequalities in \eqref{eq:+-ineqs-qc} easily follow from this. In order to prove the inequalities in \eqref{eq:avg-ineqs-qc} we integrate using polar coordinates to obtain for every $R>0$ that
 \begin{equation*}
\int_{S^{n-1}} s(v)^n\,d\hm^{n-1}(v) =  2n R^{-2n} \int_{B(0,R)} s(w)^n\,d\lm^n(w).
\end{equation*}
From this we infer that
\begin{equation*}
 \int_{S^{n-1}} s(T(v))^n\,d\hm^{n-1}(v) = 2n |\det T|^{-1}\int_{T(B(0,1))} s(v)^n\,d\lm^{n}(v).
\end{equation*}
The inequalities in \eqref{eq:avg-ineqs-qc} follow from this, the inequalities \eqref{eq:T-trans-Qqc-energy}, and the fact that $T(B(0,1))\subset B(0,\|T\|)$.
\end{proof}

\subsection{Definitions of volume in normed spaces}\label{sec:def-vol-normed}
In Euclidean space there exists essentially only one natural definition of volume, which is the Lebesgue measure. In contrast, in the realm of normed spaces, there exist several natural notions of volume. Recall from \cite{AlvT04} the following definition.

\bd\label{def:volume-def}
 A definition of volume $\mu$ is a function that assigns to each $n$-dimensional normed space $V$, $n\geq 1$, a norm $\mu_V$ on $\Lambda^n V$  such that the following properties hold:
 \begin{enumerate}
\item If $V$ is Euclidean then $\mu_V$ is induced by the Lebesgue measure;
\item If $V$, $W$ are $n$-dimensional normed spaces and $T\colon V\to W$ linear and $1$-Lipschitz then the induced map $T_*\colon\Lambda^nV \to \Lambda^nW$ is $1$-Lipschitz; 
\end{enumerate}
\ed

Well-known examples of definitions of volume are the Busemann definition $\mu^{\rm b}$, the Holmes-Thompson definition $\mu^{\rm ht}$, and the Benson (also called Gromov mass$^*$) definition $\mu^{m^*}$ of volume, see e.g.~\cite{AlvT04}. We also mention the Loewner (or intrinsic Riemannian) volume $\mu^{\rm i}$ studied by Ivanov \cite{Iva08}.

Let $\mu$ be a definition of volume. Define the Jacobian with respect to $\mu$ of a seminorm $s$ on $\R^n$ by
\begin{equation*}
 \jac_n^\mu(s):= \left\{\begin{array}{ll}
  \mu_{(\R^n, s)}(e_1\wedge \dots\wedge e_n) & \text{if $s$ is a norm}\\
  0&\text{otherwise,}
 \end{array}\right.
\end{equation*}
where $e_1, \dots, e_n$ denote the standard unit vectors in $\R^n$. Note that the function $s\mapsto \jac^\mu_n(s)$ is continuous with respect to the metric $d_{\mathfrak{S}_n}$.

Let $\Omega\subset\R^n$ be an open, bounded subset and $Y$ a finite dimensional normed space or a Finsler manifold. Define the parametrized $\mu$-volume of a Lipschitz map $u\colon\Omega\to Y$ by 
\begin{equation*}
 \Vol_\mu(u):= \int_{\Omega} \jac_n^\mu(d_z u)\,d\lm^n(z). 
\end{equation*}
When $n=2$ we will write $\Area_\mu(u)$ instead of $\Vol_\mu(u)$. We will extend this definition to Sobolev maps from $\Omega$ to an arbitrary complete metric space in Definition~\ref{def:vol-sobolev}. The notion of parametrized volume of a Lipschitz map is a particular instance of the volume of a generalized Lipschitz surface in a metric space defined in \cite{Iva08}.

Recall the following definition.

\bd\label{def:vol-def-quasi-convex}
 Let $\mu$ be a definition of volume and $n\geq 1$. Then $\mu$ is said to induce quasi-convex $n$-volume densities if for every finite dimensional normed space $Y$ and every linear map $L: \R^n\to Y$ we have
 \begin{equation*}
  \Vol_\mu(L|_B)\leq \Vol_\mu(\psi)
 \end{equation*}
 for every smooth immersion $\psi: B\to Y$ with $\psi|_{\partial B} = L|_{\partial B}$, where $B$ denotes the closed unit ball in $\R^n$. 
\ed

Other names exist for this property in the literature. For example, in \cite{Iva08} the property is termed topologically semi-elliptic. Many known definitions of volume induce quasi-convex $n$-volume densities. Indeed, if a definition of volume induces extendibly convex $n$-volume densities (see e.g.~\cite{AlvT04} for the definition) in every finite dimensional normed space then it induces quasi-convex $n$-volume densities in the sense of Definition~\ref{def:vol-def-quasi-convex}. This follows directly from \cite[Theorem 4.23]{AlvT04}. By \cite[Theorem 4.28]{AlvT04}, the Gromov mass$^*$ definition of volume $\mu^{m^*}$ induces extendibly convex $n$-volume densities in every finite dimensional normed space for every $n\geq 1$. By \cite[Theorem 6.2]{Iva08}, the same is true for the intrinsic Riemannian volume definition $\mu^{\rm i}$. By \cite{BI13}, the Busemann definition of volume $\mu^{\rm b}$ induces extendibly convex $n$-volume densities in every finite dimensional normed space for $n=2$. A well-known conjecture asserts that this be true for all $n$. The volume densities of the Holmes-Thompson definition of volume $\mu^{\rm ht}$ are not extendibly convex, see \cite{BI02}. However, $\mu^{\rm ht}$ induces quasi-convex $2$-volume densities by \cite[Theorem 1, Section 3.1]{BI02}. In \cite{Ber14}, a new definition of volume was introduced which induces extendibly convex $n$-densities in every finite dimensional normed space for all $n$ and which coincides with the Busemann definition for $n=2$.

\section{Sobolev maps from Euclidean to metric spaces}\label{sec:background-Sobolev-theory}

We briefly recall Korevaar-Schoen's definition of Sobolev maps from Riemannian domains to metric spaces given in \cite{KS93}. Since we only need Euclidean domains we will restrict to this setting. In Section~\ref{sec:Lip-prop-Sobolev} we will establish several properties of Sobolev maps which will be useful in the rest of the paper.

Throughout this section, let $\Omega\subset \R^n$ be an open, bounded subset and $(X,d)$ a complete metric space. A map $u\colon\Omega\to X$ is measurable if for every open set $V\subset X$ the preimage $u^{-1}(V)$ is Lebesgue measurable. Furthermore, $u$ is essentially separably valued if there exists a set $N\subset\Omega$ of measure zero such that $u(\Omega\setminus N)$ is separable. 
For $p\geq 1$ denote by $L^p(\Omega, X)$ the space of all measurable and essentially separably valued maps $u:\Omega\to X$ such that for some and thus every $x_0\in X$ the function $z\mapsto d(x_0, u(z))$ belongs to $L^p(\Omega)$. 
 A sequence $(u_k)\subset L^p(\Omega, X)$ is said to converge to $u\in L^p(\Omega, X)$ in $L^p(\Omega, X)$ if $$\int_\Omega d^p(u(z), u_k(z))\,d\lm^n(z)\to 0$$ as $k\to\infty$.
Given $\varepsilon>0$ define
$$e_\varepsilon^p(z, u):=(n+p) \vint_{B(z,\varepsilon)} \frac{d^p(u(z), u(z'))}{\varepsilon^p}\, d\lm^n(z')$$ for all $z\in\Omega_\varepsilon:= \{ z'\in\Omega: \dist(z', \bdry \Omega)>\varepsilon\}$ and $e_\varepsilon^p(z, u):= 0$ for $z\in\Omega\backslash \Omega_\varepsilon$.
If $\varphi\in C_c(\Omega)$ then write 
$$E_\varepsilon^p(\varphi, u):= \int_{\Omega_\varepsilon} \varphi(z) e_\varepsilon^p(z,u) \, d\lm^n(z).$$
The Korevaar-Schoen $p$-energy of a map $u\in L^p(\Omega, X)$ is defined by
$$E^p(u):= \sup_{{\varphi\in C_c(\Omega),\, 0\leq \varphi\leq 1}} \limsup_{\varepsilon\to 0} E_\varepsilon^p(\varphi, u).$$ Note that $E^p(u)$ differs by a factor of $\omega_n^{-1}$ from the $p$-energy defined in \cite{KS93}, where $\omega_n$ is the volume of the unit ball in $\R^n$. 

For $p>1$ the Sobolev space $W^{1,p}(\Omega, X)$ in the sense of Korevaar-Schoen is the set of maps $u\in L^p(\Omega, X)$ satisfying $E^p(u)<\infty$. The space $W^{1,p}_{\rm loc}(\Omega, X)$ is defined analogously. If $u\in W^{1,p}(\Omega, X)$ and if $\varphi\colon X\to Y$ is a Lipschitz map into a complete metric space $Y$ then $\varphi\circ u\in W^{1,p}(\Omega, Y)$.

It was shown in \cite[Theorem 1.5.1]{KS93} that if $u\in W^{1,p}(\Omega, X)$ then the measures $e_\varepsilon^p(\,\cdot\,,u)d\lm^n$ converge weakly  as $\varepsilon\to0$ to an energy density measure $de^p(\,\cdot\,, u)$ with total measure $E^p(u)$. Moreover,  the measure $de^p(\,\cdot\,, u)$ is absolutely continuous with respect to the Lebesgue measure by \cite[Theorem 1.10]{KS93}. Finally, if $X=\R$ then $W^{1,p}(\Omega, X)$ coincides with the classical Sobolev space $W^{1,p}(\Omega)$ and the energy density of an element $u$ satisfies $$de^p(\,\cdot\,, u) = c_{n,p}|\nabla u(\cdot)|^p\,d\lm^n,$$ where $\nabla u$ is the weak derivative of $u$ and $c_{n,p}$ is a constant depending only on $n$ and $p$, see \cite[Theorem 1.6.2]{KS93}.

In Sections~\ref{sec:higher-integrability-area-min} and \ref{sec:cont-boundary} we will use the following terminology. Let $\Gamma\subset\R^n$ be a subset biLipschitz homeomorphic to an open interval $I$, and let $u\colon \Gamma\to X$ be a map. We write $u\in W^{1,p}(\Gamma, X)$ if $u\circ\varphi\in W^{1,p}(I, X)$ for some and thus any biLipschitz homeomorphism $\varphi\colon I\to \Gamma$. This terminology naturally extends to the case when $\Gamma$ is biLipschitz homeomorphic  to $S^1$.

As was shown in \cite{Res04}, the spaces $W^{1,p}(\Omega, X)$ can be characterized using compositions with Lipschitz functions on $X$. See \cite{Amb90} for an earlier approach towards metric space valued BV functions.

\bp\label{prop:Reshetnyak-Sobolev}
 Let $p>1$ and $u\in L^p(\Omega, X)$. Then $u\in W^{1,p}(\Omega, X)$ if and only if there exists $h\in L^p(\Omega)$ such that for every $x\in X$ the function $u_x(z):= d(x, u(z))$ belongs to $W^{1,p}(\Omega)$ and its weak gradient satisfies $|\nabla u_x|\leq h$ almost everywhere in $\Omega$.
\ep

Moreover, if $u$ and $h$ are as in Proposition~\ref{prop:Reshetnyak-Sobolev} then $E^p(u) \leq C \|h\|_p^p$ for some constant $C$ only depending on $n$ and $p$, see \cite{Res04}. In fact, we will see that one may even take $C=n$, see \eqref{eq:rel-E+-Eavg}. Finally, $h$ in Proposition~\ref{prop:Reshetnyak-Sobolev} can be chosen such that $\|h\|_p^p \leq \lambda E^p(u)$ for some constant $\lambda$ only depending on $n$ and $p$, see \cite{Res04} and \eqref{eq:rel-E+-Eavg}.

Apart from Proposition~\ref{prop:Reshetnyak-Sobolev} the following characterization of Sobolev maps will be important throughout our paper. 

\bp\label{prop:equiv-Hajlasz-KS}
 Let $p>1$ and $u\in L^p(\Omega, X)$. Then $u\in W^{1,p}(\Omega, X)$ if and only if there exist $g\in L^p(\Omega)$ and $N\subset\Omega$ with $\lm^n(N)=0$ such that 
 \begin{equation}\label{eq:Hajlasz-Sobolev}
  d(u(z), u(z')) \leq |z-z'| (g(z) + g(z'))
 \end{equation} for all $z,z'\in \Omega\setminus N$ contained in some ball $B\subset\subset\Omega$.  
Moreover, if $E^p(u)<\infty$ and $\Omega$ is a Lipschitz domain then $g$ may be be chosen so that \eqref{eq:Hajlasz-Sobolev} holds for all $z,z'\in \Omega\setminus N$.
\ep

A theory of Sobolev functions based on the condition \eqref{eq:Hajlasz-Sobolev} when $\Omega$ is replaced by a metric measure space and $X=\R$ was initiated in \cite{Haj96}. The proof of Proposition~\ref{prop:equiv-Hajlasz-KS} essentially follows from arguments in \cite{HKST01}, see also \cite{HKST15}. For the convenience of the reader we give a direct and self-contained proof here. 

\begin{proof}
 Suppose first that there exist $g\in L^p(\Omega)$ and $N\subset \Omega$ negligible such that \eqref{eq:Hajlasz-Sobolev} holds for all $z,z'\in\Omega\setminus N$ contained in some ball $B\subset\subset\Omega$. If $\varepsilon>0$ then 
 \begin{equation*}
  \frac{d^p(u(z), u(z'))}{\varepsilon^p} \leq (g(z) + g(z'))^p \leq 2^{p-1}(g^p(z) + g^p(z'))
 \end{equation*}
 for all $z,z'\in\Omega_\varepsilon\setminus N$ with $|z-z'|<\varepsilon$. In particular, we have $$e_\varepsilon^p(z,u)\leq 2^{p-1}(n+p) \left(g^p(z) + \vint_{B(z,\varepsilon)}g^p(z')\,d\lm^n(z')\right)$$ for every $z\in \Omega_\varepsilon\setminus N$. Therefore, given $\varphi\in C_c(\Omega)$ with $0\leq\varphi\leq 1$, we obtain 
 \begin{equation*}
 \begin{split}
  E_\varepsilon^p(\varphi, u)&\leq 2^{p-1}(n+p) \left(\int_\Omega g^p(z)\,d\lm^n(z) + \int_{\Omega_\varepsilon}\vint_{B(z,\varepsilon)}g^p(z')\,d\lm^n(z')\,d\lm^n(z)\right)\\
  & \leq 2^p(n+p)\|g\|_p^p
  \end{split}
  \end{equation*}
  for every $\varepsilon>0$ and thus $E^p(u)\leq 2^p(n+p)\|g\|_p^p<\infty$.
  
  Conversely, suppose $E^p(u)<\infty$. By Proposition~\ref{prop:Reshetnyak-Sobolev} there exists $h\in L^p(\Omega)$ such that for every $x\in X$ the function $u_x(z):= d(x, u(z))$ belongs to $W^{1,p}(\Omega)$ and its weak gradient satisfies $|\nabla u_x|\leq h$ almost everywhere in $\Omega$. Let $\{z_i\}_{i\in\N}\subset \Omega$ be a countable dense subset. For each $i$ let $B_i\subset\Omega$ be the open ball of maximal radius centered at $z_i$. Fix $x\in X$. There then exist negligible sets $N_i\subset B_i$ such that 
   \begin{equation}\label{eq:Lip-max-op-ball}
    |u_x(z) - u_x(z')| \leq C |z-z'| (M(|\nabla u_x|)(z) + M(|\nabla u_x|)(z'))
 \end{equation} 
for all $z,z'\in B_i\setminus N_i$, by e.g.\ \cite[Lemma 7.16]{GT01} and \cite[Lemma 2.83]{Zie89}. Here, $C$ is a constant depending only on $n$, and $M(|\nabla u_x|)$ denotes the maximal operator of $|\nabla u_x|$.
Set $N':=\cup N_i$ and note that $N'$ is negligible. 
Define $g(z):= C M(h)(z)$ for all $z\in \Omega$. Since $h\in L^p(\Omega)$ it follows from the maximal function theorem that $g$ is in $L^p(\Omega)$. Moreover, by \eqref{eq:Lip-max-op-ball}, we have
 \begin{equation}
    |u_x(z) - u_x(z')| \leq |z-z'| (g(z) + g(z'))
 \end{equation} 
 for all $z,z'\in \Omega\setminus N'$ such that $z,z'\in B_i$ for some $i$. 
 
   Since $u$ is essentially separably valued it readily follows from the above that there exists a negligible set $N\subset \Omega$ such that \eqref{eq:Hajlasz-Sobolev} holds for all $z,z'\in \Omega\setminus N$ contained in some $B_i$. Since every ball $B\subset\subset\Omega$ is contained in $B_i$ for some $i$ this proves the claim.

In order to prove the last statement of the proposition, suppose that $\Omega$ is a bounded Lipschitz domain and $E^p(u)<\infty$. We begin by making the following observation. There exist finitely many open subsets $U_i$ of $\Omega$, $i=1,\dots, m$, each of which is biLipschitz homeomorphic to a ball and such that for all $z,z'\in \Omega$ there exists $i$ such that $z$ and $z'$ are both in $U_i$. This observation is used as follows to prove the last statement. As explained above, for every $x\in X$ and every open ball $B\subset\Omega$ there exists $N_B\subset B$ negligible such that \eqref{eq:Lip-max-op} holds for all $z,z'\in B\setminus N_B$ and for a constant $C$ only depending on $n$. The same is then true with $B$ replaced by a biLipschitz copy of $B$ in $\Omega$ and with a constant $C$ depending on $n$ and the biLip\-schitz constant of the homeomorphism. From this together with the observation we obtain that for every $x\in X$ there exist $C$ (possibly depending on $\Omega$ but not on $x$) and a negligible set $N\subset \Omega$ such that 
  \begin{equation}\label{eq:Lip-max-op}
    |u_x(z) - u_x(z')| \leq C |z-z'| (M(|\nabla u_x|)(z) + M(|\nabla u_x|)(z'))
 \end{equation} 
 for all $z,z'\in \Omega\setminus N$. The same arguments as above show that \eqref{eq:Hajlasz-Sobolev} holds for all $z,z'\in \Omega\setminus N$ for some negligible set $N\subset\Omega$.
\end{proof}

\bp\label{prop:Hoelder-p>n}   
 Let $u\in W^{1,p}(\Omega, X)$ with $p>n$. Then $u$ has a unique representative $\bar{u}$ satisfying
 \begin{equation}\label{eq:Hold-cont-rep}
  d(\bar{u}(z), \bar{u}(z')) \leq C |z-z'|^{1-\frac{n}{p}}
 \end{equation}
 for every ball $B\subset \Omega$ and all $z,z'\in B$, where $C$ depends only on $n$, $p$, and $E^p(u)$. Moreover, $\bar{u}$ satisfies Lusin's property (N) and the set $\bar{u}(\Omega)$ is countably $\hm^n$-rectifiable. Finally, if $\Omega$ is a Lipschitz domain then \eqref{eq:Hold-cont-rep} holds for all $z,z'\in \Omega$ with a constant $C$ depending on $n$, $p$, $E^p(u)$, and $\Omega$. 
\ep

We recall that a map $\bar{u}\colon \Omega\to X$ is said to satisfy Lusin's property (N) if $\hm^n(\bar{u}(A)) = 0$ whenever $A\subset \Omega$ has measure $0$. Moreover, a set $A\subset X$ is called countably $\hm^n$-rectifiable if there exist countably many Lipschitz maps $\varphi_i\colon K_i\subset\R^n\to X$, $i\in\N$, such that $\hm^n(A\setminus \cup \varphi_i(K_i)) = 0$. 

\begin{proof}
 The first and last statement of the proposition are a consequence of Proposition~\ref{prop:Reshetnyak-Sobolev} and the remark following it together with Morrey's inequality for classical Sobolev functions. The fact that the continuous representative $\bar{u}$ of $u$ satisfies Lusin's property (N) then follows e.g. as in Proposition 2.4 of \cite{BMT13}.
Finally, the countable $\hm^n$-rectifiability of $\bar{u}(\Omega)$ is a consequence of Proposition~\ref{prop:equiv-Hajlasz-KS} together with the fact that $\bar{u}$ satisfies Lusin's property (N).
\end{proof}

Suppose now that $\Omega\subset\R^n$ is a bounded Lipschitz domain. The trace of a Sobolev map $u\in W^{1,p}(\Omega, X)$ with  $p>1$ can be defined as follows. Set $J=(-1,1)$ and $I=(-1,0)$. Given $x\in\bdry\Omega$ there exists an open neighborhood $U\subset\R^n$ of $x$, an open set $V\subset\R^{n-1}$, and a biLipschitz homeomorphism $\varphi\colon V\times J\to U$ such that $\varphi(V\times I) = U\cap\Omega$ and $\varphi(V\times\{0\}) = U\cap \bdry\Omega$. For $\lm^{n-1}$-almost every $v\in V$ the map $t\mapsto u\circ\varphi(v,t)$ is in $W^{1,p}(I, X)$ and thus has an absolutely continuous representative, again denoted by $u\circ\varphi(v,\cdot)$. For $\hm^{n-1}$-almost every point $z\in U\cap\bdry\Omega$ the trace of $u$ at $z$ is defined by $$\trace(u)(z):= \lim_{t\to0^{-}} u\circ\varphi(v, t),$$ where $v\in V$ is such that $\varphi(v,0) = z$. It follows from \cite[Lemma 1.12.1]{KS93} that the definition of $\trace(u)$ is independent of the choice of $\varphi$ and thus, by using a finite number of biLipschitz maps, is well-defined $\hm^{n-1}$-almost everywhere on $\bdry\Omega$. Furthermore, by \cite[Theorem 1.12.2]{KS93},  the trace map $\trace(u)$ is in  $L^p(\bdry\Omega, X)$. 

The following lemma will be used in the proof of Theorem~\ref{thm:existence-qc-area-min-general-mu}.

\bl\label{lem:bound-lp-dist-basepoint}
 Let $\Omega\subset\R^n$ be a bounded Lipschitz domain, and let $x_0\in X$ and $R>0$. If $u\in W^{1,p}(\Omega, X)$ with $1<p\leq n$ and such that $$\trace(u)(z)\in B(x_0, R)$$ for almost every $z\in\partial\Omega$ then
 \begin{equation}\label{eq:Ln-norm-center-bounded}
  \int_\Omega d^p(u(z), x_0)\,d\lm^n(z) \leq C\left(R^p + E^p(u)\right),
 \end{equation} 
 where $C$ is a constant only depending on $\Omega$ and $n$ and $p$.
\el

\begin{proof}
Define a $1$-Lipschitz function $\varphi: X\to\R$ by $$\varphi(x):= \max\{0, d(x,x_0) - R\}$$ and note that $\varphi(x) = 0$ for all $x\in B(x_0, R)$. Then $\varphi\circ u$ belongs to the classical Sobolev space $W^{1,p}(\Omega)$ by \cite[Theorem 1.6.2]{KS93}. In particular, $\varphi\circ u$ is approximately differentiable almost everywhere with approximate derivative equal to the weak derivative. It thus follows that at almost every $z\in\Omega$, the weak derivative of $\varphi\circ u$ is bounded by $$|d_z(\varphi\circ u)(v)| \leq \apmd u_z(v)$$ for every $v\in\R^n$. By Lemma~\ref{lem:int-seminorm}, there exists $\lambda$ depending only on $n$ and $p$ such that $$|\nabla(\varphi\circ u)(z)|^p\leq \mathcal{I}_+^p(\apmd u_z) \leq \lambda \mathcal{I}_{\rm avg}^p(\apmd u_z)$$
for almost every $z\in\Omega$, and hence $$\int_{\Omega} |\nabla(\varphi\circ u)(z)|^p\,d\lm^n(z) \leq \lambda E^p(u).$$
Since $\trace(\varphi\circ u) =0$ we may use the Sobolev inequality together with H\"older's inequality to estimate 
\begin{equation*}
 \|\varphi\circ u\|_p \leq C' \left(\int_{\Omega} |\nabla(\varphi\circ u)(z)|^p\,d\lm^n(z)\right)^{\frac{1}{p}} \leq C'\lambda^{\frac{1}{p}} E^p(u)^{\frac{1}{p}}
\end{equation*}
for a constant $C'$ only depending on $n$ and $p$ and $\Omega$. Now, inequality~\eqref{eq:Ln-norm-center-bounded} follows.
\end{proof}

The restriction of a Sobolev map to a subdomain is a Sobolev map. Conversely, let $\Omega_1, \Omega_2\subset\R^n$ be bounded, disjoint Lipschitz domains and let $W$ be a common boundary component of $\Omega_1$ and $\Omega_2$. Then  $\Omega =
\Omega_1 \cup \Omega_2 \cup W$ is a Lipschitz domain.  If $u_i \in W^{1,p} (\Omega_i, X)$, $i=1, 2$,  are such that $\trace(u_1) = \trace(u_2)$ almost everywhere on $W$ then the map $u$ defined as $u_i$ on $\Omega _i$ is in 
$W^{1,p}(\Omega, X)$, see \cite[Theorem 1.12.3]{KS93}.

The following lemma will be needed in Section~\ref{sec:higher-integrability-area-min}.

\bl\label{lem:gluing-Sobolev}
Let $\Omega, \Omega'\subset\R^n$ be bounded Lipschitz domains with $\Omega'\subset \Omega$. Let $p>1$ and let $u\in W^{1,p}(\Omega ,X)$ and  $v\in W^{1,p} (\Omega',X)$ be such that $\trace(v) = \trace(u|_{\Omega'})$ almost everywhere. Then the map $\bar{u}\colon\Omega \to X$ which coincides with $v$ on $\Omega'$ and with $u$ on $\Omega\setminus\overline{\Omega'}$ is in $W^{1,p}(\Omega, X)$ and satisfies $\trace(\bar{u}) =\trace(u)$ almost everywhere.  
\el

We note that the lemma will only be used in the case that $\Omega$ is an open ball.

\begin{proof}
There exist a neighborhood $\Omega _0$ of $\partial \Omega$ and a biLipschitz homeomorphism $\varphi\colon\Omega_0\to\Omega_0$ with the following properties. The set $\Omega_0$ is a Lipschitz domain decomposed by $\partial \Omega$ in two connected components, and $\varphi$ fixes $\partial \Omega$ and exchanges the two connected components.
Set $\Omega ^+ := \Omega _0\setminus \overline{\Omega}$. Then the map $u^+=u\circ \varphi\colon\Omega ^+\to X$ is contained in 
$W^{1,p} (\Omega^+,X)$ and satisfies $\trace(u^+)|_{\partial\Omega} = \trace(u)$. Set $\tilde{\Omega}:= \Omega_0\cup\Omega$. 
By the paragraph preceding the lemma, the map $w\colon \tilde{\Omega} \to X$ which coincides with $u$ on $\Omega$ and with $u^+$ on $\Omega^+$  is contained in $W^{1,p}(\tilde{\Omega}, X)$. Since $\Omega'$ and $\tilde{\Omega}\setminus\overline{\Omega'}$ are Lipschitz domains and $\trace(w|_{\tilde{\Omega}\setminus\overline{\Omega'}}) = \trace(v)$ almost everywhere on $\partial\Omega'$ it follows again from the paragraph above that the map which coincides with $v$ on $\Omega'$ and with $w$ on $\tilde{\Omega}\setminus\overline{\Omega'}$ is contained in $W^{1,p}(\tilde{\Omega}, X)$. The restriction of this map to $\Omega$ is exactly $\bar{u}$ and satisfies $\trace(\bar{u}) = \trace(u)$.
\end{proof}

\section{Differentiability properties of Sobolev maps}\label{sec:Lip-prop-Sobolev}

The aim of this section is to establish some differentiability properties of Sobolev maps which will be used in the rest of the paper.
Throughout this section, $\Omega\subset \R^n$ will be an open, bounded subset and $(X,d)$ a complete metric space. 

Recall that for a map $u\colon\Omega\to X$ the metric directional derivative of $u$ at $z\in\Omega$ in direction $v\in\R^n$ is defined by
\begin{equation*}
 \md u_z(v):= \lim_{r\to 0^+} \frac{d(u(z+rv), u(z))}{r}
\end{equation*}
if the limit exists. It was shown in \cite{Kir94} that if $u$ is Lipschitz then for almost every $z\in\Omega$ the metric directional derivative $\md u_z(v)$ exists for all $v\in\R^n$ and defines a seminorm on $\R^n$. 
The following notion of approximate metric differentiability, which already appears in \cite{Kar07}, will be useful in the sequel.

\bd
 A map $u\colon \Omega\to X$ is called approximately metrically differentiable at $z\in\Omega$ if there exists a seminorm $s$ on $\R^n$ such that 
 \begin{equation*}
  \ap\lim_{z'\to z}\frac{d(u(z'), u(z)) - s(z'-z)}{|z'-z|} = 0.
 \end{equation*}
\ed
For the definition of approximate limit see e.g.\ \cite{EG92}. The seminorm, if it exists, is unique and will be called the approximate metric derivative of $u$ at $z$ and denoted by $\apmd u_z$. It is straight-forward to check that the following holds.

\br
{\rm If $u$ is Lipschitz then $u$ is approximately metrically differentiable at $z$ if and only if the metric directional derivative $\md u_z(v)$ exists for all $v\in\R^n$ and $\md u_z$ is a seminorm. In this case one has $\apmd u_z = \md u_z$.}
\er

Every classical Sobolev function $u\in W^{1,p}(\Omega)$ is approximately differentiable at almost every $z\in\Omega$ and thus also approximately metrically differentiable at $z$ with $$\apmd u_z(v) = |\ap d_zu(v)|$$ for every $v\in\R^n$. Here, $\ap d_zu$ denotes the approximate derivative of $u$ at $z$. It was proved in \cite{Kar07} that Sobolev maps to metric spaces are approximately metrically differentiable almost everywhere. We prove the following stronger result.

\bp\label{prop:Sobolev-apmd-Lip}
 Let $p>1$ and $u\in W^{1,p}(\Omega, X)$. Then 
 \begin{enumerate}
  \item $u$ is approximately metrically differentiable at almost every $z\in\Omega$ and $z\mapsto \apmd u_z$ is measurable as a map to $\mathfrak{S}_n$; moreover, the function $z\mapsto \mathcal{I}_+^1(\apmd u_z)$ is in $L^p(\Omega)$;
  \item there exist countably many compact, pairwise disjoint sets $K_i\subset\Omega$, $i\in\N$, such that $\lm^n(\Omega\setminus\cup K_i)=0$ and such that the following property holds: for every $i\in\N$ and every $\varepsilon>0$ there exists $r_i(\varepsilon)>0$ such that $u$ is approximately metrically differentiable at every $z\in K_i$ and
 \begin{equation*}
  |d(u(z+v), u(z+w)) - \apmd u_z(v-w)|\leq \varepsilon|v-w|
 \end{equation*}
 for every $z\in K_i$ and all $v,w\in \R^n$ with $|v|, |w|\leq r_i(\varepsilon)$ and such that $z+v, z+w\in K_i$.  
 \end{enumerate}
\ep

It is worth mentioning that for almost every $z\in\Omega$ and every $v\in \R^n$ we have
\begin{equation}\label{eq:apmd-direct-energy}
 \apmd u_z(v) = |u_*(v)|(z),
\end{equation}
where $|u_*(v)|(z)$ is the directional energy-density function defined in \cite[Theorem 1.9.6]{KS93}. This follows from Proposition~\ref{prop:Sobolev-apmd-Lip} together with \cite[Lemma 1.9.5, Theorem 1.8.1]{KS93}.

\begin{proof}
 After possibly writing $\Omega$ as the countable union of (closed) cubes and restricting $u$ to a fixed open cube we may assume that $\Omega$ is an open cube and thus is bounded and has Lipschitz boundary.
 By Proposition~\ref{prop:equiv-Hajlasz-KS} there exist $g\in L^p(\Omega)$ and a negligible set $N\subset\Omega$ such that 
 \begin{equation*}
   d(u(z), u(z')) \leq |z-z'| (g(z) + g(z'))
 \end{equation*} 
 for all $z,z'\in \Omega\setminus N$. For $j\geq 1$ define $A_j:= \{z\in\Omega\setminus N: g(z)\leq j\}$ and note that $u|_{A_j}$ is $(2j)$-Lipschitz. Clearly, we have $\lm^n(\Omega\setminus \cup A_j) = 0$. 
 
 Denote by $\ell^\infty(X)$ the Banach space of bounded functions on $X$, endowed with the supremum norm. Using a Kuratowski embedding, we may view $X$ as a subset of $\ell^\infty(X)$. 
 Fix $j\geq 1$ and let $\bar{u}\colon \Omega\to \ell^\infty(X)$ be a Lipschitz extension of $u|_{A_j}$. By \cite[Theorem 2]{Kir94}, the metric derivative $\md \bar{u}_z(v)$ exists for almost every $z\in\Omega$ and for all $v\in\R^n$ and $\md \bar{u}_z$ is a seminorm. Moreover, there exist compact subsets $K'_i\subset\Omega$, $i\in\N$, such that $\md\bar{u}_z$ exists and is a seminorm for all $z\in K'_i$, such that $\lm^n(\Omega\setminus\cup K'_i) = 0$ and the following holds: for every $i\in\N$ and every $\varepsilon>0$ there exists $r'_i(\varepsilon)>0$ such that 
 \begin{equation}\label{eq:Lip-prop-ubar}
  \left|\,\|\bar{u}(z+v) - \bar{u}(z+w)\|_\infty - \md\bar{u}_z(v-w)\right| \leq \varepsilon|v-w|
 \end{equation}
 for every $z\in K'_i$ and all $v,w\in\R^n$ with $|v|, |w|\leq r'_i(\varepsilon)$ such that $z+w\in K'_i$; see \cite[Theorem 2.3]{Wen08-tree} for this variant of \cite[Theorem 2]{Kir94}. From this it follows that for every $z\in A_j\cap K'_i$, every $\varepsilon>0$, and every $0<r\leq r'_i(\varepsilon)$ we have
 \begin{equation*}
  \left\{z'\in B(z,r)\cap\Omega: \frac{|d(u(z'), u(z)) - \md\bar{u}_z(z'-z)|}{|z'-z|} >\varepsilon\right\} \subset B(z,r)\setminus(A_j\cap K'_i).
 \end{equation*}
 In particular, if $K_i$ denotes the Lebesgue density points of $A_j\cap K'_i$ then $u$ is approximately metrically differentiable at every $z\in K_i$ and $\apmd u_z = \md\bar{u}_z$. Since the map $z\mapsto \md\bar{u}_z$ is measurable as a limit of measurable maps it follows that $z\mapsto \apmd u_z$ is measurable as a map from $K_i$ to $\mathfrak{S}_n$.
 By \eqref{eq:Lip-prop-ubar}, we moreover obtain that
 \begin{equation*}
   |d(u(z+v), u(z+w)) - \apmd u_z(v-w)|\leq \varepsilon|v-w|
 \end{equation*}
 for every $z\in K_i$ and all $v,w\in \R^n$ with $|v|, |w|\leq r'_i(\varepsilon)$ and such that $z+v, z+w\in K_i$. In particular, if $z\in K_i$ is such that $g$ is approximately continuous at $z$ then $$\apmd u_z(v)\leq 2g(z) |v|$$ for every $v\in\R^n$. Now, statements (i) and (ii) easily follow since $\lm^n(A_j\setminus \cup K_i) =0$ and $\lm^n(\Omega\setminus \cup A_j)=0$. Note that the $K_i$ may be taken to be compact and pairwise disjoint by passing to smaller sets.
\end{proof}

\br\label{rem:apmd-bounded-by-g} {\rm
 The proof shows, in particular, that if $u\in W^{1,p}(\Omega, X)$ and $g\in L^p(\Omega)$ is such that \eqref{eq:Hajlasz-Sobolev} holds then $$\apmd u_z(v)\leq 2g(z) |v|$$ for almost every $z\in\Omega$ and every $v\in\R^n$.}
\er

Using the approximate metric differentiability of Sobolev maps, we can extend the definition of the parametrized volume given in Section~\ref{sec:def-vol-normed} to metric space valued Sobolev maps as follows. Let $\mu$ be a definition of volume as in Definition~\ref{def:volume-def} and recall the notion of Jacobian $\jac^\mu_n(s)$ with respect to $\mu$ of a seminorm $s$ on $\R^n$. 

\bd\label{def:vol-sobolev}
 The parametrized $\mu$-volume of a map $u\in W^{1,n}(\Omega, X)$ is defined by
\begin{equation*}
 \Vol_\mu(u):= \int_{\Omega} \jac_n^\mu(\apmd u_z)\,d\lm^n(z). 
\end{equation*}
When $n=2$ we will write $\Area_\mu(u)$ instead of $\Vol_\mu(u)$.
\ed

If $\mu$ is the Busemann definition of volume and $n=2$ then $\Area_\mu(u)$ becomes the parametrized $2$-dimensional Hausdorff measure which was simply denoted by $\Area(u)$ in the introduction.

The Korevaar-Schoen energy can be represented using the approximate metric derivative as follows.

\bp\label{prop:rep-energy}
 Let $p>1$ and $u\in W^{1,p}(\Omega, X)$. Then the Korevaar-Schoen energy density measure of $u$ is given by 
 \begin{equation}\label{eq:energy-density-KS}
 de^p(\,\cdot\,, u) = \mathcal{I}^p_{avg}(\apmd u)\,d\lm^n
 \end{equation} and, in particular, the Korevaar-Schoen energy of $u$ is
 \begin{equation}\label{eq:rep-energy}
  E^p(u) =  \int_\Omega\mathcal{I}^p_{avg}(\apmd u_z)\,d\lm^n(z).
 \end{equation}
\ep
This is a direct consequence of \cite[(1.10ii)]{KS93} together with \eqref{eq:apmd-direct-energy}. For the convenience of the reader, we provide a self-contained proof which relies on Propositions~\ref{prop:equiv-Hajlasz-KS} and \ref{prop:Sobolev-apmd-Lip} instead.

\begin{proof}
 Define $f(z):= \mathcal{I}^p_{\rm avg}(\apmd u_z)$ for almost every $z\in\Omega$. We calculate
 \begin{equation*}
   \frac{ |e_r^p(z,u) - f(z)|}{n+p} = r^{-p} \left|\;\vint_{B(z,r)}d^p(u(z), u(z')) - \apmd u_z(z'-z)^p\,d\lm^n(z')\right|
 \end{equation*}
 and thus obtain with Propositions~\ref{prop:equiv-Hajlasz-KS} and \ref{prop:Sobolev-apmd-Lip} that
 \begin{equation*}
  \begin{split}
  \frac{ |e_r^p(z,u) - f(z)|}{n+p}
   & \leq p\varepsilon (2g(z) + \varepsilon)^{p-1} + (2^p+2^{p-1}) g(z)^p \frac{\lm^n(B(z,r)\setminus K_i)}{\lm^n(B(z,r))}\\
   &\quad + 2^{p-1} \frac{1}{\lm^n(B(z,r))}\int_{B(z,r)\setminus K_i} g(z')^p\,d\lm^n(z')
  \end{split}
 \end{equation*}
 for almost every $z\in K_i$, every $\varepsilon>0$, and every $0<r<r_i(\varepsilon)$. Here, $K_i$ and $r_i(\varepsilon)$ are as in Proposition~\ref{prop:Sobolev-apmd-Lip}. It follows that $e_r^p(z,u)$ converges to $f(z)$ as $\varepsilon, r\searrow 0$ for almost every Lebesgue density point $z$ of $K_i$. Vitali's convergence theorem thus yields that for every $\varphi\in C_c(\Omega)$ we have
 \begin{equation*}
  \lim_{r\to 0} E^p_r(\varphi, u) = \int_{\Omega} \varphi(z)\mathcal{I}^p_{\rm avg}(\apmd u_z)\,d\lm^n(z).
 \end{equation*}
 This proves \eqref{eq:rep-energy} and shows that the energy density measure of $u$ is given by \eqref{eq:energy-density-KS}.
\end{proof}

Now and for Lemma~\ref{lem:approx-apmd} below we assume that $X$ is moreover separable.
Fix a countable, dense subset $\{x_i\}_{i\in\N}\subset X$. For every $N\in\N$ define a map $\varphi_N\colon X\to\ell_N^\infty$ by 
\begin{equation}\label{eq:Lip-map-fin-dim-integrand}
 \varphi_N(x):= (d(x,x_1), \dots, d(x,x_N)),
\end{equation}
where $d$ is the metric on $X$. Here, $\ell_N^\infty$ denotes $\R^N$ endowed with the sup-norm $\|\cdot\|_\infty$.
 Note that $\varphi_N$ is $1$-Lipschitz for every $N\in\N$. We will need the following auxiliary result in later sections. 

\bl\label{lem:approx-apmd}
 Let $u\in W^{1,p}(\Omega, X)$ with $p>1$ and define $u_N:= \varphi_N\circ u$. Then for almost every $z\in\Omega$ we have $$\apmd (u_N)_z(v) \nearrow \apmd u_z(v)\quad\text{as $N\to\infty$}$$ uniformly in $v\in S^{n-1}$. In particular, $\apmd (u_N)_z$ converges to $\apmd u_z$ with respect to the metric $d_{\mathfrak{S}_n}$.
\el

The lemma above implies that for almost every $z\in\Omega$ we have
\begin{equation}\label{eq:apmd-vs-nabla-u}
 \apmd u_z(v) = \sup_{i\in\N} |\langle \nabla u_{x_i}(z), v\rangle|
\end{equation}
for every $v\in\R^n$, where $u_{x_i}(z):= d(x_i, u(z))$. This together with Proposition~\ref{prop:rep-energy} yields the representation
$$E^p(u) = n\int_\Omega\vint_{S^{n-1}}\sup_{i\in\N} |\langle \nabla u_{x_i}(z), v\rangle|^p\,d\hm^{n-1}(v)\,d\lm^n(z)$$ for the  Korevaar-Schoen energy, thus providing a different proof of the main result in \cite{LS12}.

\begin{proof}
 We first note that for all $z,z'\in\Omega$ we have that 
 \begin{equation*}
  \|u_N(z') - u_N(z)\|_\infty \nearrow d(u(z'), u(z))
 \end{equation*}
 as $N\nearrow\infty$. From this it follows that for almost every $z\in\Omega$ and every $v\in S^{n-1}$ the sequence $(\apmd (u_N)_z(v))$ is non-decreasing with
 \begin{equation}\label{eq:ap-u_N}
 \lim_{N\to\infty} \apmd (u_N)_z(v) \leq \apmd u_z(v).
 \end{equation}
 Let $v\in S^{n-1}$ be fixed. We show that for almost every $z\in\Omega$ equality holds in \eqref{eq:ap-u_N}. Define $f(z):= \lim_{N\to\infty} \apmd (u_N)_z(v)$. It follows from Proposition~\ref{prop:ac-1-dim} below that for almost every $z\in\Omega$ and almost all $s<t$ in $\R$ satisfying $\{z+rv: r\in[s,t]\}\subset\Omega$ we have $$\|u_N(z+tv) - u_N(z+sv)\|_\infty \leq \int_s^t\apmd(u_N)_{z+rv}(v)\,dr\leq  \int_s^tf(z+rv)\,dr$$ for every $N\in\N$. From this we obtain that $$d(u(z+tv), u(z+sv)) \leq  \int_s^tf(z+rv)\,dr.$$ Hence, from Proposition~\ref{prop:Sobolev-apmd-Lip} and the Lebesgue differentiation theorem, we conclude that $$\apmd u_z(v) \leq f(z)$$ for almost every $z\in\Omega$. This proves that for fixed $v\in S^{n-1}$ equality holds in \eqref{eq:ap-u_N} for almost every $z\in\Omega$. From this it easily follows that for almost every $z\in\Omega$ we have that $$\apmd (u_N)_z(v) \nearrow \apmd u_z(v)$$ uniformly in $v\in S^{n-1}$. This completes the proof.
 \end{proof}

Apart from the Korevaar-Schoen energy $E^p(u)$ we will also extensively use the following energy functional, which will play a crucial role in our paper. Given $u\in W^{1,p}(\Omega, X)$ we define the energy $E_+^p(u)$ of $u$ by $$E^p_+(u):= \int_\Omega \mathcal{I}^p_+(\apmd u_z)\,d\lm^n(z),$$
where $\mathcal{I}^p_+$ is as in Section~\ref{sec:prelim}. This energy will be of particular importance in Sections~\ref{sec:qc-energy-min} and \ref{sec:existence-area-min}.
It is not difficult to see that $E^p_+$ is precisely the energy defined in \cite{Res97} when $X$ is separable and thus we call $E^p_+(u)$ the Reshetnyak energy. Indeed, let $\{x_i\}_{i\in\N}\subset X$ be a countable, dense subset. Then for every $i\in\N$
\begin{equation*}
 |\nabla u_{x_i}(z)|^p\leq \mathcal{I}^p_+(\apmd u_z)
\end{equation*}
for almost every $z\in\Omega$, where $u_{x_i}(z):= d(u(z), x_i)$. Hence, \eqref{eq:apmd-vs-nabla-u} shows that
\begin{equation*}
 \mathcal{I}^p_+(\apmd u_z) = \sup_{i\in\N} |\nabla u_{x_i}(z)|^p
\end{equation*}
for almost every $z\in\Omega$ and thus the energy considered in \cite{Res97} is precisely given by $E^p_+(u)$. It can furthermore be proved that $E_+^p(u)$ is the integral of $\rho_u^p$, where $\rho_u$ is the minimal weak upper gradient (of a Newtonian representative) of $u$, see \cite[Theorem 7.1.20]{HKST15}. We finally note that $E^p_+(\cdot)$ is related to the Korevaar-Schoen energy by 
\begin{equation}\label{eq:rel-E+-Eavg}
 n^{-1} E^p(u)\leq E^p_+(u) \leq \lambda E^p(u),
\end{equation}
where $\lambda>0$ is a constant only depending on $n$ and $p$. Moreover, one has the equality $n^{-1} E^p(u)= E^p_+(u)$ if and only if $u$ is conformal. These properties follow from Lemma~\ref{lem:int-seminorm} and the remark after Definition~\ref{def:def-qc-seminorm}.

Throughout the remainder of this section, let $X$ be a complete metric space. The following lemmas will be useful. 

\bp\label{prop:ac-1-dim}   
 Let $I=(a,b)$ be an interval and let $u\in W^{1,p}(I, X)$ with $p>1$. Then $u$ has an absolutely continuous representative $\bar{u}$ which satisfies 
  \begin{equation}\label{eq:length-apmd-1-dim}
  \length_X(\bar{u}) = \int_a^b\apmd u_t(1) dt.
 \end{equation}
In particular, the length function $t\mapsto \length_X(\bar{u}|_{(a,t)})$ is contained in $W^{1,2}(I,\R)$.
\ep

\begin{proof}
 This is a  consequence of Lemmas 1.9.2 and 1.9.3 in \cite{KS93}.  Alternatively, it can be proved as follows. The existence of an absolutely continuous representative $\bar{u}$ of $u$ is a consequence of Proposition~\ref{prop:equiv-Hajlasz-KS} above and the fact that $X$ is complete. Then the  representation \eqref{eq:length-apmd-1-dim} follows from the proof of \cite[Theorem 4.1.6]{AT04} together with the fact that $\apmd u_t(1) = \md \bar{u}_t(1)$ for almost every $t\in(a,b)$.
\end{proof}

\bl\label{lem:Sobolev-bilip-apmd}
 Let $\Omega'\subset\R^n$ be a bounded, open set and $\varphi:\Omega'\to\Omega$ a biLipschitz map. If $u\in W^{1,p}(\Omega, X)$ for some $p>1$ then $u\circ\varphi\in W^{1,p}(\Omega', X)$ and 
 \begin{equation*}
  \apmd (u\circ\varphi)_z = \apmd u_{\varphi(z)}\circ d_z\varphi
 \end{equation*}
 for almost every $z\in \Omega'$.
\el

\begin{proof}
 This is a straight-forward consequence of the existence of approximate metric derivatives almost everywhere proved in Proposition~\ref{prop:Sobolev-apmd-Lip}.
\end{proof}

The following proposition, which can essentially be obtained from arguments in \cite{KS93}, will be used repeatedly. Let $\varphi\colon I\times U\to \Omega$ be a biLipschitz map, where $I=(a,b)$ is an open interval and $U\subset\R^{n-1}$ an open set.
For $r\in U$ denote by $\gamma_r$ the curve in $\Omega$ given by $\gamma_r(t):=\varphi(t, r)$ for $t\in I$.

\bp
Let $p>1$ and $u\in W^{1,p}(\Omega, X)$. Then $u\circ \gamma_r\in W^{1,p}(I, X)$ for almost every $r\in U$ and the length of the continuous representative of
$u\circ\gamma_r$ is given by $$\length_X(u\circ\gamma_r) = \int_a^b\apmd
u_{\gamma_r(t)}(\dot{\gamma}_r(t))\,dt$$ for almost every $r\in U$.
\ep

An analogous statement holds when $I$ is replaced by $S^1$. We provide a direct proof which does not rely on the results in \cite{KS93}.

\begin{proof}
By Lemma~\ref{lem:Sobolev-bilip-apmd} we have that $u\circ\varphi\in W^{1,p}(I\times U, X)$ with 
\begin{equation}\label{eq:bilip-fam-curves-apmd}
 \apmd (u\circ\varphi)_{(t,r)} = \apmd u_{\gamma_r(t)}\circ d_{(t,r)}\varphi
\end{equation}
 for almost every $t$ and $r$. Proposition~\ref{prop:equiv-Hajlasz-KS} and Fubini's theorem imply that  $u\circ \gamma_r\in W^{1,p}(I, X)$ for almost every $r\in U$. For such $r$ the absolutely continuous representative of $u\circ\gamma_r$, again denoted by $u\circ\gamma_r$, satisfies $$\length_X(u\circ\gamma_r) = \int_a^b\apmd (u\circ \gamma_r)_t(1)\,dt$$
by Proposition~\ref{prop:ac-1-dim}. Finally, Proposition~\ref{prop:Sobolev-apmd-Lip} together with \eqref{eq:bilip-fam-curves-apmd} yield that for almost all $t$ and $r$ we have $$\apmd (u\circ\gamma_r)_t(1) = \apmd (u\circ\varphi)_{(t,r)}(e_1)= \apmd u_{\gamma_r(t)}(\dot{\gamma}_r(t)),$$ where $e_1=(1,0,\dots, 0)\in\R^n$, completing the proof.
\end{proof}

We end the section with the following useful result.

\bl\label{lem:char-Sobolev-loc}
  Let $p>1$. Then $u\in W^{1,p}(\Omega, X)$ if and only if $u\in L^p(\Omega, X)$ and $u\in W^{1,p}_{\rm loc}(\Omega, X)$ with $$\int_\Omega \mathcal{I}_+^p(\apmd u_z)\,d\lm^n(z)<\infty.$$ If $\Omega$ is convex then the hypothesis that $u$ be in $L^p(\Omega, X)$ is not needed in the `if' part. 
\el

\begin{proof}
 This is a direct consequence of Proposition~\ref{prop:Reshetnyak-Sobolev} and the corresponding classical statement of the proposition when $X=\R$.
\end{proof}

\section{Weak lower semi-continuity of generalized integrands}\label{sec:weak-lsc-integrands}

The aim of this section is to establish a general lower semi-continuity result for functionals on Sobolev maps with values in a metric space, Theorem~\ref{thm:weak-lsc-qc-fctl}. This theorem will be used to show that the volume functionals of many volume definitions are lower semi-continuous, see Corollary~\ref{cor:quasi-convex-volume-lsc}. This in turn will be used in Section~\ref{sec:existence-area-min} to prove the existence of area-minimizing Sobolev maps. Theorem~\ref{thm:weak-lsc-qc-fctl} can furthermore be used to give new proofs of the lower semi-continuity of the Korevaar-Schoen and Reshetnyak energies, see Corollaries~\ref{cor:lsc-KS-energy} and \ref{cor:lsc-+-energy}.

We first recall from Section~\ref{sec:prelim} that $\mathfrak{S}_n$ denotes the space of seminorms on $\R^n$ and that $\mathfrak{S}_n$ is endowed with the metric coming from the supremum norm on $C^0(S^{n-1}, \R)$.

\bd
 A function $\mathcal{I}\colon \R^n\times \mathfrak{S}_n\to [0,\infty)$ is called generalized integrand on $\R^n$ if $\mathcal{I}(\cdot, s)$ is measurable for every $s\in\mathfrak{S}_n$ and $\mathcal{I}(z,\cdot)$ is continuous for almost every $z\in\R^n$. The function $\mathcal{I}$ is said to have bounded $p$-growth if there exist $h\in L^1_{\rm loc}(\R^n)$ and $C\geq 0$ such that $$\mathcal{I}(z, s) \leq h(z) + C\mathcal{I}_+^p(s)$$ for almost every $z\in\R^n$.
\ed

Functions $\mathcal{I}\colon \mathfrak{S}_n \to [0, \infty)$ can and will be naturally identified with functions $\R^n\times\mathfrak{S}_n\to [0,\infty)$ independent of the first variable. 
The following elementary lemma will be needed in the sequel.

\bl
 Let $\Omega\subset\R^n$ be an open, bounded subset, $X$ a complete metric space, and $p>1$. If $u\in W^{1,p}(\Omega, X)$ and if $\mathcal{I}\colon\R^n\times \mathfrak{S}_n\to [0,\infty)$ is a generalized integrand on $\R^n$ then the function $z\mapsto \mathcal{I}(z, \apmd u_z)$ is measurable. Moreover, if $\mathcal{I}$ has bounded $p$-growth then 
 \begin{equation*}
  \int_\Omega\mathcal{I}(z, \apmd u_z)\,d\lm^n(z) \leq C'(1+ E^p(u)),
 \end{equation*}
 where $C'$ is a constant depending on $\mathcal{I}$, $\Omega$, $n$, and $p$.
 \el

\begin{proof}
 Since $\mathfrak{S}_n$ is separable, one may show exactly as in the proof of \cite[Proposition VIII.1.1]{ET76} that there exists a Borel function $\tilde{\mathcal{I}}\colon \R^n\times\mathfrak{S}_n\to [0,\infty)$ such that $\tilde{\mathcal{I}}(z,\cdot) = \mathcal{I}(z,\cdot)$ for almost every $z\in\R^n$. By Proposition~\ref{prop:Sobolev-apmd-Lip}, the map $z\mapsto \apmd u_z$ is measurable. Since $\tilde{I}$ is Borel it thus follows that the function $z\mapsto \tilde{\mathcal{I}}(z, \apmd u_z)$ is measurable. Hence the function $f(z):= \mathcal{I}(z, \apmd u_z)$ is measurable as well. This proves the first part of the proposition.
 
If $\mathcal{I}$ has bounded $p$-growth then there exist $h\in L^1(\R^n)$ and $C\geq 0$ such that $$f(z) \leq h(z) + C \mathcal{I}_+^p(\apmd u_z)$$ for almost every $z\in\Omega$. Lemma~\ref{lem:int-seminorm} and Proposition~\ref{prop:rep-energy} thus imply that 
\begin{equation*}
 \begin{split}
  \int_\Omega f(z)\,d\lm^n(z) &\leq \|h\|_{L^1(\Omega)} + \lambda C \int_\Omega\mathcal{I}_{\rm avg}^p(\apmd u_z)\,d\lm^n(z)\\
  & = \|h\|_{L^1(\Omega)} + \lambda C E^p(u)
 \end{split}
\end{equation*}
for some constant $\lambda>0$ depending only on $n$ and $p$.
This completes the proof.
\end{proof}

We next introduce a variant of the classical quasi-convexity of functions which is adapted to our situation. Given a generalized integrand $\mathcal{I}$ on $\R^n$ we define
\begin{equation*}
 \mathcal{F}_{\mathcal{I}}(u):= \int_\Omega \mathcal{I}(z, \apmd u_z)\,d\lm^n(z)
\end{equation*} 
whenever $u\in W^{1,p}(\Omega, X)$, where $\Omega\subset\R^n$ is an open, bounded subset and $X$ a complete metric space. 

\bd\label{def:quasi-convex-generalized}
 A continuous function $\mathcal{I}\colon \mathfrak{S}_n\to[0,\infty)$ is called quasi-convex if for every finite dimensional normed space $Y$ and every linear map $L: \R^n\to Y$ we have
 \begin{equation}\label{eq:ineq-cond-quasi-convex}
  \mathcal{F}_{\mathcal{I}}(L|_B) \leq \mathcal{F}_{\mathcal{I}}(\psi)
 \end{equation}
 for every smooth immersion $\psi\colon B\to Y$ with $\psi|_{\partial B} = L|_{\partial B}$, where $B$ denotes the closed unit ball in $\R^n$.
\ed

Note that if $\|\cdot\|$ denotes the norm on $Y$ then, by definition, \eqref{eq:ineq-cond-quasi-convex} becomes
\begin{equation*}
 \lm^n(B)\cdot \mathcal{I}(\|\cdot\|\circ L) \leq \int_B\mathcal{I}(\|\cdot\|\circ d_z\psi)\,d\lm^n(z).
\end{equation*}

A function $F\colon W^{1,p}(\Omega, X)\to\R$ is said to be lower semi-continuous on $W^{1,p}(\Omega, X)$ with respect to weak convergence if $$F(u)\leq \liminf_{j\to\infty} F(u_j)$$ for every $u\in W^{1,p}(\Omega, X)$ and every sequence $(u_j)\subset W^{1,p}(\Omega, X)$ with $\sup_j E^p(u_j)<\infty$ and such that $u_j\to u$ in $L^p(\Omega, X)$. 

In what follows, a function $\mathcal{I}\colon\mathfrak{S}_n\to[0,\infty)$ is called monotone if $\mathcal{I}(s)\leq \mathcal{I}(s')$ for all $s,s'\in\mathfrak{S}_n$ with $s\leq s'$. The main result of the present section can be stated as follows.

\bt\label{thm:weak-lsc-qc-fctl}
 Let $\mathcal{I}$ be a generalized integrand on $\R^n$ and let $p>1$. Suppose $\mathcal{I}$ is of bounded $p$-growth and $\mathcal{I}(z,\cdot)$ is monotone for almost every $z\in\R^n$. Then $\mathcal{I}(z,\cdot)$ is quasi-convex for almost every $z\in\R^n$ if and only if for every open, bounded subset $\Omega\subset\R^n$ and every complete metric space $X$ the functional $$\mathcal{F}_{\mathcal{I}}(u):= \int_\Omega \mathcal{I}(z,\apmd u_z)\,d\lm^n(z)$$ is lower semi-continuous on $W^{1,p}(\Omega, X)$ with respect to weak convergence.
\et

The following proposition will be useful in the proof of the theorem.

\bp\label{prop:fin-dim-proj-approx}
 Let $\mathcal{I}$ be a generalized integrand on $\R^n$ of bounded $p$-growth for some $p>1$. Let $\Omega\subset\R^n$ be an open, bounded subset and $X$ a complete metric space. 
 Then for every $u\in W^{1,p}(\Omega, X)$ and every $\varepsilon>0$ there exists a finite dimensional normed space $Y$ and some $1$-Lipschitz map $\varphi\colon X\to Y$ such that $$|\mathcal{F}_{\mathcal I}(\varphi\circ u) - \mathcal{F}_{\mathcal{I}}(u)|\leq \varepsilon.$$
\ep

\begin{proof}
We first consider the case that $X$ is separable.
 Let $\{x_i\}_{i\in\N}\subset X$ be a countable dense subset and, for $N\in\N$, let $\varphi_N\colon X\to\ell_N^\infty$  be the map defined in \eqref{eq:Lip-map-fin-dim-integrand}.
 Fix $u\in W^{1,p}(\Omega, X)$ and let $f, f_N\colon\Omega\to\R$ be the functions given by $$f(z):= \mathcal{I}(z, \apmd u_z)\quad\text{and}\quad f_N(z):= \mathcal{I}(z, \apmd (\varphi_N\circ u)_z),$$ where $\varphi_N$ is the Lipschitz map defined in \eqref{eq:Lip-map-fin-dim-integrand}. It follows from Lemma~\ref{lem:approx-apmd} and the properties of $\mathcal{I}$  that $f_N(z)$ converges to $f(z)$ for almost every $z\in\Omega$ and $$f_N(z) \leq h(z) + C\mathcal{I}_+^p(\apmd u_z),$$  where $h\in L^1(\R^n)$ and $C\geq 0$ are independent of $N$. By Proposition~\ref{prop:Sobolev-apmd-Lip}, the function $z\mapsto \mathcal{I}_+^p(\apmd u_z)$ is in $L^1(\Omega)$. Thus, by the Lebesgue dominated convergence theorem, it follows that $f_N$ converges to $f$ in $L^1(\Omega)$. From this the statement of the proposition follows with $Y=\ell_N^\infty$ and $\varphi:= \varphi_N$, where $N\in\N$ is chosen large enough. This proves the proposition in the case that $X$ is separable. 
 
 We now treat the general case. After possibly changing $u$ on a set of measure zero we may assume that $u$ has separable image. Let $X'$ denote the closure of the image of $u$. Let $\varepsilon>0$. By the first case, there exists a $1$-Lipschitz map $\varphi'\colon X'\to\ell_N^\infty$ such that $$|\mathcal{F}_{\mathcal I}(\varphi'\circ u) - \mathcal{F}_{\mathcal{I}}(u)|\leq \varepsilon.$$ Since $\ell_N^\infty$ is an injective metric space there exists a $1$-Lipschitz extension $\varphi\colon X\to\ell_N^\infty$ of $\varphi'$. Since $\mathcal{F}_{\mathcal I}(\varphi\circ u) = \mathcal{F}_{\mathcal I}(\varphi'\circ u)$ this proves the general case.
\end{proof}

We turn to the proof of Theorem~\ref{thm:weak-lsc-qc-fctl}.

\begin{proof}[Proof of Theorem~\ref{thm:weak-lsc-qc-fctl}]
 Suppose first that $\mathcal{I}(z,\cdot)$ is quasi-convex for almost every $z\in\R^n$. Let  $\Omega\subset\R^n$ be an open, bounded subset and let $X$ be a complete metric space. Let $u\in W^{1,p}(\Omega, X)$ and let $(u_j)\subset W^{1,p}(\Omega, X)$ be such that $u_j\to u$ in $L^p(\Omega, X)$ and $\sup_j E^p(u_j)<\infty$. 
We claim that it is enough to show that 
\begin{equation}\label{eq:lsc-pushforward}
 \mathcal{F}_{\mathcal{I}}(\varphi\circ u)\leq \liminf_{j\to\infty} \mathcal{F}_{\mathcal{I}}(\varphi\circ u_j)
\end{equation}
 for every finite dimensional normed space $Y$ and every $1$-Lipschitz map $\varphi\colon X\to Y$. Indeed, let $\varepsilon>0$ and let $\varphi$ be as in Proposition~\ref{prop:fin-dim-proj-approx}. Since $\mathcal{I}(z,\cdot)$ is monotone for almost every $z\in\Omega$ it follows that $\mathcal{F}_{\mathcal{I}}(\varphi\circ u_j)\leq \mathcal{F}_{\mathcal{I}}(u_j)$ for every $j\in\N$. Therefore, if \eqref{eq:lsc-pushforward} holds then we obtain
 \begin{equation*}
   \mathcal{F}_{\mathcal{I}}(u) - \varepsilon \leq  \mathcal{F}_{\mathcal{I}}(\varphi\circ u)  \leq \liminf_{j\to\infty}  \mathcal{F}_{\mathcal{I}}(\varphi\circ u_j) \leq \liminf_{j\to\infty}  \mathcal{F}_{\mathcal{I}}(u_j).
 \end{equation*}
 Since $\varepsilon>0$ was arbitrary it follows that $$\mathcal{F}_{\mathcal{I}}(u) \leq  \liminf_{j\to\infty}  \mathcal{F}_{\mathcal{I}}(u_j).$$
 It remains to be proven that \eqref{eq:lsc-pushforward} holds. However, identifying $Y$ with $\R^N$ via any  linear isomorphism the statement  translates into the classical sequential weak lower semicontinuity statement of quasi-convex functionals as it is stated in Theorem II.4 of \cite{AF84}. This proves the if part of the theorem.
 
The only if part of the theorem follows from Theorem II.5 of \cite{AF84} via identifying a given $N$-dimensional normed space $Y$ with $\R^N$ via any isomorphism and by using Rellich's theorem. 
\end{proof}

 Throughout the rest of this section, let $\Omega\subset \R^n$ be an open, bounded subset and $X$ a complete metric space. Using Theorem~\ref{thm:weak-lsc-qc-fctl} we can give a new proof of the lower semi-continuity statement in \cite[Theorem 1.6.1]{KS93}.

\bc\label{cor:lsc-KS-energy}
 Let $p>1$. Then the Korevaar-Schoen energy $E^p(\cdot)$ is lower semi-continuous on $W^{1,p}(\Omega, X)$ with respect to weak convergence. 
\ec

\begin{proof}
 The function $\mathcal{I}^p_{\rm avg}$ defined in Section~\ref{sec:prelim} is continuous and hence defines a generalized integrand on $\R^n$. Furthermore, $\mathcal{I}^p_{\rm avg}$ is monotone, of bounded $p$-growth, and satisfies $$E^p(u) = \int_\Omega\mathcal{I}^p_{\rm avg}(\apmd u_z)\,d\lm^n(z)$$ for every $u\in W^{1,p}(\Omega, X)$ by Proposition~\ref{prop:rep-energy}. Finally, it is not difficult to see that $\mathcal{I}^p_{\rm avg}$ is quasi-convex in the sense of Definition~\ref{def:quasi-convex-generalized}. Indeed, let $(Y, \|\cdot\|)$ be a finite dimensional normed space and let $B$ denote the closed unit ball in $\R^n$. Let $L\colon B\to Y$ be the restriction of a linear map and let $\psi\colon B\to Y$ be a smooth immersion such that $\psi|_{\partial B} = L|_{\partial B}$. Fix $v\in S^{n-1}$ and denote by $W\subset\R^n$ the subspace orthogonal to $v$. The triangle inequality and Jensen's inequality imply that for every $y\in W$ 
 \begin{equation*}
  \int_{\R} 1_{B}(y+tv)\cdot \|L(v)\|^p\,dt \leq \int_{\R}1_{B}(y+tv)\cdot \|d_{y+tv}\psi(v)\|^p\,dt.
 \end{equation*}
 Hence, Fubini's theorem yields
 \begin{equation}\label{eq:quasi-convexity-fixed-vector}
  \lm^n(B) \|L(v)\|^p \leq \int_{B} \|d_z\psi(v)\|^p\,d\lm^n(z)
 \end{equation}
 and thus
 \begin{equation*}
 \lm^n(B) \cdot \mathcal{I}^p_{\rm avg}(\| \cdot \|\circ L) \leq \int_B \mathcal{I}^p_{\rm avg}(\|\cdot\|\circ d_z \psi)\,d\lm^n(z).
 \end{equation*}
 This shows that $\mathcal{I}^p_{\rm avg}$ is quasi-convex. It thus follows from Theorem~\ref{thm:weak-lsc-qc-fctl} that $E^p(\cdot)$ is lower semi-continuous on $W^{1,p}(\Omega, X)$ with respect to weak convergence.
\end{proof}

In the same way one proves the weak lower semi-continuity of $E_+^p(\cdot)$ and thus partly recovers \cite[Theorem 4.2]{Res97}. 

\bc\label{cor:lsc-+-energy}
Let $p>1$. Then the Reshetnyak energy $E^p_+(\cdot)$ is lower semi-continuous on $W^{1,p}(\Omega, X)$ with respect to weak convergence. 
\ec

As a consequence of Theorem~\ref{thm:weak-lsc-qc-fctl} we have the following result which will be used to prove the existence of area minimizers in Section~\ref{sec:existence-area-min}.
Recall Definition~\ref{def:vol-def-quasi-convex} for the notion of quasi-convex volume densities.

\bc\label{cor:quasi-convex-volume-lsc}
 Let $\mu$ be a definition of volume and $n\geq 1$. If $\mu$ induces quasi-convex $n$-volume densities then $\Vol_\mu(\cdot)$ is lower semi-continuous on $W^{1,n}(\Omega, X)$ with respect to weak convergence. 
\ec

\begin{proof}
The function $\mathcal{I}\colon \mathfrak{S}_n\to[0,\infty)$ given by $\mathcal{I}(s):= \jac_n^\mu(s)$ defines a generalized integrand on $\R^n$ which is monotone, of bounded $n$-growth, and satisfies $$\Vol_\mu(u) = \int_\Omega\mathcal{I}(\apmd u_z)\,d\lm^n(z)$$ for every $u\in W^{1,n}(\Omega, X)$.
Furthermore, $\mathcal{I}$ is quasi-convex in the sense of Definition~\ref{def:quasi-convex-generalized}. Thus, the claim follows from Theorem~\ref{thm:weak-lsc-qc-fctl}.
\end{proof}

Corollary~\ref{cor:quasi-convex-volume-lsc} will be used in Section~\ref{sec:existence-area-min} in order to prove the existence of area minimizers.

\section{Quasi-conformality of energy minimizers}\label{sec:qc-energy-min}

The main purpose of this section is to prove Theorem~\ref{thm:intro-qc-energy-min} from the introduction, which is restated below as Theorem~\ref{thm:qc-domain-minimizers} for convenience and which shows that every energy minimizing maps is (weakly) quasi-conformal with a universal constant. This is well-known when $X$ is Euclidean space, however, the classical proof of this result does not to carry over to the general setting of metric spaces. This comes from the fact that it seems impossible to obtain a good description of variation of the energy if non-Euclidean norms appear as approximate metric derivatives.

Let $X$ be a complete metric space. Recalling from Section~\ref{sec:prelim} the definition of quasi-conformality of seminorms on $\R^n$ we now define:

\bd
A map $u\in W^{1,n}(\Omega, X)$, where $\Omega\subset\R^n$ is open and bounded, is called $Q$-quasi-conformal if $\apmd u_z$ is $Q$-quasi-conformal for almost every $z\in \Omega$. 
\ed

Moreover, $1$-quasi-conformal Sobolev mappings will be called conformal. We note that ($\R^N$-valued) conformal maps according to our definition are called weakly conformal by some authors. 

Denote by $D$ the open unit disc in $\R^2$. The main result of this section is:

\bt\label{thm:qc-domain-minimizers}
 Let $X$ be a complete metric space. Suppose that $u\in W^{1,2}(D, X)$ is such that $$E^2_+(u)\leq E^2_+(u\circ\psi)$$ for every biLipschitz homeomorphism $\psi\colon D\to D$. Then $u$ is $\sqrt{2}$-quasi-conformal.
\et

The proof will furthermore show the following. Suppose that $u$ is as in the theorem and, in addition, $\apmd u_z$ is induced by an inner product for almost every $z\in D$ for which $\apmd u_z$ is non-degenerate. Then $u$ is conformal. The quasi-conformality constant $\sqrt{2}$ is optimal in general as the following remark shows.

\br\label{rem:qc-optimal-constant}
 {\rm Let $\ell^\infty_2$ be the $2$-dimensional plane endowed with the supremum norm. If $u\in W^{1,2}(D,\ell^\infty_2)$ is non-constant then $u$ cannot be better $\sqrt{2}$-quasi-conformal. Indeed, there exists a set $A\subset D$ of positive measure such that $u$ is approximately differentiable with $T:=\ap d_zu$ non-degenerate at each $z\in A$. Let $r>0$ be the largest number so that $T(D)$ contains $rB$, where $B$ is the unit ball of $\ell^\infty_2$. By John's theorem (see Section 3 of \cite{Bal97}) we have $T(D)\not\subset \lambda rB$ for every $\lambda<\sqrt{2}$; thus $\apmd u_z$ cannot be better than $\sqrt{2}$-quasi-conformal.}
\er

In Theorem~\ref{thm:qc-domain-minimizers-KS-energy} we will obtain an analogue of Theorem~\ref{thm:qc-domain-minimizers} for the energy $E^2(\cdot)$. In this case, however, we can only bound the quasi-conformality constant by $2\sqrt{2} + \sqrt{6}$, which is probably not optimal.

Before proving the theorem we establish several auxiliary results. We start with the following easy observation.

\bl\label{lem:qc-biLiphomeo-energy}
 Let $\Omega$ and $\Omega'$ be bounded, open subsets of $\R^n$ and $\varphi\colon\Omega'\to\Omega$ a conformal biLipschitz homeomorphism. Then for every $u\in W^{1,n}(\Omega, X)$ we have $E^n_+(u\circ\varphi) = E^n_+(u)$ and $E^n(u\circ\varphi) = E^n(u)$.
\el

\begin{proof}
 By Lemma~\ref{lem:Sobolev-bilip-apmd}, we have $u\circ\varphi\in W^{1,n}(\Omega', X)$ and 
 \begin{equation*}
   \apmd (u\circ\varphi)_z(v) = \apmd u_{\varphi(z)}(d_z\varphi(v))
 \end{equation*}
 for almost every $z\in\Omega'$ and every $v\in\R^n$. Since $v\mapsto |d_z\varphi(v)|$ is conformal, Lemma~\ref{lem:qc-transform-energy-seminorm}  implies that
 \begin{equation*}
\mathcal{I}^n_+(\apmd (u\circ\varphi)_z)= |\det d_z\varphi|\cdot \mathcal{I}^n_+(\apmd u_{\varphi(z)})
 \end{equation*}
and
  \begin{equation*}
  \mathcal{I}^n_{\rm avg}(\apmd (u\circ\varphi)_z) = |\det d_z\varphi|\cdot  \mathcal{I}^n_{\rm avg}(\apmd u_{\varphi(z)})
 \end{equation*}
 for almost every $z\in\Omega$. The lemma now follows from the change of variables formula.
\end{proof}

The following lemma proves the infinitesimal version of Theorem~\ref{thm:qc-domain-minimizers}. It is a reformulation of Theorem~\ref{thm:qc-domain-minimizers} for linear maps to normed vector spaces $X$.

\bl\label{lem:qc-seminorm}
 Let $s$ be a seminorm on $\R^2$ such that for every $T\in {\rm SL}_2(\R)$ we have 
 \begin{equation}\label{eq:min-+-energy-seminorm}
  \mathcal{I}^2_+(s)\leq \mathcal{I}^2_+(s\circ T).
 \end{equation} 
 Then $s$ is $\sqrt{2}$-quasi-conformal. Moreover, if $s$ induced by an inner product then $s$ is conformal.
\el

The proof will show that a norm $s$ satisfying \eqref{eq:min-+-energy-seminorm} is isotropic in the following sense. The ellipse of maximal area contained in the unit ball with respect to the norm $s$ is a Euclidean disc.

\begin{proof}
 If $s$ is degenerate then it follows from \eqref{eq:min-+-energy-seminorm} that $s\equiv 0$. We may therefore assume that $s$ is non-degenerate. After rescaling $s$, we may also assume that $\mathcal{I}^2_+(s) = 1$. Denote by $B$ the open unit ball with respect to $s$, that is, $$B:= \{v\in\R^2: s(v)<1\}.$$  Since $\mathcal{I}_+ ^2 (s) =1$ we have $D\subset B$. We will show that $D$ is the ellipse of largest area contained in $B$. Arguing by contradiction we assume that there exists $L\in {\rm GL}_2(\R)$ with $|\det L|>1$ and such that $L(D) \subset B$. Set $\lambda:= |\det L|$ and define $T:= \lambda^{-\frac{1}{2}} L$. Then $T\in{\rm SL}_2(\R)$ and $T$ satisfies $$s(T(v))^2 = \lambda^{-1} s(L(v))^2 \leq \lambda^{-1}$$ for every $v\in D$. Thus we have $\mathcal{I}^2_+(s\circ T) \leq \lambda^{-1} < \mathcal{I}^2_+(s)$, contradicting \eqref{eq:min-+-energy-seminorm}. Therefore, no such $L$ exists. It follows from this that $D$ is the ellipse of largest area contained in $B$ and thus, by definition, $D$ is the Loewner ellipse for $B$. Therefore, by John's theorem (see e.g.~\cite[Theorem 2.18]{AlvT04}), we have that $B\subset\sqrt{2} D$ and thus $s(v)\geq 1/\sqrt{2}$ for every $v\in S^1$. This shows that $s$ is indeed $\sqrt{2}$-quasi-conformal. Finally, if $s$ is induced by an inner product then $B$ is itself an ellipse. Since $D\subset B$ and $D$ is the ellipse of largest area contained in $B$ it follows that $B=D$ and so $s$ is conformal.
\end{proof}

The next simple lemma, essentially a consequence of the Lebesgue differentiation theorem, allows us to obtain from an infinitesimal a local energy-decreasing variation.

\bl\label{lem:approx-continuity-I-apmd}
Let $\Omega\subset\R^n$ be an open, bounded subset and $X$ a complete metric space. Let $p>1$ and $u\in W^{1,p}(\Omega, X)$.  Let $\mathcal{I}\colon\mathfrak{S}_n\to[0,\infty)$ be continuous with bounded $p$-growth. Then for almost every $z_0\in\Omega$ 
 \begin{equation}\label{eq:approx-cont-int-apmd}
  \vint_{B(z_0, r)}\mathcal{I}(\apmd u_z\circ T)\,d\lm^n(z) \longrightarrow \mathcal{I}(\apmd u_{z_0}\circ T)
 \end{equation}
as $r\to0$ for every linear map $T\colon\R^n\to\R^n$. 
\el

\begin{proof}
 It is enough to show that for every $\varepsilon>0$ there exists a measurable set $A\subset\Omega$ such that $\lm^n(A)<\varepsilon$ and such that \eqref{eq:approx-cont-int-apmd} hold for every $z_0\in\Omega\setminus A$. 
 
 Let therefore $\varepsilon>0$. Let $f\colon\Omega\to \mathfrak{S}_n$ be the function given by $f(z):= \apmd u_z$ if $\apmd u_z$ exists and $f(z)=0$ otherwise. By Lusin's theorem \cite[2.3.5]{Fed69}, there exists $A\subset \Omega$ measurable with $\lm^n(A)<\varepsilon$ and such that $f|_{\Omega\setminus A}$ is continuous. 
 
We first show that for almost every $z_0\in\Omega\setminus A$ we have 
 \begin{equation}\label{eq:int-conv-0-apmd-I}
  \frac{1}{\lm^n(B(z_0, r))} \int_{B(z_0, r)\cap A}\mathcal{I}(\apmd u_z\circ T)\,d\lm^n(z) \longrightarrow 0
 \end{equation}
 as $r\to 0$ for every linear map $T\colon\R^n\to\R^n$. By Proposition~\ref{prop:Sobolev-apmd-Lip}, there exists $g\in L^p(\Omega)$ such that $\apmd u_z(v)\leq 2g(z) |v|$ for almost every $z\in\Omega$ and every $v\in\R^n$. Since $\mathcal{I}$ is of bounded $p$-growth there exists $C\geq 0$ such that 
 \begin{equation}\label{eq:bound-I-circ-T}
  \mathcal{I}(\apmd u_z\circ T)\leq C + 2^pC g(z)^p \|T\|^p
 \end{equation}
  for almost every $z\in\Omega$ and every $T\colon\R^n\to\R^n$ linear. Here, $\|T\|$ denotes the operator norm of $T$. The Lebesgue differentiation theorem together with \eqref{eq:bound-I-circ-T} immediately yields \eqref{eq:int-conv-0-apmd-I} for almost every $z_0\in\Omega\setminus A$.
  
  Let $z_0\in \Omega\setminus A$ be a Lebesgue density point of $\Omega\setminus A$ and such that \eqref{eq:int-conv-0-apmd-I} holds. We show that \eqref{eq:approx-cont-int-apmd} holds for $z_0$. For this, let $T\colon\R^n\to\R^n$ be linear. Let $\delta>0$. Since $\mathcal{I}$ and $f|_{\Omega\setminus A}$ are continuous there exists $r_0>0$ such that $B(z_0,r_0)\subset\Omega$ and $$|\mathcal{I}(\apmd u_z\circ T) - \mathcal{I}(\apmd u_{z_0}\circ T)|\leq \delta$$ for every $z\in\Omega\setminus A$ with $|z-z_0|\leq r_0$. It thus follows that
 \begin{equation*}
  \begin{split}
  \int_{B(z_0, r)}&|\mathcal{I}(\apmd  u_z\circ T) - \mathcal{I}(\apmd u_{z_0}\circ T)|\,d\lm^n(z) \\
   &\quad \leq \delta\cdot \lm^n(B(z_0,r)\setminus A) + \lm^n(B(z_0, r)\cap A)\cdot \mathcal{I}(\apmd u_{z_0}\circ T)\\
   &\quad\quad +  \int_{B(z_0, r)\cap A}\mathcal{I}(\apmd u_z\circ T)\,d\lm^n(z)\\
  \end{split}
 \end{equation*}
 for every $0<r<r_0$ and hence $$\limsup_{r\to0}\vint_{B(z_0, r)}|\mathcal{I}(\apmd  u_z\circ T) - \mathcal{I}(\apmd u_{z_0}\circ T)|\,d\lm^n(z)\leq \delta.$$ Since $\delta>0$ was arbitrary we conclude that \eqref{eq:approx-cont-int-apmd} holds for $z_0$. This concludes the proof.
\end{proof}

The following elementary lemma is, together with Lemma~\ref{lem:qc-biLiphomeo-energy}, the key to the localization of the variational argument. 

\bl\label{lem:local-deformation}
 Let $T\in {\rm GL}_2(\R)$ and let $z_0\in\R^2$ and $r>0$. Then there exists a biLipschitz homeomorphism $\varrho\colon\R^2\to\R^2$ such that $$\varrho(z) = z_0 + T(z-z_0)$$ for every $z\in \bar{B}(z_0,r)$ and such that $\varrho$ is smooth and conformal on $\R^2\setminus \bar{B}(z_0, r)$.
\el

\begin{proof}
 After a translation and a dilation we may assume that $z_0=0$ and $r=1$. We may furthermore assume that $T$ is diagonal with strictly positive entries. Indeed, by the polar decomposition theorem and by diagonalization, we may write $T=K\cdot \hat{T}\cdot L$ for suitable $K,L\in {\rm O}_2(\R)$ and a diagonal matrix $\hat{T}$ with strictly positive entries.
 
We identify $\R^2$ with $\C$ in the usual way.  We may thus assume that $T$ is given by $T(x+iy) = ax + iby$ for some $a,b>0$ and that $z_0=0$ and $r=1$. Set $D:= B(0,1)\subset\C$ and define a holomorphic function $\varrho\colon \C\setminus D \to \C$ by $\varrho(z)= cz + dz^{-1}$, where $c,d\in\R$ are such that $c+d = a$ and $c-d = b$. Note that $c>0$ and $|d|<c$. Then $\varrho$ satisfies $\varrho(z) = T(z)$ for every $z\in\C$ with $|z|=1$. Moreover, $\varrho$ is injective and satisfies
\begin{equation}\label{eq:derivative-bounds}
c-|d| \leq |\varrho'(z)|\leq c+|d|
\end{equation}
 for all $z$. Since $\varrho$ maps the circle $\{|z|=r\}$ for $r\geq 1$ surjectively onto the ellipse 
\begin{equation*}
\left\{x+iy: \frac{x^2}{(cr+d/r)^2} + \frac{y^2}{(cr-d/r)^2} = 1\right\}
\end{equation*} 
it follows that the image of $\varrho$ is all of $\C\setminus T(D)$. This together with \eqref{eq:derivative-bounds} implies that $\varrho$ is a biLipschitz homeomorphism from $\C\setminus D$ onto $\C\setminus T(D)$. We now extend $\varrho$ to all of $\C$ by setting $\varrho(z) = T(z)$ for $z\in D$. It follows that $\varrho$ satisfies all the desired properties.
\end{proof}

We are finally ready to prove the main theorem of this section.

\begin{proof}[Proof of Theorem~\ref{thm:qc-domain-minimizers}]
 In order to prove that $u$ is $\sqrt{2}$-quasi-conformal it is enough, by Lemma~\ref{lem:qc-seminorm}, to show that for almost every $z_0\in D$ we have 
 \begin{equation}\label{eq:E-ineq-for-qc}
  \mathcal{I}^2_+(\apmd u_{z_0}) \leq \mathcal{I}^2_+(\apmd u_{z_0}\circ T)
 \end{equation}
 for every $T\in {\rm SL}_2(\R)$. 
  By Lemma~\ref{lem:approx-continuity-I-apmd}, we have that for almost every $z_0\in D$
 \begin{equation}\label{lem:approx-continuity-I-apmd-again}
  \vint_{B(z_0, r)}\mathcal{I}^2_+(\apmd u_z\circ T)\,d\lm^n(z) \longrightarrow \mathcal{I}^2_+(\apmd u_{z_0}\circ T)
 \end{equation}
as $r\to0$ for every linear map $T\colon\R^2\to\R^2$. We prove by contradiction that \eqref{eq:E-ineq-for-qc} holds for every $z_0\in D$ for which \eqref{lem:approx-continuity-I-apmd-again} holds. Assume therefore that $z_0\in D$ is such that \eqref{lem:approx-continuity-I-apmd-again} holds but  $\mathcal{I}^2_+(\apmd u_{z_0}\circ T) <\mathcal{I}^2_+(\apmd u_{z_0})$ for some $T\in {\rm SL}_2(\R)$. 
 Let $\delta>0$ be so small that $$\mathcal{I}^2_+(\apmd u_{z_0}\circ T) +3\delta \leq \mathcal{I}^2_+(\apmd u_{z_0}).$$ By \eqref{lem:approx-continuity-I-apmd-again}, there exists $r>0$ such that $B(z_0,r)\subset\subset D$ and such that
 \begin{equation*}
  \vint_{B(z_0, r)}\mathcal{I}^2_+(\apmd u_z\circ T)\,d\lm^2(z) \leq \mathcal{I}^2_+(\apmd u_{z_0}\circ T) +\delta
 \end{equation*}
and
\begin{equation*}
  \vint_{B(z_0, r)}\mathcal{I}^2_+(\apmd u_z)\,d\lm^2(z) \geq \mathcal{I}^2_+(\apmd u_{z_0}) - \delta.
 \end{equation*}
It follows that
 \begin{equation*}
   \vint_{B(z_0, r)}\mathcal{I}^2_+(\apmd u_{z}\circ T)\,d\lm^2(z)\leq \vint_{B(z_0, r)}\mathcal{I}^2_+(\apmd u_z)\,d\lm^2(z) - \delta.
 \end{equation*}
 By Lemma~\ref{lem:local-deformation}, there exists a biLipschitz homeomorphism $\varrho\colon \R^2\to\R^2$ such that $\varrho(z) = z_0 + T^{-1}(z-z_0)$ for all $z\in \bar{B}(z_0, r)$ and such that $\varrho$ is smooth and conformal outside $\bar{B}(z_0,r)$. 
 Let $\varphi\colon D\to \varrho(D)$ be a conformal diffeomorphism. Since $\varrho(S^1)$ is smooth,  $\varphi$ and $\varphi^{-1}$ are smooth up to the boundary. In particular, $\varphi$ is a biLipschitz homeomorphism. Thus, $\psi:= \varrho^{-1}\circ\varphi$ is a biLip\-schitz homeomorphism from $D$ onto itself.
We calculate using Lemma~\ref{lem:qc-biLiphomeo-energy} and the properties of $\varrho$ that
\begin{equation*}
\begin{split}
 E^2_+&(u\circ\psi) \\
 & = \int_D\mathcal{I}^2_+(\apmd (u\circ\psi)_z)\,d\lm^2(z) \\
  &= \int_{\varrho(D)}\mathcal{I}^2_+(\apmd (u\circ\varrho^{-1})_z)\,d\lm^2(z)\\
  &= \int_{D\setminus \bar{B}(z_0, r)}\mathcal{I}^2_+(\apmd u_z)\,d\lm^2(z) + \int_{\varrho(B(z_0, r))}\mathcal{I}^2_+(\apmd u_{\varrho^{-1}(z)}\circ T)\,d\lm^2(z)\\
  &= \int_{D\setminus \bar{B}(z_0, r)}\mathcal{I}^2_+(\apmd u_z)\,d\lm^2(z) + \int_{B(z_0, r)}\mathcal{I}^2_+(\apmd u_z\circ T)\,d\lm^2(z)\\
  &\leq \int_D\mathcal{I}^2_+(\apmd u_z)\,d\lm^2(z) - \delta\lm^2(B(z_0,r))\\
  &= E^2_+(u) - \delta\lm^2(B(z_0,r)).
\end{split}
\end{equation*}
 This is in contradiction with the hypothesis of the theorem. We therefore conclude that \eqref{eq:E-ineq-for-qc} holds for almost every $z_0\in D$. This completes the proof.
 \end{proof}

We have the following analog of Theorem~\ref{thm:qc-domain-minimizers} for the energy considered by Korevaar-Schoen.

\bt\label{thm:qc-domain-minimizers-KS-energy}
 Let $X$ be a complete metric space. Suppose that $u\in W^{1,2}(D, X)$ is such that $$E^2(u)\leq E^2(u\circ\psi)$$ for every biLipschitz homeomorphism $\psi\colon D\to D$. Then $u$ is $Q$-quasi-conformal with $Q= 2\sqrt{2}+\sqrt{6}$.
\et

The proof of Theorem~\ref{thm:qc-domain-minimizers-KS-energy} is analogous to that of Theorem~\ref{thm:qc-domain-minimizers} but uses the following lemma instead of Lemma~\ref{lem:qc-seminorm}. Since the result will not be used in the sequel and the constant we obtain is worse than that for the $E_+$-energy we leave the details of
the proof to the reader.

\bl\label{lem:qc-seminorm-KS-energy}
 Let $s$ be a seminorm on $\R^2$ such that for every $T\in {\rm SL}(2, \R)$ we have $$\mathcal{I}^2_{\rm avg}(s)\leq \mathcal{I}^2_{\rm avg}(s\circ T).$$ Then $s$ is $Q$-quasi-conformal with $Q= 2\sqrt{2}+\sqrt{6}$. Moreover, if $s$ is induced by an inner product then $s$ is conformal.
\el

\begin{proof}
We use the following fact, which can be proved by a straight-forward calculation. Let $\|\cdot\|$ be a norm on $\R^2$ induced by an inner product and let $\bar{Q}\geq 1$.  If $\|\cdot\|$ satisfies $$2\bar{Q} \mathcal{I}^2_{\rm avg}(\|\cdot\|)\leq (\bar{Q}^2+1)\mathcal{I}^2_{\rm avg}(\|\cdot \|\circ T)$$ for every $T\in {\rm SL}_2(\R)$ then $\|\cdot\|$ is $\bar{Q}$-quasi-conformal.

 Let now $s$ be as in the lemma. It is straight-forward to see that if $s$ is degenerate then $s\equiv 0$. We may therefore suppose that $s$ is non-degenerate. If $s$ is induced by an inner product, then $s$ is conformal by the above fact. In general, by John's theorem (see e.g.~\cite[Theorem 2.18]{AlvT04}), there exists a norm $\|\cdot\|$ on $\R^2$ induced by an inner product such that $$\|v\|\leq s(v)\leq \sqrt{2}\|v\|$$ for every $v\in\R^2$. This together with the hypothesis yields that
  \begin{equation*}
  \mathcal{I}^2_{\rm avg}(\|\cdot\|) \leq 2 \mathcal{I}^2_{\rm avg}(\|\cdot\|\circ T)
 \end{equation*}
 for every $T\in{\rm SL}(2, \R)$. The fact above thus implies that $\|\cdot\|$ is $\bar{Q}$-quasi-conformal with $\bar{Q}=2+\sqrt{3}$. Hence, $s$ is $\sqrt{2}\bar{Q}$-quasi-conformal. 
\end{proof}

\section{Existence of area minimizers and quasi-conformality}\label{sec:existence-area-min}

Given a complete metric space $X$ and a Jordan curve $\Gamma\subset X$ we define $\Lambda(\Gamma, X)$ to be the set of all maps $u\in W^{1,2}(D, X)$ such that $\trace(u)$ has a continuous representative which is a weakly monotone parametrization of $\Gamma$. We refer to Section~\ref{sec:basic-notation} for the notion of weakly monotone parametrization.

The main result of the present section can be stated as follows.

\bt\label{thm:existence-qc-area-min-general-mu}
 Let $\mu$ be a definition of volume which induces quasi-convex $2$-volume densities. Let $X$ be a proper metric space and $\Gamma\subset X$ a Jordan curve. If $\Lambda(\Gamma, X)\not=\emptyset$ then there exists $u\in \Lambda(\Gamma, X)$ which satisfies $$\Area_\mu(u) = \inf\left\{\Area_\mu(u'): u'\in \Lambda(\Gamma, X)\right\}$$ and which is $\sqrt{2}$-quasi-conformal.
 \et

Remark~\ref{rem:qc-optimal-constant} shows that the quasi-conformality factor $\sqrt{2}$ cannot be improved in general.
We refer to Definition~\ref{def:vol-sobolev} for the parametrized $\mu$-area $\Area_\mu(u)$ and to Section~\ref{sec:def-vol-normed} for examples of definitions of volume inducing quasi-convex $2$-volume densities. Since the Busemann (Hausdorff) definition of volume induces quasi-convex $2$-volume densities, see Section~\ref{sec:def-vol-normed},  Theorem~\ref{thm:existence-qc-area-min-general-mu}, in particular, implies Theorem~\ref{thm:intro-exist-area-min}.
It is well-known that if $X$ is Euclidean space then every energy minimizer is an area minimizer. We will show in Section~\ref{sec:ET-case} that this is no longer true in the setting of general metric spaces.

The proof of Theorem~\ref{thm:existence-qc-area-min-general-mu} will be given after establishing several auxiliary results.

\bl\label{lem:comparison-energy-volume}
 Let $\mu$ be a definition of volume and $\Omega\subset\R^n$ an open, bounded subset. Let $X$ be a complete metric space and $u\in W^{1,n}(\Omega, X)$. Then
 \begin{equation*}
  \Vol_\mu(u) \leq E^n_+(u)\quad \text{ and }\quad \Vol_\mu(u)\leq C E^n(u),
 \end{equation*}
where $C$ only depends on $n$. If $u$ is $Q$-quasi-conformal then
 \begin{equation*}
  E^n_+(u) \leq  Q^n \Vol_\mu(u) \quad\text{ and }\quad E^n(u) \leq n Q^n \Vol_\mu(u).
 \end{equation*}
 For $n=2$ the constant $C$ can be taken to be $1$.  
\el

In particular, if $u$ is conformal then $$E^n(u) = nE^n_+(u)= n \Vol_\mu(u).$$

\begin{proof}
 Since $\mu$ is monotone it follows that $$\jac^\mu_n(\apmd u_z) \leq \mathcal{I}^n_+(\apmd u_z)$$ for almost every $z\in\Omega$ and thus  $ \Vol_\mu(u) \leq E^n_+(u)$ upon integration. Moreover, integrating the above inequality and using Lemma~\ref{lem:int-seminorm} and Proposition~\ref{prop:rep-energy} yields $\Vol_\mu(u) \leq C E^n(u)$ for some constant $C$ depending only on $n$. If $n=2$ then $C$ be taken to be $1$ by Lemma~\ref{lem:int-seminorm}. This proves the first part of the lemma.
 
If $u$ is $Q$-quasi-conformal then, by the monotonicity of $\mu$, we obtain $$\mathcal{I}^n_+(\apmd u_z) \leq Q^n \jac^\mu_n(\apmd u_z)$$ for almost every $z\in\Omega$ and thus $E^n_+(u) \leq  Q^n \Vol_\mu(u)$ upon integration. Moreover, Lemma~\ref{lem:int-seminorm}, Proposition~\ref{prop:rep-energy} and the inequality above yields $$E^n(u) \leq nE_+^n(u)\leq nQ^n \Vol_\mu(u).$$ This concludes the proof.
\end{proof}

We have the following variant of the Courant-Lebesgue Lemma which is valid for general complete metric spaces. 

\bl\label{lem:Courant-Lebesgue-general-X}
 Let $(X,d)$ be a complete metric space and $u\in W^{1,2}(D, X)$. Let $z_0\in \overline{D}$ and $\delta\in(0,1)$. For each $r\in (0,1)$ let $\gamma_r$ be an arc-length parametrization of  $\{z\in D: |z-z_0|=r\}$. Then there exists $A\subset(\delta, \sqrt{\delta})$ of strictly positive measure such that $u\circ\gamma_r$ has an absolutely continuous representative of length
 \begin{equation}\label{eq:length-Courant-Lebegue-statement}
  \length_X(u\circ\gamma_r) \leq \pi\left(\frac{2E^2(u)}{|\log \delta|}\right)^{\frac{1}{2}}
 \end{equation}
for every $r\in A$. In particular, if $|z|\geq 1-\delta$ then
 \begin{equation}\label{eq:dist-Courant-Lebegue-statement}
  d(\trace(u)(y_r), \trace(u)(z_r))\leq \pi\left(\frac{2E^2(u)}{|\log \delta|}\right)^{\frac{1}{2}}
 \end{equation}
 for almost every $r\in A$, where $y_r$ and $z_r$ are the points in $S^1$ at distance $r$ from $z_0$.
\el

The proof is a straight-forward adaptation of the classical proof for Euclidean spaces. For the classical proof see e.g.~\cite{Dierkes-et-al10}.

Let $X$ be a complete metric space and $\Gamma\subset X$ a Jordan curve. Fix three distinct points $p_1,p_2,p_3\in S^1$ and three distinct points $\bar{p}_1, \bar{p}_2, \bar{p}_3\in \Gamma$. A map $u\in\Lambda(\Gamma, X)$ is said to satisfy the $3$-point condition with respect to $\{p_1,p_2,p_3\}$ and $\{\bar{p}_1,\bar{p}_2, \bar{p}_3\}$ if the continuous representative of $\trace(u)$, again denoted by $\trace(u)$, satisfies 
\begin{equation}\label{eq:3-pt-cond}
\trace(u)(p_i) = \bar{p}_i\quad\text{for $i=1,2,3$.}
\end{equation}

Using Lemma~\ref{lem:Courant-Lebesgue-general-X} one may establish exactly as in the Euclidean case the proposition below.

\bp\label{prop-equi-cont-3-point}
 Let $X$, $\Gamma$, $p_i$, and $\bar{p}_i$ be as above and let $M>0$. Then the family 
 \begin{equation*}
  \left\{\trace(u): \text{$u\in \Lambda(\Gamma, X)$ satisfies the $3$-point condition \eqref{eq:3-pt-cond} and $E^2(u)\leq M$}\right\}
 \end{equation*}
 is equi-continuous.
\ep

In the above, $\trace(u)$ refers to the continuous representative of $\trace(u)$.

\begin{proof}
 This follows as in the Euclidean case except that the classical Courant-Lebesgue Lemma is replaced by Lemma~\ref{lem:Courant-Lebesgue-general-X}. We refer e.g.~to \cite[pp.~257--258]{Dierkes-et-al10} for the proof in the classical case.
\end{proof}

Using the proposition above we can prove:
 
\bp\label{prop:seq-Jordan-equi-bdd-energy}
 Let $X$ be a proper metric space, $\Gamma\subset X$ a Jordan curve, and $(u_j)\subset \Lambda(\Gamma, X)$ a sequence such that $$\sup_j E_+^2(u_j)<\infty.$$ Then there exist $v\in \Lambda(\Gamma, X)$, a subsequence $(u_{j_k})$ and Moebius transformations $\psi_k\colon \overline{D}\to\overline{D}$ such that $u_{j_k}\circ\psi_k$ converges to $v$ in $L^2(D, X)$.
\ep

\begin{proof}
 For each $j\in\N$ let $\psi_j\colon \overline{D}\to\overline{D}$ be a Moebius transformation such that $v_j:= u_j\circ\psi_j$ satisfies the $3$-point condition \eqref{eq:3-pt-cond}. By Lemmas~\ref{lem:int-seminorm} and \ref{lem:qc-biLiphomeo-energy} we have $$E^2(v_j)\leq 2E^2_+(v_j) = 2E_+^2(u_j)$$ and hence $\sup_j E^2(v_j)<\infty$. Fix $x_0\in \Gamma$. By Lemma~\ref{lem:bound-lp-dist-basepoint}, we have
 \begin{equation*}
\sup_{j\in\N}  \int_D d^2(v_j(z), x_0)\,d\lm^2(z)<\infty
 \end{equation*}
 and hence \cite[Theorem 1.13]{KS93} implies that there exist $v\in W^{1,2}(D, X)$ and a subsequence $(v_{j_k})$ which converges to $v$ in $L^2(D, X)$. It remains to show that $v\in\Lambda(\Gamma, X)$. By Proposition~\ref{prop-equi-cont-3-point}, the sequence $(\trace(v_{j_k}))$ is equi-continuous and thus we may assume, after possibly passing to a further subsequence, that $(\trace(v_{j_k}))$ converges uniformly to a continuous map $c\colon S^1\to X$. Then $c$ is a weakly monotone parametrization of $\Gamma$. By \cite[Theorem 1.12.2]{KS93}, the traces $\trace(v_{j_k})$ converge to $\trace(v)$ in $L^2(S^1, X)$; hence $\trace(v)=c$ almost everywhere on $S^1$. This shows that $v\in\Lambda(\Gamma, X)$ and concludes the proof.
\end{proof}

We are ready to prove the main result of the present section.

\begin{proof}[Proof of Theorem~\ref{thm:existence-qc-area-min-general-mu}]
 We first claim that for every $u\in\Lambda(\Gamma, X)$ there exists $v\in \Lambda(\Gamma, X)$ which is $\sqrt{2}$-quasi-conformal and satisfies $$\Area_\mu(v)\leq \Area_\mu(u).$$ For this let $u\in\Lambda(\Gamma, X)$ and define $$\Lambda_u:= \{v\in \Lambda(\Gamma, X): \Area_\mu(v) \leq \Area_\mu(u)\},$$ which is non-empty since $u\in\Lambda_u$. Let $(u_j)\subset\Lambda_u$ be a sequence such that $E^2_+(u_j) \rightarrow m$ as $j\to\infty$, where $$m:=\inf\left\{E^2_+(u'): u'\in\Lambda_u\right\}.$$ 
 By Proposition~\ref{prop:seq-Jordan-equi-bdd-energy} that there exist $v\in\Lambda(\Gamma, X)$, a subsequence $(u_{j_k})$, and Moebius transformations $\psi_k$ such that $v_k:= u_{j_k}\circ\psi_k$ converges to $v$ in $L^2(D, X)$. By Lemma~\ref{lem:qc-biLiphomeo-energy} and Corollary~\ref{cor:lsc-+-energy} we have $$E^2_+(v) \leq \liminf_{k\to\infty} E_+^2(v_k) = \liminf_{k\to\infty} E_+^2(u_{j_k}) = m$$ and Corollary~\ref{cor:quasi-convex-volume-lsc} implies that $$\Area_\mu(v) \leq \liminf_{j\to\infty}\Area_\mu(u_{n_j}) \leq \Area_\mu(u).$$ In particular, $v\in \Lambda_u$ and $E_+^2(v) = m$. For every biLipschitz homeomorphism $\psi\colon D\to D$ we have $v\circ\psi\in \Lambda_u$ and therefore $$E^2_+(v) = m \leq E^2_+(v\circ\psi).$$ Theorem~\ref{thm:qc-domain-minimizers} thus implies that $v$ is $\sqrt{2}$-quasi-conformal. This proves the claim.
 
Let now $(u_j)\subset \Lambda(\Gamma, X)$ be a sequence with $\Area_\mu(u_j)\rightarrow m'$ as $j\to\infty$, where $$m':= \inf\{\Area_\mu(u'): u'\in\Lambda(\Gamma, X)\}.$$ By the first part of the proof, we may assume that each $u_j$ is $\sqrt{2}$-quasi-conformal. In particular, by Lemma~\ref{lem:comparison-energy-volume}, we have $E_+^2(u_j) \leq 2\Area_\mu(u_j)$ for every $j$ and thus $\sup_j E_+^2(u_j)<\infty$.  Thus we obtain as above that after possibly pre-composing with a Moebius transformation and passing to a subsequence, $(u_j)$ converges to some $u\in\Lambda(\Gamma, X)$ in $L^2(D, X)$ and $$\Area_\mu(u) = m'.$$ Finally, by the first part of the proof, we may assume that $u$ is $\sqrt{2}$-quasi-conformal. This concludes the proof. 
\end{proof}

We note that one can use the same methods to obtain the existence of energy minimizers as follows.

\bt\label{thm:existence-energy-min}
 Let $X$ be a proper metric space and $\Gamma\subset X$ a Jordan curve. If $\Lambda(\Gamma, X)\not=\emptyset$ then there exists $u\in \Lambda(\Gamma, X)$ satisfying $$E_+^2(u) = \inf\left\{E_+^2(u'): u'\in\Lambda(\Gamma, X)\right\}.$$ Every such $u$ is $\sqrt{2}$-quasi-conformal.
\et

An analogous result holds for $E_+^2$ replaced by $E^2$ and $\sqrt{2}$ replaced by $2\sqrt{2} + \sqrt{6}$. 

\begin{proof}
 The first statement follows from Proposition~\ref{prop:seq-Jordan-equi-bdd-energy}, Lemma~\ref{lem:qc-biLiphomeo-energy}, and Corollary~\ref{cor:lsc-+-energy}. The second statement is a consequence of Theorem~\ref{thm:qc-domain-minimizers}.
\end{proof}

\section{Interior regularity of area minimizing discs}\label{sec:higher-integrability-area-min}

The aim of this section is to prove Theorem~\ref{thm:int-reg-summary} below, which establishes interior regularity of $\mu$-area minimizers for arbitrary $\mu$ and which generalizes the interior regularity results stated in the introduction. 

We begin by extending Definition~\ref{def:isop-intro} to arbitrary volumes and by giving classes of spaces satisfying the definition.

\bd
 Let $\mu$ be a definition of volume, and $C, l_0>0$. A complete metric space $X$ is said to admit a uniformly $l_0$-local quadratic isoperimetric inequality with constant $C$ for $\mu$ if for every Lipschitz curve $c\colon S^1\to X$  of length $\length_X(c)\leq l_0$ there exists $u\in W^{1,2}(D, X)$ with $$\Area_\mu(u) \leq C \length_X(c)^2$$ and such that $\trace(u)(t) = c(t)$ for almost every $t\in S^1$.
\ed

If the above holds for Lipschitz curves of arbitrary length then $X$ is said to admit a (global) quadratic isoperimetric inequality with constant $C$ for $\mu$.

In what follows, if a choice of definition of volume $\mu$ has been fixed, a uniformly $l_0$-local quadratic isoperimetric inequality with constant $C$ for $\mu$ will simply be called a $(C, l_0)$-isoperimetric inequality. We observe that if $X$ admits a $(C, l_0)$-isoperimetric inequality for some definition of volume $\mu$ then $X$ admits a $(2C, l_0)$-isoperimetric inequality for any other definition of volume because any two definitions of volume induce areas of Sobolev maps which differ by a factor of at most $2$. Note, however, that area minimizers with respect to two different definitions of volume, spanning the same curve, need not have anything to do with each other, see Proposition~\ref{prop:area-min-diff}.

Many interesting classes of spaces admit a uniformly local quadratic isoperimetric inequality. This includes homogeneously regular Riemannian manifolds in the sense of \cite{Mor48}, compact Lipschitz manifolds and, in particular, all compact Finsler mani\-folds. It furthermore includes complete metric spaces all of whose balls of radius at most $l_0$ are $\gamma$-Lipschitz contractible with fixed $\gamma$ in the sense of \cite{Wen07}. In particular, this applies to complete ${\rm CAT}(\kappa)$ spaces, $\kappa\in\R$, and compact Alexandrov spaces by \cite{PP93}, and, in fact, also to non-compact volume non-collapsed Alexandrov spaces, cf.~\cite{PP93}. It moreover applies to complete metric spaces with a convex bicombing in the sense of \cite{Wen05} or \cite{Lan13} and, in particular, to all Banach spaces and all injective metric spaces. Further examples of spaces admitting a uniformly local quadratic isoperimetric inequality are given by the Heisenberg groups $\mathbb{H}^n$ of topological dimension $2n+1$ for $n\geq 2$, endowed with a Carnot-Carath\'eodory distance. This follows e.g.~from \cite{Alc98}. In all the spaces mentioned above the isoperimetric filling of a Lip\-schitz curve $c$ is given by a Lip\-schitz map of $\mu$-area bounded by $C\length_X(c)^2$ for a suitable constant $C$. In the case of spaces satisfying the local $\gamma$-Lipschitz contractibility condition mentioned above, the constant $C$ depends only on $\gamma$.

The following theorem summarizes our main results concerning the interior regularity of area minimizing discs.

\bt\label{thm:int-reg-summary}
 Let $X$ be a complete metric space admitting a uniformly local quadratic isoperimetric inequality with constant $C$. Let $\mu$ be a definition of volume and let $Q\geq 1$. If $u\in W^{1,2}(D, X)$ is $Q$-quasi-conformal and satisfies 
   \begin{equation*}
   \Area_\mu(u) = \inf\left\{\Area_\mu(v): v\in W^{1,2}(D,X), \trace(v) = \trace(u)\text{ a.e.}\right\}
  \end{equation*}
  then the following statements hold:
 \begin{enumerate}
  \item There exists $p>2$ such that $u\in W^{1,p}_{\rm loc}(D, X)$; in particular, $u$ has a continuous representative $\bar{u}$ which moreover satisfies Lusin's property (N).
  \item The representative $\bar{u}$ is locally $\alpha$-H\"older continuous with $\alpha = (4\pi Q^2C)^{-1}$.
 \end{enumerate}
\et

Note that no assumption is made on $\mu$ and no local compactness condition is made on $X$. Statement (i) and the first part of statement (ii) of Theorem~\ref{thm:reg-intro} are consequences of Theorem~\ref{thm:int-reg-summary}. The following example shows that, in the case $Q=1$, the H\"older exponent $\alpha=\frac{1}{4\pi C}$ is optimal.

\begin{example}\label{example:Hoelder-exp-opt}
 Let $S\subset S^2$ be a round circle of radius $r\in(0,1]$ in the unit sphere $S^2\subset\R^3$ and let $X$ be the cone over $S$, endowed with the intrinsic metric. Then $X$ is a complete metric space admitting a global quadratic isoperimetric inequality with constant $C= \frac{1}{4\pi r}$ for any definition of volume $\mu$, see e.g.~\cite{MR02}. Let $\varphi\colon S^1\to S$ be a natural identification (one which stretches lengths by a constant factor). Then the map $u\colon D\to X$ given by $u(0) = 0$ and $u(z) = |z|^{r}\varphi(z/|z|)$ if $z\not=0$ is in the Sobolev space $W^{1,2}(D, X)$, it is conformal and satisfies
 \begin{equation*}
   \Area_\mu(u) = \inf\left\{\Area_\mu(v): v\in W^{1,2}(D,X), \trace(v) = \trace(u)\text{ a.e.}\right\}.
  \end{equation*}
 Moreover, $u$ is $r$-H\"older continuous but not $s$-H\"older continuous for any $s>r$.
\end{example}

Statement (i) of Theorem~\ref{thm:int-reg-summary} will be proved in Proposition~\ref{prop:hoelder-cont-min} while statement (ii) follows from Proposition~\ref{prop:strong-cont?}. The proof of Proposition~\ref{prop:hoelder-cont-min} uses the isoperimetric inequality in conjunction with a strengthening of Gehring's lemma. The proof of Proposition~\ref{prop:strong-cont?} follows the classical approach of Morrey and uses, in particular, Morrey's growth lemma.

Throughout the remainder of this section, let $\mu$ be a definition of volume, let $C>0$, $l_0>0$, $Q\geq 1$, and let $X$ be a complete metric space $X$ admitting a $(C, l_0)$-isoperimetric inequality. Unless otherwise stated, let $u\in W^{1,2}(D,X)$ be $Q$-quasi-conformal and such that
  \begin{equation*}
   \Area_\mu(u) = \inf\left\{\Area_\mu(v): v\in W^{1,2}(D,X), \trace(v) = \trace(u)\text{ a.e.}\right\}.
  \end{equation*}

Our first proposition establishes higher integrability of $u$.

\bp\label{prop:hoelder-cont-min}
There exists $p>2$ such that $u\in W^{1,p}_{\rm loc}(D, X)$. In particular, $u$ has a locally H\"older continuous representative $\bar{u}$, and $\bar{u}$ satisfies Lusin's property (N).
\ep

The proof is based on the local isoperimetric inequality. We first establish two lemmas, the first of which shows that short curves with a $W^{1,2}$-parametrization have an isoperimetric filling.

\bl\label{lem:isop-Sobolev-curves}
 Let $c\colon S^1\to X$ be a continuous curve with $\length_X(c)\leq l_0$. If $c\in W^{1,2}(S^1, X)$ then there exists $v\in W^{1,2}(D, X)$ with $$\Area_{\mu}(v)\leq C\length_X(c)^2$$ and such that $\trace(v) = c$ almost everywhere on $S^1$.
\el

Regarding the notation $c\in W^{1,2}(S^1, X)$ we refer to the terminology introduced in the paragraph preceding Proposition~\ref{prop:Reshetnyak-Sobolev}.

\begin{proof}
 We may assume that $l:=\length_X(c)>0$. 
 Let $\bar{c}\colon S^1\to X$ be the constant speed parametrization of $c$. By the local isoperimetric inequality there exists $w\in W^{1,2}(D, X)$ with $$\Area_\mu(w)\leq C\length_X(c)^2$$ and such that $\trace(w) = \bar{c}$ almost everywhere on $S^1$. Let $v\colon \bar{B}(0,2)\to X$ be the map which coincides with $w$ on $D$ and which gives a `linear' reparametrization from $c$ to $\bar{c}$ on the annulus $\bar{A}:=\bar{B}(0,2)\setminus D$. 
 More precisely, view $c$ and $\bar{c}$ as curves parametrized on $[0,1]$ by composing with the map $s\mapsto e^{2\pi is}$ and let $\varrho\colon[0,1]\to [0,1]$ be the normalized length function given be $\varrho(s):= l^{-1}\cdot \length_X(c|_{[0,s]})$. Define $v$ on $D$ by $v:=w$ and define $v$ on $\bar{A}$ by $$v(r e^{2\pi is}):= \bar{c}\left((2-r)s + (r-1)\varrho(s)\right).$$ Since $\varrho\in W^{1,2}((0,1))\cap C^0([0,1])$ by Proposition~\ref{prop:ac-1-dim} and since $c = \bar{c} \circ \varrho$ it follows that $v|_{\bar{A}}\in W^{1,2}(A, X)\cap C^0(\bar{A}, X)$ and that $v$ coincides with $c$ on the outer boundary of $A$ and with $\bar{c}$ on the inner boundary of $A$. Thus, Lemma~\ref{lem:gluing-Sobolev} implies that $v\in W^{1,2}(B(0,2), X)$ and $\trace(v)(2z) = c(z)$ for almost every $z\in S^1$. Since $\Area_\mu(v|_A) = 0$ we moreover have that $$\Area_\mu(v)\leq C\length_X(c)^2.$$ Identifying $B(0,2)$ with $D$ via the scaling map we obtain the desired isoperimetric filling of $c$.
\end{proof}

For the next lemma, let $u$ satisfy the hypotheses stated in the paragraph preceding Proposition~\ref{prop:hoelder-cont-min}. Actually, the quasi-conformality condition on $u$ is not needed for this lemma.

\bl\label{lem:isop-u-area-min}
 Let $\Omega\subset D$ be a domain enclosed by some biLipschitz curve in $\overline{D}$. If $\trace(u|_\Omega)$ has a continuous representative, denoted by $u|_{\partial \Omega}$, such that $u|_{\partial \Omega}\in W^{1,2}(\partial \Omega, X)$ and if $\length_X(u|_{\partial \Omega})\leq l_0$ or $\Area_\mu(u|_\Omega)\leq Cl_0^2$ then $$\Area_\mu(u|_\Omega) \leq C\length_X(u|_{\partial \Omega})^2.$$
\el

\begin{proof}
Let $\Omega$ be as in the statement and let $\beta\colon S^1\to\partial\Omega$ be a  biLipschitz homeomorphism. Suppose that $\trace(u|_\Omega)$ has a continuous representative, which we denote by $u|_{\partial\Omega}$, such that $u|_{\partial \Omega}\circ\beta\in W^{1,2}(S^1, X)$ and that $\length_X(u|_{\partial \Omega})\leq l_0$ or $\Area_\mu(u|_\Omega)\leq Cl_0^2$. 

We may assume that $\length_X(u|_{\partial \Omega})\leq l_0$ since otherwise the statement is trivially true.
By Lemma~\ref{lem:isop-Sobolev-curves}, there exists $v\in W^{1,2}(D, X)$ such that $$\Area_\mu(v)\leq C \length_X(u|_{\partial \Omega})^2$$ and such that $\trace(v) = u|_{\partial \Omega}\circ\beta$ almost everywhere. Let $\varphi\colon \overline{D}\to \overline{\Omega}$ be a biLip\-schitz map extending $\beta$. Such $\varphi$ exists by \cite[Theorem A]{Tuk80}.  
By Lemma~\ref{lem:gluing-Sobolev}, the map $\bar{u}\colon D\to X$ which agrees with $u$ on $D\setminus \Omega$ and with $v\circ\varphi^{-1}$ on $\Omega$ is contained in $W^{1,2}(D, X)$ and satisfies $\trace(\bar{u}) = \trace(u)$ almost everywhere. Since $u$ is an area-minimizer it follows that
\begin{equation*}
\Area_\mu(u) \leq \Area_\mu(\bar{u}) \leq \Area_\mu(u|_{D\setminus \Omega}) + C\length_X(u|_{\partial \Omega})^2
\end{equation*}
and hence $\Area_\mu(u|_{\Omega})\leq C\length_X(u|_{\partial \Omega})^2$. This proves the proposition.
\end{proof}

Using Lemma~\ref{lem:isop-u-area-min} we can give the proof of Proposition~\ref{prop:hoelder-cont-min} as follows.

\begin{proof}[Proof of Proposition~\ref{prop:hoelder-cont-min}]
Let $r_0>0$ be such that $\Area_\mu(u|_{D\cap B(z_0,2r_0)})\leq Cl_0^2$ for every $z_0\in D$.
We first show that the function $f(z):= \mathcal{I}_+^1(\apmd u_z)$ satisfies the local weak reverse H\"older inequality
\begin{equation}\label{eq:reverse-Hoelder-apmd}
 \left(\vint_W f^2(z)\,d\lm^2(z)\right)^{\frac{1}{2}} \leq C_1\vint_{2W}f(z)\,d\lm^2(z)
\end{equation}
for some constant $C_1$ and for every square $W$ of edge length at most $2r_0$ such that $2W\subset D$, where $2W$ denotes the square with same center as $W$ but twice the edge length. For this, fix a square $W$ centered at some point $z_0\in D$ and of edge length $2r$, where $r\leq r_0$, in such a way that $2W\subset D$. 
For almost every $0<s<2r$ the map $\trace(u|_{B(z_0,s)})$ has an absolutely continuous representative, denoted by $u|_{\partial B(z_0,s)}$, such that  $u|_{\partial B(z_0,s)}\in W^{1,2}(\partial B(z_0,s), X)$, and such that
 \begin{equation*}
  \length_X(u|_{\partial B(z_0,s)}) = \int_{\partial B(z_0,s)} \apmd u_z(v(z))\,d\hm^1(z) \leq \int_{\partial B(z_0,s)} f(z)\,d\hm^1(z),
 \end{equation*}
 where $v(z)\in S^1$ is the vector orthogonal to $z-z_0$. Hence, Lemma~\ref{lem:isop-u-area-min} shows that
\begin{equation*}
 \Area_\mu(u|_W)\leq \Area_\mu(u|_{B(z_0, s)})\leq C\left(\int_{\partial B(z_0,s)} f(z)\,d\hm^1(z)\right)^2
\end{equation*}
for almost every $\sqrt{2}r<s<2r$ and thus, Lemma~\ref{lem:comparison-energy-volume} yields 
\begin{equation*}
\begin{split}
 \left(\int_Wf^2(z)\,d\lm^2(z)\right)^{\frac{1}{2}} &\leq Q\Area_{\mu}(u|_W)^{\frac{1}{2}}\\
  &\leq Q\sqrt{C} \vint_{\sqrt{2}r}^{2r}\int_{\partial B(z_0,s)} f(z)\,d\hm^1(z)\,d\lm^1(s)\\
  &\leq Q\sqrt{C}(2-\sqrt{2})^{-1} r^{-1} \int_{2W}f(z)\,d\lm^2(z).
 \end{split}
\end{equation*}
Thus, inequality \eqref{eq:reverse-Hoelder-apmd} holds with a constant $C_1$ depending only on $C$ and $Q$. Since $f$ satisfies \eqref{eq:reverse-Hoelder-apmd} a strengthening of Gehring's lemma, see e.g. Theorem~1.5 in \cite{Kin94}, implies that there exists $p>2$ such that $f\in L^p_{\rm loc}(D)$.  This together with Proposition~\ref{prop:Reshetnyak-Sobolev} and the Sobolev inequality implies that $u\in W^{1,p}_{\rm loc}(D, X)$. 
The remaining statements of the proposition follow from Proposition~\ref{prop:Hoelder-p>n}. This concludes the proof.
\end{proof}

The following proposition shows that the  continuous representative of $u$ is locally $\alpha$-H\"older continuous with $\alpha = \frac{1}{4\pi Q^2C}$. As already mentioned in the introduction, the H\"older continuity is obtained by finding curves in $D$  whose images in $X$ have small lengths.
Before stating the proposition we define for $z_1, z_2\in D$
$$A(z_1, z_2):= B(z_1, |z_1-z_2|)\cap B(z_2, |z_1-z_2|)$$ and denote the closure of $A(z_1, z_2)$ by $\bar{A}(z_1, z_2)$. Note that 
\begin{equation}\label{eq:vol-Az-disc}
\diam(A(z_1,z_2))=\sqrt{3} |z_1-z_2|\quad\text{ and }\quad \lm^2(A(z_1,z_2)) \geq \frac{\pi}{3}|z_1-z_2|^2.
\end{equation}

As before, we assume that $u$ satisfies the hypotheses stated in the paragraph preceding Proposition~\ref{prop:hoelder-cont-min}. 

\bp\label{prop:strong-cont?}
If $u\in  C^0(D, X)$ then for every $0<\delta<1$ and all $z_1,z_2\in \bar{B}(0, \delta)$ there exists a piecewise affine curve $\gamma$ in $\bar{A}(z_1, z_2)\cap \bar{B}(0, \delta)$ from $z_1$ to $z_2$ such that $$\length_X(u\circ\gamma)  \leq L\cdot |z_1-z_2|^\alpha,$$
where $\alpha = \frac{1}{4\pi Q^2C}$ and where $L$ does not depend on $z_1$ and $z_2$.  In particular, for every $0<\delta<1$ the restriction of $u$ to $\bar{B}(0,\delta)$ is $\alpha$-H\"older continuous.
\ep

If $X$ admits a global quadratic isoperimetric inequality or if $\Area_{\mu}(u)\leq C l_0^2$ then $L$ depends only on $E^2(u)$, $\alpha$, and $\delta$ and is increasing in $E^2(u)$ and $\delta$.

For the proof of Proposition~\ref{prop:strong-cont?} we need the following two lemmas, essentially due to Morrey. The first lemma gives a  bound on energy growth on balls in $D$ and the second lemma relates the energy of balls with lengths of some curves.

\bl\label{lem:growth-int-conf-factor}
If $z_0\in D$ and $0< r_0\leq 1-|z_0|$ are such that $\Area_\mu(u|_{B(z_0,r_0)})\leq Cl_0^2$ then
\begin{equation*}
 \int_{B(z_0,r)} \mathcal{I}^1_+(\apmd u_z)\,d\lm^2(z)\leq \left[\pi E^2_+(u|_{B(z_0, s)})\right]^{\frac{1}{2}}s^{-\alpha} r^{1+\alpha}
\end{equation*}
for all $0\leq r\leq s\leq r_0$, where $\alpha:= \frac{1}{4\pi Q^2 C}$.
\el

In particular, for every $0\leq r\leq r_0$ we have
\begin{equation*}
 \int_{B(z_0,r)} \mathcal{I}^1_+(\apmd u_z)\,d\lm^2(z)\leq \frac{\left[\pi E^2_+(u)\right]^{\frac{1}{2}}}{r_0^{\alpha}} \cdot r^{1+\alpha}.
\end{equation*}

\begin{proof}
 As in the proof of Proposition~\ref{prop:hoelder-cont-min} for almost every $0<r<r_0$ the map $\trace(u|_{B(z_0,r)})$ has an absolutely continuous representative, denoted by $u|_{\partial B(z_0,r)}$, such that  $u|_{\partial B(z_0,r)}\in W^{1,2}(\partial B(z_0,r), X)$, and such that
 \begin{equation*}
  \length_X(u|_{\partial B(z_0,r)}) \leq \int_{\partial B(z_0,r)} \mathcal{I}_+^1(\apmd u_z)\,d\hm^1(z).
 \end{equation*}
 For such $r$, Jensen's inequality yields
 \begin{equation}\label{eq:length-curve-prop-regularity}
  \length_X(u|_{\partial B(z_0,r)})^2 \leq 2\pi r \int_{\partial B(z_0,r)} \mathcal{I}^2_+(\apmd u_z)\,d\hm^1(z).
 \end{equation}
 It follows from Lemma~\ref{lem:isop-u-area-min} that
 \begin{equation*}
  \begin{split}
   \Area_\mu(u|_{B(z_0,r)}) &\leq C\length_X(u|_{\partial B(z_0,r)})^2\leq 2\pi C\cdot  r \cdot \int_{\partial B(z_0,r)} \mathcal{I}^2_+(\apmd u_z)\,d\hm^1(z)\\
    &= 2\pi C \cdot r \cdot \frac{d}{dr} E^2_+(u|_{B(z_0,r)})
  \end{split}
 \end{equation*}
   for almost every $0<r<r_0$. From Lemma~\ref{lem:comparison-energy-volume} we thus obtain
  \begin{equation}
   E^2_+(u|_{B(z_0,r)}) \leq Q^2 \cdot \Area_\mu(u|_{B(z_0,r)})\leq Q^2\cdot 2\pi C \cdot r \cdot \frac{d}{dr} E^2_+(u|_{B(z_0,r)}).
  \end{equation}
 Hence, upon integration, we get
   \begin{equation*}
    E^2_+(u|_{B(z_0,r)}) \leq \frac{E^2_+(u|_{B(z_0,s)})}{s^{2\alpha}} \cdot r^{2\alpha}
   \end{equation*}
   for all $0< r \leq s\leq r_0$. By H\"older's inequality,
   \begin{equation*}
     \int_{B(z_0,r)} \mathcal{I}^1_+(\apmd u_z)\,d\lm^2(z) \leq  \left[\pi E^2_+(u|_{B(z_0,r)})\right]^{\frac{1}{2}} r \leq \left[\pi E^2_+(u|_{B(z_0, s)})\right]^{\frac{1}{2}}s^{-\alpha} r^{1+\alpha}.
   \end{equation*}
This completes the proof.
\end{proof}

\bl\label{lem:length-bound-riesz-pot}
 Let $v\in W^{1,2}(D,X)\cap C^0(D,X)$ and let $A\subset D$ be a convex subset with $\lm^2(A)>0$. Then for all $z_1, z_2\in A$ there exists a piecewise affine curve $\gamma$ in $A$, joining $z_1$ with $z_2$, and such that
 \begin{equation}\label{eq:length-upper-bound-rieszpot}
  \length_X(v\circ\gamma)\leq 2^{-1}\diam(A)^2 \vint_{A}\left(\frac{\mathcal{I}_+^1(\apmd v_z)}{|z-z_1|}+ \frac{\mathcal{I}_+^1(\apmd v_z)}{|z-z_2|}\right)\,d\lm^2(z).
 \end{equation}
\el

\begin{proof}
 Given $z\in A$, denote by $\gamma_z$ the piecewise affine curve in $A$ from $z_1$ to $z_2$ going through $z$. Then for almost every $z\in A$ we have
 \begin{equation*}
  \begin{split}
   \length_X(v\circ \gamma_z) &= \sum_{i=1}^2 \int_0^1\apmd v_{z_i+t(z-z_i)}(z-z_i)\,dt\\
    &\leq \sum_{i=1}^2\int_0^1\mathcal{I}_+^1(\apmd v_{z_i+t(z-z_i)})\cdot  |z-z_i|\,dt.
   \end{split}
   \end{equation*}
 Set $d:= \diam(A)$. For $i=1,2$ we have
  \begin{equation*}
   \begin{split}
    \vint_A&\int_0^1\mathcal{I}_+^1(\apmd v_{z_i+t(z-z_i)})\cdot |z-z_i|\,dt\,d\lm^2(z)\\
     &= \frac{1}{\lm^2(A)} \int_{S^1}\int_0^d1_A(z_i+ sw) s\int_0^s\mathcal{I}_+^1(\apmd v_{z_i+tw})\,dt\,ds\,d\hm^1(w)\\
     &\leq \frac{1}{\lm^2(A)} \int_{S^1}\int_0^d s\int_0^d1_A(z_i+ tw) \mathcal{I}_+^1(\apmd v_{z_i+tw})\,dt\,ds\,d\hm^1(w)\\
     &= \frac{d^2}{2\lm^2(A)} \int_{S^1}\int_0^d1_A(z_i+ tw) \mathcal{I}_+^1(\apmd v_{z_i+tw})\,dt\,d\hm^1(w)\\
     &= \frac{d^2}{2}\vint_A\frac{\mathcal{I}_+^1(\apmd v_z)}{|z-z_i|}\,d\lm^2(z).
   \end{split}
  \end{equation*}
 There thus exists a subset $B\subset A$ of positive measure such that \eqref{eq:length-upper-bound-rieszpot} holds for every $z\in B$. This completes the proof.
\end{proof}

Using the lemmas above we prove Proposition~\ref{prop:strong-cont?}.

\begin{proof}[Proof of Proposition~\ref{prop:strong-cont?}]
Let $r_0>0$ be such that $\Area_\mu(u|_{D\cap B(z_0,r_0)})\leq Cl_0^2$ for every $z_0\in D$.
Let $0<\delta<1$ and let $z_1,z_2\in \bar{B}(0,\delta)$. Define $A:= \bar{A}(z_1,z_2)\cap \bar{B}(0,\delta)$. Suppose first that $|z_1-z_2|\leq \eta:=\min\{1-\delta, r_0\}$. By Lemma~\ref{lem:length-bound-riesz-pot} there exists a piecewise affine curve $\gamma$ in $A$, joining $z_1$ with $z_2$, and such that 
 \begin{equation}\label{eq:length-bound-curve-for-Holder}
  \length_X(u\circ\gamma)\leq 2^{-1}\diam(A)^2 \vint_A\left(\frac{\mathcal{I}_+^1(\apmd u_z)}{|z-z_1|} + \frac{\mathcal{I}_+^1(\apmd u_z)}{|z-z_2|}\right)\,d\lm^2(z).
 \end{equation}
Set $r:=|z_1-z_2|$. By Lemma~\ref{lem:growth-int-conf-factor}, we have for $i=1,2$ that
 \begin{equation*}
  \int_{B(z_i, s)} \mathcal{I}_+^1(\apmd u_z)\,d\lm^2(z) \leq \left[\pi E^2_+(u)\right]^{\frac{1}{2}}\eta^{-\alpha} s^{1+\alpha}
 \end{equation*}
and hence
\begin{equation*}
 \int_{B(z_i, s)\setminus B(z_i, 2^{-1}s)} \frac{\mathcal{I}_+^1(\apmd u_z)}{|z-z_i|}\,d\lm^2(z) \leq 2\left[\pi E^2_+(u)\right]^{\frac{1}{2}}\eta^{-\alpha} s^{\alpha}
\end{equation*}
for every $0\leq s\leq r$. By summing over annuli we obtain
 \begin{equation*}
  \int_{B(z_i, r)}\frac{\mathcal{I}_+^1(\apmd u_z)}{|z-z_i|}\,d\lm^2(z) \leq \frac{ 2\left[\pi E^2_+(u)\right]^{\frac{1}{2}}}{(1-2^{-\alpha})\eta^{\alpha}}\cdot r^\alpha.
 \end{equation*}
 Since $2\lm^2(A) \geq \lm^2(A(z_1, z_2))$ it follows with \eqref{eq:length-bound-curve-for-Holder} and \eqref{eq:vol-Az-disc} that 
 \begin{equation*}
  \length_X(u\circ\gamma)\leq \frac{36E^2_+(u)^{\frac{1}{2}}}{\sqrt{\pi}(1-2^{-\alpha})\eta^{\alpha}}\cdot |z_1-z_2|^\alpha.
 \end{equation*}
 This proves the proposition in the special case that $|z_1-z_2|\leq \min\{1-\delta, r_0\}$. The general case follows from the special case by subdividing the segment from $z_1$ and $z_2$.
 \end{proof}

\section{Continuity up to the boundary of area minimizing discs}\label{sec:cont-boundary}

The main results of this section are Theorems~\ref{thm:bdry-cont-u-classical} and \ref{thm:bdry-Hoelder}. They imply, in particular, the second part of statement (ii) as well as statement (iii) of Theorem~\ref{thm:reg-intro}.

Let $\mu$ be a definition of volume, $C, l_0>0$, and let $X$ be a complete metric space admitting a $(C, l_0)$-isoperimetric inequality.

\bt\label{thm:bdry-cont-u-classical}
 Suppose $u\in W^{1,2}(D,X)\cap C^0(D, X)$ is quasi-conformal and satisfies
  \begin{equation*}
   \Area_\mu(u) = \inf\left\{\Area_\mu(v): v\in W^{1,2}(D,X), \trace(v) = \trace(u)\text{ a.e.}\right\}.
  \end{equation*}
 If $\trace(u)$ has a continuous representative then the map $\bar{u}\colon\overline{D}\to X$ defined by $$\bar{u}(x):= \left\{\begin{array}{ll}
   u(x) & x\in D\\
   \trace(u)(x) & x\in S^1
  \end{array}\right.$$
  is continuous.
\et

The proof relies on the following estimate, which will be applied to a repara\-metrized piece of the area minimizer.

\bl\label{lem:unif-cont-lower-bound-energy}
 Let $(X, d)$ be a complete metric space, $\varrho\in(0,1)$ and $\varepsilon>0$. Then for every $v\in W^{1,2}(D, X)\cap C^0(D, X)$ satisfying
  \begin{enumerate}
   \item $d(\trace(v)(z), v(0))\geq \varepsilon$ for almost every $z\in S^1$, and 
   \item $d(v(z), v(0))< \varepsilon/2$ for all $z\in D$ with $|z|<\varrho$
  \end{enumerate} we have $$E^2(v|_{v^{-1}(B(v(0), \varepsilon))})\geq \frac{\pi\varrho \varepsilon^2}{8}.$$
\el

\begin{proof}
 By continuity of $v$ and hypothesis (i) there exists for almost every $w\in S^1$ some number $\bar{\varrho}\in(0,1)$ such that $d(v(0), v(\bar{\varrho}w))\geq \frac{3\varepsilon}{4}$.  Let $\bar{\varrho}(w)$ be the smallest such $\bar{\varrho}$ and observe that $\bar{\varrho}(w)>\varrho$ and $d(v(0), v(\bar{\varrho}(w)w))=\frac{3\varepsilon}{4}$. Thus for every such $w$ we have 
 \begin{equation}\label{eq:unif-cont-energy-dist-length}
 \frac{\varepsilon}{4} \leq d(v(\varrho w), v(\bar{\varrho}(w)w)) \leq \length_X(u\circ\gamma_w),
 \end{equation}
where $\gamma_w\colon [\varrho, \bar{\varrho}(w)] \to D$ is the affine curve given by $\gamma_w(r):= rw$. For almost every $w\in S^1$ the curve $u\circ\gamma_w$ is absolutely continuous and satisfies $$\length_X(v\circ\gamma_w) = \int_\varrho^{\bar{\varrho}(w)} \apmd v_{\gamma_w(r)}(w)\,dr \leq \int_\varrho^{\bar{\varrho}(w)} \mathcal{I}_+^1(\apmd v_{rw})\,dr.$$ This together with \eqref{eq:unif-cont-energy-dist-length} and H\"older's inequality implies $$\int_\varrho^{\bar{\varrho}(w)} \mathcal{I}_+^2(\apmd v_{rw})\,dr \geq  \left(\frac{\varepsilon}{4}\right)^2.$$
 Integrating in polar coordinates and using the fact that $v\circ\gamma_w(r)\in B(v(0), \varepsilon)$ for almost every $w\in S^1$ and all $r\in[0,\bar{\varrho}(w))$ we conclude that $$E_+^2(v|_{v^{-1}(B(v(0), \varepsilon))})\geq \int_{S^1}\int_\varrho^{\bar{\varrho}(w)} r\cdot \mathcal{I}_+^2(\apmd v_{rw})\,dr\,d\hm^1(w)\geq  \frac{\pi\varrho\varepsilon^2}{8}.$$ The claim now follows with Lemma~\ref{lem:int-seminorm}.
\end{proof}

\begin{proof}[Proof of Theorem~\ref{thm:bdry-cont-u-classical}]
 Let $\bar{u}$ be defined as in the statement of the theorem. Note that $\bar{u}|_D$ and $\bar{u}|_{S^1}$ are continuous. In order to prove that $\bar{u}$ is continuous on $\overline{D}$ it thus suffices to show that the auxiliary function $\bar{r}\colon D\setminus\{0\}\to[0,\infty)$ given by $$\bar{r}(x):= d(\bar{u}(x), \bar{u}(x/|x|))$$ satisfies $\bar{r}(x)\to 0$ as $|x|\to 1$. 
 
In order to show this we will apply Lemma~\ref{lem:unif-cont-lower-bound-energy} to a suitable map $v$ defined below. Let $Q\geq 1$ be such that $u$ is $Q$-quasi-conformal. Set $\alpha:= \frac{1}{4\pi Q^2C}$ and let $L$ be the constant from Proposition~\ref{prop:strong-cont?} in the case noted after the proposition with the parameters $E^2(u)$, $\alpha$, and where the $\delta$ appearing there is to be taken to equal $\frac{1}{2}$. Let $\varepsilon\in(0,1)$ and set $\varrho:= \min\left\{[(2L)^{-1} \varepsilon]^{\frac{1}{\alpha}}, \frac{1}{2}\right\}$. Choose $\delta\in(0,1)$ so small that $$\pi\cdot\left(\frac{2E^2(u)}{|\log\delta|}\right)^{\frac{1}{2}} < \varepsilon$$ and such that $d(\bar{u}(z), \bar{u}(z'))< \varepsilon$ for all $z,z'\in S^1$ with $|z-z'|<\sqrt{\delta}$ and $$E^2\left(u|_{D\cap B(z, \sqrt{\delta})}\right)< \min\left\{\frac{\pi\varrho\varepsilon^2}{8}, Cl_0^2\right\}$$ for every $z\in S^1$. We claim that $\bar{r}(x)<3\varepsilon$ for every $x\in D$ with $|x|>1-\delta$. Suppose this is wrong and fix an $x$ for which this fails. By Lemma~\ref{lem:Courant-Lebesgue-general-X}, there exists $r\in(\delta, \sqrt{\delta})$ such that the curve $\gamma\colon(a_1, a_2)\to D$ parametrizing $\{z\in D: |z - x/|x||=r\}$ satisfies $$\length_X(u\circ\gamma)\leq \pi\cdot\left(\frac{2E^2(u)}{|\log\delta|}\right)^{\frac{1}{2}}<\varepsilon$$ and $\lim_{t\to a_i} u\circ\gamma(t) = \bar{u}(\lim_{t\to a_i}\gamma(t))$ for $i=1,2$.
 Set $\Omega:= D\cap B(x/|x|, r)$. From the choice of $\delta$ and $r$ it follows that 
 \begin{equation}\label{eq:energy-u-omega-boundary}
  E^2(u|_\Omega) \leq E^2(u|_{D\cap B(x/|x|, \sqrt{\delta})}) < \min\left\{\frac{\pi\varrho\varepsilon^2}{8}, Cl_0^2\right\}
 \end{equation}
  and that $\trace(u|_\Omega)$ has a continuous representative, simply given by $\bar{u}|_{\partial \Omega}$, whose image is contained in the ball $B(x_1,  2\varepsilon)$ with center $x_1:= \bar{u}(x/|x|)$.
  Let $\varphi\colon D\to \Omega$ be a conformal diffeomorphism which maps the origin to $x$. Then the map $v:= u\circ\varphi$ is continuous and satisfies $v\in W^{1,2}(D, X)$ by Lemma~\ref{lem:char-Sobolev-loc} and is $Q$-quasi-conformal with $E^2(v) = E^2(u|_\Omega)$. If we can show that $v$ satisfies the hypotheses of Lemma~\ref{lem:unif-cont-lower-bound-energy} with $\varepsilon$ and $\varrho$ given as above then Lemma~\ref{lem:unif-cont-lower-bound-energy} and \eqref{eq:energy-u-omega-boundary} yield $$\frac{\pi\varrho\varepsilon^2}{8}\leq E^2(v|_{v^{-1}(B(v(0), \varepsilon))}) \leq E^2(v)= E^2(u|_\Omega) <\frac{\pi\varrho\varepsilon^2}{8},$$ which is impossible. Therefore, we must have $\bar{r}(x)\leq 3\varepsilon$ as claimed.

 It thus remains to show that $v$ satisfies the hypotheses of Lemma~\ref{lem:unif-cont-lower-bound-energy}. As already mentioned, $v$ is continuous and satisfies $v\in W^{1,2}(D, X)$. In order to establish property (i) of Lemma~\ref{lem:unif-cont-lower-bound-energy} note first that $\varphi$ extends to a homeomorphism from $\overline{D}$ to $\overline{\Omega}$ which is locally biLipschitz away from the preimage of the two `corners' of $\partial\Omega$. It follows that $\trace(v)$ has a continuous representative, denoted by the same symbol, satisfying $$\trace(v) = \trace(u|_\Omega)\circ\varphi|_{S^1} = \bar{u}|_{\partial \Omega}\circ\varphi|_{S^1}$$ everywhere. Since $\bar{r}(x)\geq 3\varepsilon$ we obtain that $$d(v(0), \trace(v)(z)) \geq d(u(x), x_1) - d(x_1, \trace(v)(z)) \geq \bar{r}(x) - 2\varepsilon \geq \varepsilon$$ for every $z\in S^1$, showing (i). We will use Proposition~\ref{prop:strong-cont?} to establish (ii). For this we first claim \begin{equation}\label{eq:area-min-v_n}
   \Area_\mu(v) = \inf\left\{\Area_\mu(w): w\in W^{1,2}(D,X), \trace(w) = \trace(v)\text{ a.e.}\right\}.
  \end{equation}
In order to see this, let $w\in W^{1,2}(D, X)$ be such that $\trace(w) = \trace(v)$ almost everywhere. Then $w\circ\varphi^{-1}\in W^{1,2}(\Omega, X)$ and, moreover, $$\trace(w\circ\varphi^{-1}) = \trace(w)\circ\varphi^{-1}|_{\partial\Omega} = \trace(u|_{\Omega})$$ almost everywhere on $\partial \Omega$. The area minimizing property of $u$ and the discussion at the end of Section~\ref{sec:background-Sobolev-theory} then yield $$\Area_\mu(v) = \Area_\mu(u|_{\Omega}) \leq \Area_\mu(w\circ\varphi^{-1}) = \Area_\mu(w),$$ which proves \eqref{eq:area-min-v_n}. Lemma~\ref{lem:comparison-energy-volume} and \eqref{eq:energy-u-omega-boundary} imply $$\Area_\mu(v) \leq E^2(v)< Cl_0^2$$ and hence Proposition~\ref{prop:strong-cont?} implies that the restriction of $v$ to the closed ball $\bar{B}(0, \frac{1}{2})$ is $\alpha$-H\"older continuous with constant $L$. It follows that $$d(v(z), v(0))\leq L|z|^\alpha <  L\varrho^\alpha \leq \frac{\varepsilon}{2}$$ for every $|z|<\varrho$, establishing (ii). We conclude that $v$ satisfies the hypotheses of Lemma~\ref{lem:unif-cont-lower-bound-energy}. This completes the proof.
\end{proof}

As for the second main result of this section recall that a rectifiable Jordan curve $\Gamma \subset X$ is called a chord-arc curve if there exists $\lambda\geq 1$ such that  for any $x, y\in\Gamma$ the length of the shorter of the two segments in $\Gamma$ connecting $x$ and $y$ is bounded from above by $\lambda\cdot d(x,y)$.

\bt\label{thm:bdry-Hoelder}
  Let $\Gamma\subset X$ be a chord-arc curve and suppose $u\in \Lambda(\Gamma, X)$ is quasi-conformal and satisfies 
   \begin{equation}\label{eq:area-min-Jordan-bdry-Holder}
   \Area_\mu(u) = \inf\left\{\Area_\mu(v): v\in \Lambda(\Gamma, X)\right\}.
  \end{equation}
 Then the continuous representative of $u$ is H\"older continuous on all of $\overline{D}$.
 \et
 
 In fact, after possibly pre-composing with a Moebius transformation the continuous representative $\bar{u}$ of $u$ is $\beta$-H\"older on all of $\overline{D}$ with $$\beta= \frac{1}{4\pi Q^2 C(1+2\lambda)^2}.$$ 
Here, $Q$ is the quasi-conformality factor, $\lambda$ the parameter in the chord-arc condition for $\Gamma$, and $C$ is the isoperimetric constant. If $\Area_\mu(u)\leq Cl_0^2$ then it will furthermore follow that the H\"older constant $L$ of $\bar{u}$ is given by $L = M(\beta)\cdot\length_X(\Gamma)$ for some decreasing function $M(\beta)$. Note that Theorem~\ref{thm:bdry-Hoelder} implies statement (iii) of Theorem~\ref{thm:reg-intro}.

\begin{proof}
By Theorems~\ref{thm:int-reg-summary} and \ref{thm:bdry-cont-u-classical} we may assume that $u$ is continuous on all of $\overline{D}$.
Fix three points $p_1, p_2, p_3\in S^1$ at equal distance from each other and let $q_1,q_2,q_3\in\Gamma$ be three points such that the three segments into which they divide $\Gamma$ have equal length. After possibly pre-composing $u$ with a Moebius transformation we may assume that $u$ satisfies the $3$-point condition $u(p_i)=q_i$ for $i=1,2,3$.
Let $0<r_0\leq \frac{1}{2}$ be such that $\Area_\mu(u|_{D\cap B(z,r_0)})\leq Cl_0^2$ for every $z\in \overline{D}$. We claim that it is enough to prove that 
 \begin{equation}\label{eq:area-ball-length-bdry}
  \Area(u|_{D\cap B(z,r)}) \leq C(1+2\lambda)^2\cdot \length_X(u|_{D\cap \partial B(z,r)})^2
 \end{equation}
for every $z\in\overline{D}$ and almost every $r\in(0,r_0)$. 
 Indeed, if \eqref{eq:area-ball-length-bdry} is true then one argues as in the proof of Lemma 8.8 to obtain 
  \begin{equation*}
 \int_{B(z,r)} \mathcal{I}^1_+(\apmd u_w)\,d\lm^2(w)\leq \left[\pi E^2_+(u)\right]^{\frac{1}{2}}r_0^{-\beta} r^{1+\beta}
\end{equation*}
for every $z\in \overline{D}$ and every $r\in(0,r_0)$. Finally, the proof of Proposition 8.7 shows that for all $z_1,z_2\in \overline{D}$ one has $$d(u(z_1), u(z_2))\leq K\cdot\frac{r_0^{2\beta-1}}{1-2^{-\beta}}\cdot E_+^2(u)^{\frac{1}{2}}\cdot |z_1-z_2|^\beta$$ for some universal constant $K$. From this the statement of the theorem follows. Note that if $\Area_\mu(u)\leq Cl_0^2$ then $r_0$ can be taken to be $\frac{1}{2}$. In this case we have $E_+^2(u)\leq Q^2\Area(u)\leq Q^2C\length_X(\Gamma)^2$ and hence $$d(u(z_1), u(z_2))\leq M(\beta)\cdot \length_X(\Gamma)\cdot |z_1-z_2|^\beta$$ with $M(\beta) =K\left[\sqrt{\beta}(1-2^{-\beta})\right]^{-1}$ for some universal constant $K$, yielding the remark after the theorem.

It remains to show that \eqref{eq:area-ball-length-bdry} holds. For this, let $z\in \overline{D}$. Then $u|_{D\cap \partial B(z,r)}\in W^{1,2}(D\cap \partial B(z,r), X)$ for almost every $r\in(0,r_0)$. Fix such $r$. If $0<r<1-|z|$ then $$\Area(u|_{B(z,r)})\leq C\length_X(u|_{\partial B(z,r)})^2$$ by Lemma 8.6 and hence \eqref{eq:area-ball-length-bdry} in this case. If $r>1-|z|$ then denote by $a$ and $b$ the intersection points of $S^1$ with $\partial B(z,r)$. Since $u|_{S^1}$ satisfies the three-point condition and weakly monotonically parametrizes $\Gamma$ it follows that $$\length_X(u|_{S^1\cap B(z,r)})\leq \frac{2}{3}\cdot \length_X(\Gamma) \leq 2\cdot \length_X(u|_{S^1\setminus B(z,r)})$$ and hence the chord-arc property implies 
\begin{equation}\label{eq:estimate-bdry-curve-interior}
\length_X(u|_{S^1\cap B(z,r)}) \leq 2\lambda\cdot d(u(a),u(b)) \leq 2\lambda\cdot \length_X(u|_{D\cap \partial B(z,r)}).
\end{equation}
Set $\Omega:= D\cap B(z, r)$, and let $\beta\colon S^1\to\partial\Omega$ be an orientation preserving biLipschitz homeomorphism. By \eqref{eq:estimate-bdry-curve-interior} we have 
\begin{equation}\label{eq:estimate-u-bdry-omega}
\length_X(u \circ\beta) \leq (1+2\lambda)\cdot \length_X(u|_{D\cap\partial B(z,r)}).
\end{equation}
 Let $J\subset S^1$ be the segment that gets mapped by $\beta$ to $S^1\cap B(z, r)$. Define a homeomorphism $\psi\colon S^1\to S^1$ such that $\psi|_{S^1\setminus J}$ is the identity and such that on $J$ the map $\psi$ is a homeomorphism of $J$ as in Lemma~\ref{lem:reparam-Lip} for the curve $u\circ\beta$. It follows that $u\circ\beta\circ\psi\in W^{1,2}(S^1, X)$. 
Now, employing an argument similar to that in the proof of Lemma~\ref{lem:isop-u-area-min}, and using \eqref{eq:estimate-u-bdry-omega} as well as the fact that $u$ satisfies \eqref{eq:area-min-Jordan-bdry-Holder}, one shows that \eqref{eq:area-ball-length-bdry} holds. Indeed, by Lemma~\ref{lem:isop-Sobolev-curves}, there exists $v\in W^{1,2}(D, X)$ such that $$\Area_\mu(v)\leq C \length_X(u\circ\beta)^2$$ and such that $\trace(v) = u\circ\beta\circ\psi$ almost everywhere. Let $\varrho\colon S^1\to S^1$ be the homeomorphism which is given by $\beta\circ\psi\circ\beta^{-1}$ on $\beta(J)$ and which is the identity on $S^1\setminus \beta(J)$. Let $\varphi\colon \overline{D}\to \overline{\Omega}$ be a biLip\-schitz map extending $\beta$. Such $\varphi$ exists by \cite[Theorem A]{Tuk80}.   Let $\bar{u}\colon D\to X$ be the map which agrees with $u$ on $D\setminus\Omega$ and with $v\circ\varphi^{-1}$ on $\Omega$. A similar argument as in the proof of Lemma~\ref{lem:gluing-Sobolev} shows that $\bar{u}\in W^{1,2}(D, X)$ and $\trace(\bar{u}) = u\circ\varrho$ almost everywhere. Since $u$ satisfies \eqref{eq:area-min-Jordan-bdry-Holder} it follows that
\begin{equation*}
\Area_\mu(u) \leq \Area_\mu(\bar{u}) \leq \Area_\mu(u|_{D\setminus\Omega}) + C\length_X(u\circ\beta)^2
\end{equation*}
and hence $$\Area_\mu(u|_\Omega)\leq C\length_X(u\circ\beta)^2\leq C(1+2\lambda)^2\cdot \length_X(u|_{D\cap\partial B(z,r)})^2.$$ This proves \eqref{eq:area-ball-length-bdry} and completes the proof of the theorem. 
\end{proof}

\section{Proofs of Corollaries~\ref{cor:Abs-Plateau} and \ref{cor:tree}}\label{sec:abs-Plateau-corollaries}

Throughout this section, let $\mu$ be a definition of volume which induces quasi-convex $2$-volume densities.

Corollary~\ref{cor:tree} is a special case of the following result.

\bt\label{thm:tree-gen}
Let $X$ be a proper, geodesic metric space admitting a global quadratic isoperimetric inequality with some constant $C$ for $\mu$. If $C<\frac{1}{8\pi}$ then $X$ is a metric tree, that is, every geodesic triangle in $X$ is isometric to a tripod.
\et
 
 \begin{proof}
 lt suffices to show that $X$ does not contain any rectifiable Jordan curve. Suppose to the contrary that there exists a rectifiable Jordan curve $\Gamma$ in $X$. The global quadratic isoperimetric inequality implies that $\Lambda(\Gamma, X)\not=\emptyset$. By Theorem~\ref{thm:existence-qc-area-min-general-mu} there exists $u\in \Lambda(\Gamma, X)$ which is $\sqrt{2}$-quasi-conformal and minimizes the $\mu$-area among all maps in $\Lambda(\Gamma, X)$. By Theorem~\ref{thm:int-reg-summary}, we may assume $u$ to be locally $\alpha$-H\"older continuous with $\alpha = \frac{1}{8\pi C}$. Since $\alpha>1$ it follows that $u$ is constant on $\overline{D}$, thus contradicting the fact that $\trace(u)$ is a weakly monotone parametrization of the Jordan curve $\Gamma$. If $X$ satisfies property (ET) then $u$ may be chosen to be conformal and hence locally $\alpha$-H\"older continuous on $D$ with $\alpha = \frac{1}{4\pi C}$. Therefore, if $C<\frac{1}{4\pi}$ then $\alpha>1$ and it follows that $u$ is constant.
\end{proof}

Recall that a metric space $X$ is injective if it is an absolute $1$-Lipschitz retract. Equivalently, $X$ is injective if for every metric space $Y$, every subset $A\subset Y$, and every Lipschitz map from $A$ to $X$ there exists a Lipschitz extension to all of $Y$ with the same Lipschitz constant. Examples of such spaces include metric trees, $\ell^\infty(W)$ for every set $W$, and $L^\infty(Z, \mu)$ for any measure space $(Z, \mu)$. Every injective space is complete and geodesic, see e.g.~\cite{Lan13}.

\bt\label{thm:Plateau-injective-metric}
 Let $\Gamma\subset X$ be a rectifiable Jordan curve. Then there exists $u\in \Lambda(\Gamma, X)$ such that 
    \begin{equation*}
   \Area_\mu(u) = \inf\left\{\Area_\mu(v): v\in \Lambda(\Gamma, X)\right\}
  \end{equation*}
 and such that $u$ is $\sqrt{2}$-quasi-conformal. Moreover, $u\in W^{1,p}_{\rm loc}(D, X)$ for some $p>2$, and $u$ has a representative which is continuous on $\overline{D}$ and locally $\frac{1}{4}$-H\"older continuous on $D$.
\et

No assumption on local compactness of $X$ is needed. It applies, in particular, to all $L^\infty$-spaces.

We begin with the following lemma.

\bl\label{lem:isop-injective-space}
 Let $X$ be an injective metric space. Then $X$ admits a global isoperimetric inequality with constant $\frac{1}{2\pi}$ for $\mu$.
\el

\begin{proof}
 Let $c\colon S^1\to X$ be a Lipschitz curve such that $r:=\length_X(c)>0$. Let $\bar{c}\colon S^1\to X$ be the constant speed parametrization of $c$. Endow $S^1$ with the length metric. Then $c$ is $(2\pi)^{-1}r$-Lipschitz and, since $X$ is injective,  it has a $(2\pi)^{-1}r$-Lipschitz extension $\varphi$ to the standard upper hemisphere $S^2_+$. Pre-composing $\varphi$ with a bijective Lipschitz map from $\overline{D}$ to $S^2_+$ which restricts to the identity on $S^1$ we obtain a Lipschitz extension $w\colon \overline{D}\to X$ of $\overline{c}$ whose $\mu$-area is bounded above by $$\Area_\mu(w) \leq \lip(\varphi)^2\cdot\Area(S^2_+)\leq (2\pi)^{-1}r^2.$$ Now, one constructs exactly as in Lemma~\ref{lem:isop-Sobolev-curves} a map $v\colon \bar{B}(0,2)\to X$ which coincides with $w$ on $D$ and which gives a `linear' reparametrization from $c$ to $\bar{c}$ on the annulus $\bar{B}(0,2)\setminus D$. Since $c$ is a Lipschitz map it follows from the construction that $v$ is Lipschitz; moreover, $\Area_\mu(v|_{A})=0$. Identifying $\bar{B}(0,2)$ with $\overline{D}$ via the scaling map we thus obtain a Lipschitz extension of $c$ to $\overline{D}$ whose $\mu$-area is bounded by $(2\pi)^{-1}r^2$. This completes the proof.
\end{proof}

\begin{proof}[{Proof of Theorem~\ref{thm:Plateau-injective-metric}}]
 By \cite{Isb64}, there exists an injective hull $Y$ of $\Gamma$, which is moreover compact and isometrically embeds into $X$, see also \cite{Lan13}. Since $Y$ is injective there exists a $1$-Lipschitz retraction $r\colon X\to Y$. By Lemma~\ref{lem:isop-injective-space}, $Y$ admits a global isoperimetric inequality with constant $C=\frac{1}{2\pi}$ for $\mu$. In particular, we have that $\Lambda(\Gamma, Y)\not=\emptyset$.
 By Theorem~\ref{thm:existence-qc-area-min-general-mu}, there exists $u\in \Lambda(\Gamma, Y)$ which is $\sqrt{2}$-quasi-conformal and minimizes the $\mu$-area among all maps in $\Lambda(\Gamma, Y)$. Since $r$ is a $1$-Lipschitz retraction it follows that $u$ also minimizes the $\mu$-area among all maps in $\Lambda(\Gamma, X)$.
 Theorem~\ref{thm:int-reg-summary} shows that $u\in W^{1,p}_{\rm loc}(D, X)$ for some $p>2$ and that $u$ has a representative which is locally $\alpha$-H\"older with $\alpha= \frac{1}{8\pi C} = \frac{1}{4}$. By Theorem~\ref{thm:bdry-cont-u-classical}, the continuous representative of $u$ extends continuously to $\overline{D}$. This completes the proof.
\end{proof}

Given a metric space $\Gamma$ homeomorphic to $S^1$ and of finite length, define
 $$m(\Gamma, \mu):= \inf\{\Area_\mu(v): \text{$Y$ complete, $\iota\colon\Gamma\hookrightarrow Y$ isometric, $u\in\Lambda(\iota(\Gamma), Y)$}\}.$$ 
Corollary~\ref{cor:Abs-Plateau} is a special case of the following result.

\bc\label{cor:Abs-Plateau-general}
There exist  a compact metric space $X$, an isometric embedding $\iota\colon \Gamma\hookrightarrow X$, and a map $u\in \Lambda(\iota(\Gamma), X)$ such that $$\Area_\mu(u) = m(\Gamma, \mu).$$ Moreover, $u$ is $\sqrt{2}$-quasi-conformal and has a representative which is continuous on $\overline{D}$ and locally $\frac{1}{4}$-H\"older continuous on $D$.
\ec

\begin{proof}
 Let $X$ be an injective hull of $\Gamma$, see \cite{Isb64}. By Theorem~\ref{thm:Plateau-injective-metric} there exists $u\in\Lambda(\Gamma, X)$ with minimal area among maps in $\Lambda(\Gamma, X)$ and which satisfies the regularity properties required in Corollary~\ref{cor:Abs-Plateau-general}. Finally, since $X$ is an injective metric space and since the area does not increase under compositions with $1$-Lipschitz maps,  we have $$m(\Gamma, \mu)= \inf\{\Area_\mu(v) : v\in \Lambda(\Gamma, X)\}$$ and hence $\Area_\mu(u)=m(\Gamma, \mu)$. This completes the proof.
\end{proof}

\section{The infinitesimally Euclidean case}\label{sec:ET-case}

In what follows let $(X,d)$ be a complete metric space.

\bd
 We say that $X$ has property (ET) if for every $u\in W^{1,2}(D, X)$ the approximate metric derivative $\apmd u_z$ is induced by a possibly degenerate inner product at almost every $z\in D$.
\ed

Many geometrically interesting classes of spaces have property (ET).  For instance, this is the case for Riemannian manifolds with continuous metric tensor, metric spaces of curvature bounded from above or below in the sense of Alexandrov, and equiregular sub-Riemannian manifolds. In order to see this, we only need to observe that  in every such space no metric blow-up (tangent cone) at any point  may contain non-Euclidean normed metric spaces. Then the result follows from
the proposition below. We refer to \cite{BrH99} for basics on ultralimits and to \cite{Lyt04} for more about blow-ups and tangent cones.

\bp\label{prop:suff-cond-ET}
 Let $\omega$ be a non-principal ultrafilter on $\N$. Suppose that for every $x\in X$, there is some sequence $r_j\to\infty$ such that  the ultralimit $X_\omega$ of the sequence $(X, r_jd, x)$ does not contain  isometrically embedded $2$-dimensional non-Euclidean normed spaces. Then $X$ has property (ET).
\ep

\begin{proof}
 This follows directly from Proposition~\ref{prop:Sobolev-apmd-Lip}.
\end{proof}

Another interesting class of spaces with property (ET) is given by infinitesimally Hilbertian metric spaces with (synthetic) Ricci curvature bounded below. More precisely, if $(X,d, \mathfrak{m})$ is an infinitesimally Hilbertian ${\rm CD}^*(K,N)$ space for some $K\in\R$ and $N\in[1,\infty)$, see e.g.~\cite{GMR} for the terminology, then $(X,d)$ has property (ET). Indeed, for each $x\in X$ the collection of (measured) tangents of $(X, d, \mathfrak{m})$ is non-empty and each tangent $(Y, d_Y, \mathfrak{n})$ is an infinitesimally Hilbertian ${\rm CD}^*(0,N)$ space, see (2.7) of \cite{GMR}. In particular, by \cite{BS10}, the support of the measure $\mathfrak{n}$ is all of $Y$. Thus, if $Y$ contains a normed plane $V$ then the Splitting Theorem \cite{Gig} implies that $V$ must be Euclidean. From this and Proposition~\ref{prop:suff-cond-ET} it follows that $X$ has property (ET).

The validity of property (ET) simplifies many results and formulas. 

\bt\label{thm:inner-var-min-conformal-XEucl}
  Let $X$ satisfy property (ET) and let $u\in W^{1,2}(D,X)$. If $$E_+^2(u)\leq E_+^2(u\circ\psi)$$ for every biLipschitz homeomorphism $\psi\colon D\to D$ then $u$ is conformal. The same statement holds when $E_+^2$ is replaced by $E^2$.
\et

\begin{proof}
This follows from the same arguments as in the proof of Theorem~\ref{thm:qc-domain-minimizers}. Indeed, since $\apmd u_z$ comes from an inner product for almost every $z\in D$ for which $\apmd u_z$ is non-degenerate, it follows from \eqref{eq:E-ineq-for-qc} and Lemma~\ref{lem:qc-seminorm} that $\apmd u_z$ is conformal for almost every $z\in D$. Hence, $u$ is conformal.  Using Lemma~\ref{lem:qc-seminorm-KS-energy} instead of Lemma~\ref{lem:qc-seminorm} one obtains the second statement.
\end{proof}

If $X$ satisfies property (ET) then for every $u\in W^{1,2}(D, X)$ and any two definitions of volume $\mu_1$ and $\mu_2$ one has $\Area_{\mu_1}(u) = \Area_{\mu_2}(u)$ by property (i) of Definition~\ref{def:volume-def}. We will therefore simply write $\Area(u)$ in this case.

\bt\label{thm:energy-min-is-area-min-ET}
 Let $X$ satisfy property (ET) and let $\Gamma\subset X$ be a Jordan curve. If $u\in \Lambda(\Gamma, X)$ satisfies 
 \begin{equation*}
  E_+^2(u) = \inf\left\{E_+^2(u'): u'\in\Lambda(\Gamma, X)\right\}
 \end{equation*}
 then $u$ is conformal and an area minimizer, that is, 
 \begin{equation*}
  \Area(u) = \inf\left\{\Area(u'): u'\in \Lambda(\Gamma, X)\right\}.
 \end{equation*}
 The same statement holds when $E_+^2$ is replaced by $E^2$.
\et

\begin{proof}
 The fact that $u$ is conformal is a direct consequence of Theorem~\ref{thm:inner-var-min-conformal-XEucl}. We show that $u$ is an area minimizer. Arguing by contradiction we assume there exists $v\in\Lambda(\Gamma, X)$ such that $$\Area(v)< \Area(u).$$
 Arguing exactly as in the first part of the proof of Theorem~\ref{thm:existence-qc-area-min-general-mu} but using Theorem~\ref{thm:inner-var-min-conformal-XEucl} instead of Theorem~\ref{thm:qc-domain-minimizers}, one shows that there exists $w\in \Lambda(\Gamma, X)$ which is conformal and satisfies $$\Area(w)\leq \Area(v).$$ Together with Lemma~\ref{lem:comparison-energy-volume} one thus obtains that $$E_+^2(w) = \Area(w) <\Area(u) = E_+^2(u)$$ which contradicts the fact that $u$ minimizes $E_+^2$. It follows that $u$ is an area minimizer. The proof for $E^2$ is analogous.
\end{proof}

Combining Theorem~\ref{thm:energy-min-is-area-min-ET} with Theorem~\ref{thm:existence-energy-min} we obtain the existence of conformal area minimizers.

\bc\label{cor:ex-conf-area-min-ET}
Let $X$ be a proper metric space satisfying property (ET) and let $\Gamma\subset X$ be a Jordan curve. If $\Lambda(\Gamma, X)\not=\emptyset$ then there exists $u\in \Lambda(\Gamma, X)$ which minimizes the $E^2_+$-energy among all maps in $\Lambda(\Gamma, X)$. Every such $u$ is conformal and minimizes the area among all maps in $\Lambda(\Gamma, X)$. 
\ec

The same holds with $E^2_+$ replaced by the Korevaar-Schoen energy $E^2$.

We now show that in spaces without property (ET) area minimizers with respect to two different definitions of area are in general different.

\bp\label{prop:area-min-diff}
 Let $\mu$ and $\bar{\mu}$ be quasi-convex definitions of volume such that $\mu_V \not=\bar{\mu}_V$ for some normed plane $V$. Then there exist a metric space $X$ biLipschitz homeomorphic to $S^2$ and a closed biLipschitz curve $\Gamma$ in $X$ such that for every $u\in\Lambda(\Gamma, X)$ with $$\Area_\mu(u) =  \inf\left\{\Area_\mu(v): v\in \Lambda(\Gamma, X)\right\}$$ there is some $v\in\Lambda(\Gamma, X)$ with $\Area_{\bar{\mu}}(v)<\Area_{\bar{\mu}}(u)$.
\ep

It follows, in particular, that energy minimizers with respect to a fixed definition of energy (for example the Reshetnyak or Korevaar-Schoen energy) can in general only be area minimizers with respect to at most one definition of area. In \cite{LW-energy-area} we show that Reshetnyak energy minimizers are in fact area minimizers with respect to the intrinsic Riemannian volume $\mu^{\rm i}$ and that, more generally, for every suitable notion of quasi-convex energy $\tilde{E}$ there is an induced quasi-convex definition of area $\tilde{\mu}$ such that $\tilde{E}$-energy minimizers are $\tilde{\mu}$-area minimizers.

\begin{proof}
Let $\|\cdot\|$ be a norm on $\R^2$ such that $\mu_V \not=\bar{\mu}_V$ for $V=(\R^2, \|\cdot\|)$. We may assume that $\mu_V<\bar{\mu}_V$, the proof for the other case being analogous. Let $\lambda>0$ be such that $$\mu_V(D) < \lambda^2\pi <\bar{\mu}_V(D),$$ where $D$ denotes the Euclidean unit disc as usual. Let $\overline{D}_1$ and $\overline{D}_2$ be two copies of $\overline{D}$. Endow $\overline{D}_1$ with the metric coming from the norm $\|\cdot\|$ and $\overline{D}_2$ with $\lambda$ times the Euclidean metric. 
Let $X$ be the metric space obtained by gluing $\overline{D}_1$ and $\overline{D}_2$ along their boundaries, endowed with the quotient metric. Then $X$ is biLipschitz homeomorphic to the standard sphere $S^2$ and, in particular, admits a $(C,l_0)$-isoperimetric inequality for some $C, l_0>0$ for every definition of volume. Embed $\overline{D}_i$ into $X$ via the natural inclusion and denote by $\Gamma\subset X$ the boundary of $\overline{D}_i$. Then $\Gamma$ is a closed biLipschitz curve. For $j=1$, $2$, let $u_j\colon \overline{D}\to \overline{D}_j \hookrightarrow X$  be the natural inclusion. Then $u_j\in\Lambda(\Gamma, X)$ and $\Area_{\bar{\mu}}(u_2)=\lambda^2\pi$ and $\Area_{\mu}(u_1)=\mu_V(D)$.

Let $u\in\Lambda(\Gamma, X)$ be such that $$\Area_\mu(u) = \inf\left\{\Area_\mu(u'): u'\in\Lambda(\Gamma, X)\right\}.$$ Since $\Lambda(\Gamma, X)$ is not empty such $u$ exists by Theorem~\ref{thm:existence-qc-area-min-general-mu}. We claim that 
\begin{equation}\label{eq:area-too-big}
\Area_{\bar{\mu}}(u) \geq \bar{\mu}_V(D)
\end{equation}
and thus $\Area_{\bar{\mu}}(u) \geq \bar{\mu}_V(D)>\lambda^2\pi = \Area_{\bar{\mu}}(u_2)$, which shows that $u$ is not an area minimizer in $\Lambda(\Gamma, X)$ for $\bar{\mu}$.

It remains to prove \eqref{eq:area-too-big}. Due to the quasi-convexity of $\mu$, $\bar{\mu}$, and $E_+^2$ and Proposition~\ref{prop:seq-Jordan-equi-bdd-energy} we find a map $\hat{u}\in \Lambda (\Gamma ,X)$  which has minimal $E_+^2$-energy among all maps $v\in \Lambda (\Gamma ,X)$ satisfying
$\Area _{\mu} (v) = \Area _{\mu} (u)$ and $\Area _{\bar{\mu}} (v) \leq \Area _{\bar{\mu}} (u)$.  By Theorem~\ref{thm:qc-domain-minimizers}, such a map $\hat{u}$ is quasi-conformal.

By Propositions~\ref{prop:hoelder-cont-min} and \ref{thm:bdry-cont-u-classical}, we may assume that $\hat{u}$ is continuous on $\overline{D}$ and satisfies Lusin's property (N). In particular, there exists $i\in\{1,2\}$ such that $\overline{D}_i\subset \hat{u}(\overline{D})$. Since
 \begin{equation*}
  \mu_V(D)= \Area_\mu(u_1) \geq \Area_{\mu}(\hat{u})
 \end{equation*}
 the area formula thus implies that $i=1$ and hence that $\Area_{\bar{\mu}}(\hat{u}) \geq \bar{\mu}_V(D)$. Since $\Area_{\bar{\mu}}(u)\geq \Area_{\bar{\mu}}(\hat{u})$ this proves \eqref{eq:area-too-big} and completes the proof.
\end{proof}

We conclude the paper by noting that property (ET) implies a corresponding property in all dimensions.

\bp
Let $X$ satisfy property (ET), let $n\in\N$, and let $\Omega\subset\R^n$ be an open, bounded subset. Then for any $u\in W^{1,2}(\Omega, X)$ and almost every point $z\in\Omega$ the approximate metric derivative $\apmd u_z$ is induced by a possibly degenerate inner product.
\ep

\begin{proof}
We may assume that $\Omega$ is a ball since the claim is local. We now argue by contradiction and assume that $u\in W^{1,2}(\Omega, X)$ is such that, on a set of strictly positive measure, $\apmd u_z$ is not induced by a possibly degenerate inner product. Fix a countable dense set of $2$-planes $V_i$ in $\R^n$. Slicing $\Omega$ by translates of the $V_i$, we find a point $z\in \Omega$ at which $\apmd u_z$ exists and is non-degenerate, is not induced by an inner product, but is such that the restriction 
of $\apmd u_z$ to each $V_i$ is given by an inner product. However, by the parallelogram identity, a norm comes from an inner product if and only if its restriction to each $2$-plane comes from an inner product. By the density of the planes $V_i$ this leads to a contradiction.
\end{proof}

\def\cprime{$'$} \def\cprime{$'$}
\providecommand{\bysame}{\leavevmode\hbox to3em{\hrulefill}\thinspace}
\providecommand{\MR}{\relax\ifhmode\unskip\space\fi MR }
\providecommand{\MRhref}[2]{%
  \href{http://www.ams.org/mathscinet-getitem?mr=#1}{#2}
}
\providecommand{\href}[2]{#2}

\end{document}